\documentclass[11pt]{article}
\usepackage[english]{babel} 
\usepackage[utf8]{inputenc} 
\usepackage[a4paper,left=2.5cm,bottom=2.5cm,top=2cm]{geometry} 
\usepackage{amsmath} 
\usepackage{amssymb}
\usepackage{amsfonts}
\usepackage{bm} 
\usepackage{graphicx} 
\usepackage[export]{adjustbox} 
\usepackage{pdflscape} 
\usepackage{fancyhdr} 
\usepackage{color} 
\usepackage{wrapfig} 
\usepackage{lipsum} 
\usepackage[normalem]{ulem} 
\usepackage{blindtext}
\usepackage{titlesec}
\usepackage{algpseudocode}
\usepackage{enumitem}
\usepackage{caption}
\usepackage{subcaption}
\usepackage{tcolorbox}
\usepackage{glossaries}
\usepackage{csvsimple}
\usepackage{multicol}
\usepackage{multirow}
\usepackage{adjustbox}
\usepackage{booktabs}
\usepackage[bookmarks=false]{hyperref}
\usepackage{natbib}
\usepackage{colortbl}
\usepackage[toc,page]{appendix}
\usepackage{float}
\usepackage{url}

\newcolumntype{C}[1]{>{\centering\let\newline\\\arraybackslash\hspace{0pt}}m{#1}}

\floatstyle{ruled}
\newfloat{model}{thp}{lop}
\floatname{model}{Eq.}

\usepackage{authblk}
\usepackage{rotating}
\usepackage{makecell} 

\RequirePackage{etoolbox}
\allowdisplaybreaks

\newcommand{\buses}{\mathcal{N}}
\newcommand{\tracts}{\mathcal{C}}

\definecolor{ForestGreen}{RGB}{34,139,34}

\newacronym{mip}{MIP}{mixed-integer program}
\newacronym[longplural = public safety power shutoffs]{psps}{PSPS}{public safety power shutoff}
\newacronym{ercot}{ERCOT}{Electric Reliability Council of Texas}
\newacronym{wfpi}{WFPI}{Wind-Enhanced Fire Potential Index}
\newacronym{usgs}{USGS}{United States Geological Survey}
\newacronym{pge}{PG\&E}{Pacific Gas \& Electric Company}
\newacronym{svi}{SVI}{Social Vulnerability Index}
\newacronym{eji}{EJI}{Environmental Justice Index}
\newacronym{j40}{J40}{Justice40}
\newacronym{cdc}{CDC}{Center for Disease Control and Prevention}
\newacronym{cejst}{CEJST}{Climate \& Economic Justice Screening Tool}
\newacronym{atsdr}{ATSDR}{Agency for Toxic Substances and Disease Registry}
\newacronym{mmf}{MMF}{Min-Max Fairness}
\definecolor{red}{rgb}{0,0,0}
\renewcommand{\sout}[1]{\unskip}

\title{Equitably Allocating Wildfire Resilience Investments for Power Grids -- The Curse of Aggregation and Vulnerability Indices}

\author[a]{Madeleine Pollack}
\author[b]{Ryan Piansky}
\author[a]{Swati Gupta}
\author[b]{Daniel Molzahn}

\affil[a]{Massachusetts Institute of Technology\protect \\ {\small \tt \{pollack9,swatig\}@mit.edu}}
\affil[b]{Georgia Institute of Technology\protect \\ {\small \tt \{rpiansky3, molzahn\}@gatech.edu}}

\date{\today}

\begin{document}

\maketitle

\begin{abstract}
    Social vulnerability indices have increased traction for guiding infrastructure investment decisions to prioritize communities that need these investments most. One such plan is the Biden-Harris Justice40 initiative, which aims to guide equitable infrastructure investments by ensuring that disadvantaged communities defined by the \gls{cejst} receive 40\% of the total benefit realized by the investment. However, there is limited research on the practicality of applying vulnerability indices like the \gls{cejst} to real-world decision-making for policy outcomes. In this paper, we study this gap by examining the effectiveness of vulnerability indices in a case study focused on power shutoff and undergrounding decisions in wildfire-prone regions. Using a mixed-integer program and a high-fidelity synthetic transmission network in Texas, we model resource allocation policies inspired by Justice40 and evaluate their impact on reducing power outages and mitigating wildfire risk for vulnerable groups. 
    Our analysis reveals that the Justice40 framework may fail to protect certain communities facing high wildfire risk. In our case study, we show that Indigenous groups are particularly impacted. We posit that this outcome is likely due to information losses from data aggregation and the use of generalized vulnerability indices.
    Through the use of explicit group-level protections, we provide bounds on the best possible outcome for population groups that are proportionally most affected.
\end{abstract}

\pagenumbering{arabic}
\glsresetall

\section{Introduction}\label{sec: introduction}

The increased frequency and severity of extreme climate events in the 21st century necessitates policy and infrastructure change to mitigate property damage, environmental impacts, and loss of life from these extreme events.
While climate change has widespread negative effects, the degree of impact is non-uniform across communities. In particular, research shows that historically marginalized and low-income communities are most susceptible to environmental hazards such as water contamination \citep{campbell2016case,schaider2019environmental}, excessive heat \citep{gronlund2014racial,jesdale2013racial}, and air pollution \citep{ard2016two,woo2019residential}. Hence, it is important that legislators prioritize the protection of these groups when allocating public funds for climate infrastructure.

With these goals in mind, the Biden-Harris administration issued United States Executive Order 14008 \citep{exec_order2021}, which established the Justice40 initiative in its overarching plan of ``tackling the climate crisis." The Justice40 initiative has the stated goal that 40\% of the overall benefits of federal investments in environmental and energy infrastructure should flow to disadvantaged communities. In November 2022, the US Council on Environmental Quality launched the \gls{cejst} to serve as a \textit{vulnerability index} for identifying census tracts that should be prioritized for investment \citep{ScreeningTool_2022}.

The \gls{cejst} is the result of years of dedicated research on climate and socioeconomic vulnerability; however, there is currently very little research analyzing the practical effectiveness of specific policies deployed in line with the Justice40 initiative's goals. This is not unique to the \gls{cejst} or the Justice40 initiative; various vulnerability indices have been scrutinized for lack of demonstrated efficacy in their intended contexts. 
\citet{rufat2019valid} evaluates the CDC's \gls{svi} \citep{CDC_SVI_data,flanagan2011social} and Cutter's Social Vulnerability Index (SoVI) \citep{cutter2003social} in their ability to identify which groups and communities were most detrimentally impacted by Hurricane Sandy in 2012, and concludes that these metrics lack ``construct validity"\footnote{Measures that have construct validity relate to other measures 
 that are known to quantify a phenomenon \citep{strauss2009construct}.}. \citet{rufat2019valid} and \citet{spielman2020evaluating} each advise expert development of indices for the specific case study as opposed to blindly applying the SoVI or SVI. However, all of these studies recognize the importance of measuring, monitoring, and prioritizing social vulnerability when planning for hazardous events 
 \citep{spielman2020evaluating,rufat2019valid}, and they call for better development of benchmarks for developing and evaluating appropriate indices~\citep{fekete2019social}. 
These studies underscore the need for detailed analyses that examine how these indices influence policy decisions and resource distribution and whether index-based policies truly benefit the targeted populations. 

Building on this critical gap, our paper examines the outcomes of applying the Justice40 initiative to budget allocation policies relevant to infrastructure investment for wildfire resilience. Specifically, we focus on optimizing budget allocations for undergrounding transmission lines in high wildfire risk areas that currently rely on emergency power shutoffs to mitigate both wildfire risk and power outages. Our analysis centers on evaluating the effectiveness of two vulnerability indices—the CDC's \gls{svi} and the \gls{cejst}—in directing funds to vulnerable groups, which we define as communities more susceptible to wildfire risk due to social, political, economic, and institutional disenfranchisement~\citep{fekete2019social}.

To assess the impact of the Justice40 allocation policy, we (1) examine the ability of the \gls{cejst} and \gls{svi} to identify high-risk, high-need populations for undergrounding investments, and (2) compute the \textit{load shed}, or quantity of power loss, experienced by different populations on the network after undergrounding decisions have been implemented under different interpretations of the Justice40 initiative. Additionally, we compare load shedding outcomes from these allocation policies with those output from a \gls{mmf} framework, a widely used approach in fair operations research literature. Although this analysis is applied to a case study in Texas, the proposed methodology is generalizable to transmission networks in other regions.

In this paper, we demonstrate that allocating budget according to Justice40-style constraints based on vulnerability indices fails to reduce power outages for Indigenous communities in the representative region managed by the \gls{ercot} and for areas at high risk of wildfire ignition. This finding is based on analysis using a high-fidelity synthetic power grid dataset that replicates the key characteristics of \gls{ercot}. We explore how aggregating community features into a single vulnerability index and employing a ``one-size-fits-all'' approach to defining vulnerability may contribute to a misalignment between the intended goals of these indices and their actual impact in this context. 
Then, using an \gls{mmf} framework to include group-level protections on the percentage of load shed across each population group of interest, we provide a lower bound on the percentage load shed experienced by the proportionally most-affected population group.
Our work thus adds to the ongoing discourse on the construct validity of vulnerability indices, offering empirical insights into their performance in a specific application. This approach aids in identifying the strengths and limitations of these tools and helps guide future refinements to enhance their effectiveness and relevance in policy-making.

The rest of the paper is organized as follows. First, \S\ref{sec: literature review} provides a brief literature review and background on wildfire risk management, social vulnerability indices, and our Texas case study. Then, in \S\ref{sec:methods}, we discuss our optimization framework for making line de-energization and undergrounding decisions under various objectives and constraints for equitable budget allocation. In \S\ref{sec:case_study}, we show the network configuration and demographic profile of our Texas case study. Finally, in \S\ref{sec:results}, we show the results of our optimization models, and in \S\ref{sec: policy insights}, we discuss the limitations of vulnerability indices and the possible causes.

\section{Literature Review}\label{sec: literature review}
The relevant literature for this paper can be divided into three categories: (1) information on the known or perceived efficacy of social vulnerability indices, (2) sources pertaining to the context of our case study on wildfire risk and its management, and (3) background on Texas' environmental risk and demographic profile, validating its use for our case study.
\subsection{The Role and Application of Vulnerability Indices in Policy-Making}

Vulnerability indices have emerged as a critical tool for identifying and prioritizing communities that are disproportionately impacted by environmental hazards and socioeconomic inequalities. In this work, we consider both the \gls{cejst} and \gls{svi} in our analyses.

The \gls{svi}, developed by the \gls{cdc} and \gls{atsdr}, assesses social vulnerability based on factors such as socioeconomic status, household composition, minority status, housing, and transportation \citep{CDC_SVI_data}. This index has been widely applied in various domains, including disaster management \citep{flanagan2011social}, healthcare resource allocation \citep{ganatra2022impact, wolkin2022comparison}, and environmental justice assessments \citep{ramesh2022flooding}.

Similarly, the \gls{cejst} was created by the US Council on Environmental Quality via Executive Order 14008 to identify disadvantaged census tracts for use by the Justice40 initiative. This initiative aims to direct 40\% of the benefits of federal infrastructure investments to these disadvantaged communities \citep{exec_order2021, White_House_2023}. The \gls{cejst} incorporates a range of indicators, including environmental, climate, health, and socioeconomic factors, to identify communities eligible for Justice40 benefits. 

There is much optimism regarding the ability of vulnerability indices to direct investments toward at-risk groups. In the context of energy systems, recent studies have begun to integrate social vulnerability assessments into infrastructure planning and emergency response strategies. \citet{Taylor2023Managing} explore the use of the \gls{svi} in configuring microgrids to reduce load shedding impacts on vulnerable communities during wildfire-related outages. \citet{ganz2023socioeconomic} examine the socioeconomic impacts of renewable energy deployments, highlighting the importance of incorporating vulnerability considerations in energy transitions. Other research applications making use of vulnerability indices include healthcare~\citep{ganatra2022impact,wolkin2022comparison}, and disaster relief~\citep{flanagan2011social,ramesh2022flooding}.

However, empirical evaluations assessing how effectively these indices inform specific infrastructure investment decisions, such as those related to wildfire risk mitigation, remain limited. \citet{fekete2019social} and \citet{rufat2019valid} emphasize the need for more comprehensive research to validate these tools and ensure that they accurately represent vulnerability and effectively guide resource allocation. Concerns have been raised regarding the potential oversimplification of complex social dynamics and the risk of misidentifying or overlooking certain vulnerable populations due to methodological limitations. This study seeks to help address this gap by analyzing the application and outcomes of vulnerability index-informed resource allocation strategies in a power grid wildfire resilience context.

\subsection{Mitigating Risk of Wildfire from Power Lines}
While power line-sparked ignitions account for only a small fraction of total wildfires in the United States, these ignitions often occur under extreme environmental conditions that lead to catastrophic spread and destruction \citep{mitchell2013power}. Notable incidents, such as the Smokehouse Creek (Texas), Maui Morning Fire (Hawaii), Echo Mountain (Oregon), and Camp (California) fires, have been linked to power infrastructure ignitions \citep{Penn_2024}.

To address this issue, high wildfire-risk states employ emergency \gls{psps} events, where power lines are proactively de-energized during high-risk conditions to prevent ignitions \citep{huang2023review}. While effective in reducing wildfire occurrences, \gls{psps} events can result in power outages, referred to as \emph{load shed}, adversely affecting communities, especially those with higher vulnerability due to factors such as food insecurity, poverty, and medical needs \citep{Fuller_2019, Ham_Lee_2022}.

While \gls{psps} events provide short-term solutions to high wildfire risk, a promising long-term solution is to ``underground" (i.e., bury) high-risk power lines, which significantly reduces the likelihood of wildfire ignitions. However, the high costs associated with undergrounding—ranging from \$5 to \$10 million per mile—and the extended timelines required for such projects pose substantial challenges \citep{hall2012out, yang2022resilient, PGE, Beam_2023}. Consequently, determining how to allocate limited budgets for undergrounding projects becomes a complex resource allocation problem that must balance efficiency and equity considerations.

Recent research efforts have focused on optimizing PSPS implementations to balance wildfire risk reduction and minimizing adverse impacts on communities. \citet{rhodes2020balancing} and \citet{huang2023review} utilize mixed-integer programming models to determine optimal combinations of power line shutoffs that minimize both fire risk and the extent of load shedding. Other works have optimized \gls{psps} operations to reduce disruptions, ensure reliability, and incorporate renewable resources \citep{rhodes2022co, kody2022load, rhodes2024security, kody2022optimizing}.

Advancements in modeling techniques have incorporated stochastic elements to account for uncertainties in wildfire conditions. \citet{su2024quasi} develop a stochastic mixed-integer nonlinear programming approach that considers variable environmental factors to optimize PSPS plans. Additionally, the integration of renewable energy storage solutions, such as long-duration batteries, has been investigated as a supplementary strategy to maintain power supply during high-risk periods while mitigating ignition risks \citep{piansky2024long}. Despite these advancements, there is limited research literature on the integration of social vulnerability considerations into wildfire risk management and infrastructure investment strategies. There is a critical need for frameworks that explicitly incorporate equity considerations, ensuring that mitigation efforts adequately protect and benefit the most vulnerable populations.

Power line undergrounding is not the only long-term approach for mitigating wildfires sparked from power lines; other options include vegetation management, covered conductors, and fast trip settings. This paper focuses on undergrounding power lines for wildfire risk mitigation because utilities are already actively doing this \citep{PGE,Blunt_2023}. In previous work, \cite{kody2022optimizing} showed undergrounding may be a good choice when balancing load shed and wildfire risk reduction. Furthermore, because  undergrounding provides a permanent and comprehensive solution to many of the risk factors contributing to ignitions, a power line undergrounding campaign is well-aligned with the goals of the Justice40 initiative.


\subsection{Social Vulnerability and Infrastructure Investments: The Case of Texas}

Texas offers a compelling context for examining the intersection of social vulnerability, infrastructure investments, and wildfire risk. The state's diverse demographic composition and varied socioeconomic landscapes contribute to differential exposures and capacities to cope with environmental hazards, including wildfires.

Studies have highlighted the inequitable impacts of environmental disasters and infrastructure failures across Texas communities. For instance, \citet{shah2023inequitable} document the disproportionate effects of power outages on minority populations during Winter Storm Uri in 2021, underscoring systemic vulnerabilities exacerbated by inadequate infrastructure resilience. 
In a similar report of Winter Storm Uri, these authors show that, adjusted for income level, census block groups with higher proportions of minority populations are more likely to experience a power outage \cite{Hsu_Taneja_Carvallo_Shah_2023}. This suggests that Texas residents' vulnerability to power outages is related to racial minority status, but interestingly somewhat unrelated to income. Our experiments with a simulated dataset corroborate this finding (discussed in \S\ref{subsection: ls under no intervention} referencing Appendix~\ref{appendix:income and ls}), which is why the case study described in \S\ref{sec:case_study} considers racial groups primarily instead of focusing on income groups.

While electric utilities in Texas have not historically used \gls*{psps} events, Texas has experienced more than 4,000 wildfires caused by power lines between 2011 and 2014 \citep{powerfiretexas}, 
making \gls*{psps} events a possible future solution. Most recently, the 2024 Smokehouse Creek fire, the largest wildfire in Texas history, was likely ignited by power infrastructure \citep{rose2024utility}, emphasizing the state's susceptibility to such hazards and the critical importance of effective risk mitigation strategies \citep{Penn_2024}. Furthermore, the increase in power outages due to ice storms and other extreme weather in Texas has led to increased popularity of line undergrounding as a method of preventing these blackouts \citep{macmillan2021longer}. The combination of climatic conditions conducive to wildfires and expansive power infrastructure networks necessitates proactive measures to prevent ignitions and protect communities.

\section{Methods}\label{sec:methods}
We first introduce the framework used to model line de-energization from \gls*{psps} events on an electric power transmission network. Then, we discuss modifications to this model to incorporate line undergrounding and equity considerations. 

\subsection{Network Description}
We consider an electric transmission network comprised of buses (nodes) connected by power lines (edges). 
%
For a given network, let $\mathcal{N}$ be the set of buses, $\mathcal{L}$ be the set of transmission lines, and $\mathcal{G}$ be the set of generators.
Let $\mathcal{T} = \{1, \ldots, T \}$ be the considered set of time indices over the period of a day. Let $\mathcal{D} = \{1, \ldots, D \}$ be the considered set of days. We define a 100~MVA per unit (p.u.)~base power. 
The following network parameters are provided for all lines $\ell \in \mathcal{L}$:
\begin{itemize}[noitemsep,nolistsep]
    \item $b_{\ell}$, line susceptance in p.u.,
    \item $\overline{f}_{\ell}$, the power flow limit in p.u., 
    \item $r_{\ell,d}$, the wildfire risk incurred if line $\ell$ is energized on day $d$,
    \item $n_{\ell, \text{to}}$ and $n_{\ell, \text{fr}}$,  \textit{to} and \textit{from} buses, respectively, denoting direction of positive power flow,
    \item $\overline{\delta}_{\ell}$ and $\underline{\delta}_{\ell}$, upper and lower voltage angle difference limits in radians, respectively,
\end{itemize}
The set $\mathcal{L}$ is further divided in to $\mathcal{L}^{high}_d$, $\mathcal{L}^{med}_d$, and $\mathcal{L}^{low}_d$ to indicate the set of lines that have high, medium, or low wildfire risk on day $d$, respectively. These categories are further described in Appendix~\ref{appendix: usgs data}.
For all generators $i \in \mathcal{G}$, define the parameters:
\begin{itemize}[noitemsep, nolistsep]
    \item $\overline{p}_i^g$ and $\underline{p}_i^g$, upper and lower power generation limits, respectively, in {p.u.}~for generator $i$,
\end{itemize}
For all buses $n \in \mathcal{N}$, define the parameters:
\begin{itemize}[noitemsep, nolistsep]
    \item $p_{n,d,t}^l$, power demand (load) at bus $n$ at time period $t \in \mathcal{T}$ on day $d \in \mathcal{D}$ in p.u.,
    \item $\mathcal{G}^n$, the set of generators located at bus $n$,
    \item $\mathcal{L}^{n, {\text{to}}}$ and $\mathcal{L}^{n, {\text{fr}}}$, the subset of lines $\ell \in \mathcal{L}$ with bus $n$ as the designated \textit{to} bus and bus $n$ as the designated \textit{from} bus, respectively.
\end{itemize}
The operation of the network during a multi-time-period PSPS event is characterized by the following set of variables using the B$\Theta$ representation of the DC power flow model:
\begin{itemize}[noitemsep, nolistsep]
    \item $p_{i,d,t}^g$, power generated at unit $i \in \mathcal{G}$ at time period $t \in \mathcal{T}$ on day $d \in \mathcal{D}$ in p.u.,
    \item $\theta_{n,d,t}$, voltage angle at bus $n \in \mathcal{N}$ at time period $t \in \mathcal{T}$ on day $d \in \mathcal{D}$ in radians,
    \item $p_{n,d,t}^{s}$, load shedding at bus $n \in \mathcal{N}$ at time period $t \in \mathcal{T}$ on day $d \in \mathcal{D}$ in p.u.,
    \item $f_{\ell,d,t}$, power flowing 
    along line $\ell \in \mathcal{L}$ at time period $t \in \mathcal{T}$ on day $d \in \mathcal{D}$ in p.u.,
    \item $z_{\ell,d} \in \{0,1\}$, state of energization of line $\ell \in \mathcal{L}^{\text{high}}_d$ and $\ell \in \mathcal{L}^{\text{med}}_d$ on day $d \in \mathcal{D}$. If $z_{\ell,d}=0$, then line $\ell$ is de-energized, and  if $z_{\ell,d}=1$, then line $\ell$ is energized. Note that the line's energization state is constant for all $t \in \mathcal{T}$ on day $d \in \mathcal{D}$. Transmission lines that are de-energized for safety reasons, such as a \gls*{psps} event, need to be visually inspected by a grounds crew before being re-energized, so intra-day de-energizations are not desirable. Models for transmission line re-energization have been explored in ~\cite{rhodes2022co}. For all $\ell \in \mathcal{L}^{\text{low}}_d$, $z_{\ell,d} = 1$. Let $\mathcal{L}^{\text{switch}}_d = \mathcal{L}^{\text{high}}_d \cup \mathcal{L}^{\text{med}}_d$.
\end{itemize}

\subsection{Operational and Physical Constraints}\label{appendix:base model}
We define the constraints for the DC Optimal Transmission Switching Problem (DC-OTS) in Equation \eqref{eq: dcots}, $\ \forall d \in \mathcal{D}, \ \forall t \in \mathcal{T}$.
\begin{model}[t]
\caption{DC-OTS Constraints}
\label{eq: dcots}
\begin{subequations} 
\vspace{-0.2cm}
    \begin{align}
        &\underline{p}_{i}^g \leqslant p_{i,d,t}^g \leqslant \overline{p}_{i}^g, &&\hspace{-6.9em}\!\forall i \in \mathcal{G}, \label{subeq: gen limits} 
        \\
        &0 \leqslant p_{n,d,t}^{s} \leqslant p_{n,d,t}^l, 
        &&\hspace{-6.9em}\!\forall n \in \mathcal{N}, \label{subeq: load shed limits} 
        \\ 
        &-\overline{f}_\ell z_{\ell,d} \leqslant f_{\ell,d,t} \leqslant \overline{f}_\ell z_{\ell,d}, 
        &&\hspace{-6.9em}\!\forall \ell \in \mathcal{L}^{\text{switch}}_d, \label{subeq: power flow limits switching}
        \\
        &-\overline{f}_\ell \leqslant f_{\ell,d,t} \leqslant \overline{f}_\ell, 
        &&\hspace{-6.9em}\!\forall \ell \in \mathcal{L} \setminus \mathcal{L}^{\text{switch}}_d, \label{subeq: power flow limits} 
        \\
        &\underline{\delta}_\ell  \leqslant \theta_{{n_{\ell, \text{fr}}},d,t} - \theta_{{n_{\ell, \text{to}}},d,t} \leqslant \overline{\delta}_\ell, 
        &&\hspace{-6.9em}\!\forall \ell \in \mathcal{L}\setminus\mathcal{L}^{\text{switch}}_d, \label{subeq: voltage angle} 
        \\
        &\theta_{{n_{\ell, \text{fr}}},d,t} \!-\! \theta_{{n_{\ell, \text{to}}},d,t} \!\geqslant\! \underline{\delta}_\ell z_{\ell,d} \!+\! \underline{M}(1\!-\!z_{\ell,d}),
        &&\hspace{-5.15em}\!\forall \ell \in \mathcal{L}^{\text{switch}}_d, \label{subeq: voltage angle switching lower} 
        \\
        &\theta_{{n_{\ell, \text{fr}}},d,t} \!-\! \theta_{{n_{\ell, \text{to}}},d,t} \!\leqslant\! \overline{\delta}_\ell z_{\ell,d} \!+\! \overline{M}(1\!-\!z_{\ell,d}), 
        &&\hspace{-5.15em}\!\forall \ell \in \mathcal{L}^{\text{switch}}_d, \label{subeq: voltage angle switching upper} 
        \\
        & f_{\ell,d,t} \!\geqslant\! -b_\ell(\theta_{{n_{\ell, \text{fr}}},d,t} \!-\! \theta_{{n_{\ell, \text{to}}},d,t}) \!+\! |b_\ell|\underline{M}(1\!-\!z_{\ell,d}), \nonumber \\
        &\hspace{14.1em}\forall \ell \in \mathcal{L}^{\text{switch}}_d, \label{subeq: power flow switching lower} 
        \\
        &f_{\ell,d,t} \!\leqslant\! -b_\ell(\theta_{{n_{\ell, \text{fr}}},d,t} \!-\! \theta_{{n_{\ell, \text{to}}},d,t}) \!+\! |b_\ell|\overline{M}(1\!-\!z_{\ell,d}), \nonumber \\
        &\hspace{14.1em}\forall \ell \in \mathcal{L}^{\text{switch}}_d, \label{subeq: power flow switching upper} 
        \\
        &f_{\ell,d,t} = -b_\ell(\theta_{{n_{\ell, \text{fr}}},d,t} - \theta_{{n_{\ell, \text{to}}},d,t}),
        &&\hspace{-6.9em}\!\forall \ell \in \mathcal{L}\setminus\mathcal{L}^{\text{switch}}_d, \label{subeq: power flow}
        \\
        &\sum_{\ell \in \mathcal{L}^{n, \text{fr}}} f_{\ell,d,t} \!-\! \sum_{\ell \in \mathcal{L}^{n, \text{to}}} f_{\ell,d,t} \!=\! \sum_{i \in \mathcal{G}^n} p_{i,d,t}^g \!-\! p_{n,d,t}^l \!+\! p_{n,d,t}^{s},\nonumber \\
        &\hspace{14.1em}\forall n \in \mathcal{N}, \label{subeq: power balance}
    \end{align}
\end{subequations}
\end{model}

In \eqref{eq: dcots}, \eqref{subeq: gen limits} enforces lower and upper generation limits, \eqref{subeq: load shed limits} constrains any load shedding to be nonnegative and less than the load demanded at that time at that bus,~\eqref{subeq: power flow limits switching} and \eqref{subeq: power flow limits} enforce line flow limits, \eqref{subeq: voltage angle} and \eqref{subeq: voltage angle switching lower}/\eqref{subeq: voltage angle switching upper} constrain angle differences across lines, \eqref{subeq: power flow switching lower}/\eqref{subeq: power flow switching upper} and \eqref{subeq: power flow} model the DC power flow approximation, and \eqref{subeq: power balance} ensures power balance at all buses in the network. In equations \eqref{subeq: voltage angle switching lower}, \eqref{subeq: voltage angle switching upper}, \eqref{subeq: power flow switching lower} and \eqref{subeq: power flow switching upper}, $\overline{M}$ and $\underline{M}$ are big-M constants set to $2\pi$ and $-2\pi$ respectively for the numerical experiments in this paper. Note that an energized line (i.e., $z^{\ell}=1$), simplifies \eqref{subeq: voltage angle switching lower}/\eqref{subeq: voltage angle switching upper} to \eqref{subeq: voltage angle} and \eqref{subeq: power flow switching lower}/\eqref{subeq: power flow switching upper} to \eqref{subeq: power flow}.

We then constrain the total risk from all energized above-ground lines to be below a given threshold, $R_{PSPS}$, on each day:
\begin{align} \label{const: switching threshold}
    \sum_{\ell \in \mathcal{L}} z_{\ell,d} r_{\ell,d} \leq R_{PSPS} \qquad \forall d \in \mathcal{D}.
\end{align}

\subsection{Line Undergrounding Formulation}\label{sec:undergrounding}
Let the subset $\mathcal{L}^{\text{harden}} \subseteq \mathcal{L}$ be the set of lines that are candidates for hardening/maintenance where $\mathcal{L}^\text{harden} = \mathcal{L}^\text{high}_d \bigcup_{d \in \mathcal{D}} \mathcal{L}^\text{med}_d$.
Define $\phi_\ell^{\text{ug}}$ as the cost of undergrounding line $\ell$ in millions of dollars per line.
For all $\ell \in \mathcal{L}^{\text{harden}}$, we introduce the variable $y_\ell \in \{0,1\}$, which indicates whether a line has been undergrounded ($y_\ell=1$) or not ($y_\ell=0$). Note for all $\ell \in \mathcal{L}\setminus\mathcal{L}^{\text{harden}}$, we set $y_\ell=0$. In our model, we assume that the entire length of line $\ell$ is undergrounded.

There is no benefit to simultaneously undergrounding and de-energizing a line since undergrounding a line already reduces the line's risk to zero. Thus, we impose the following constraints to prevent simultaneously de-energizing and hardening a line: 
\begin{align} 
    z_{\ell,d} = y_\ell, \qquad &\forall \ell \in \mathcal{L}^{\text{high}}_d, \ \forall d \in \mathcal{D}.\label{const: Line high}
    \\
    (1-z_{\ell,d}) + y_\ell \leqslant 1, \qquad &\forall \ell \in \mathcal{L}^{\text{med}}_d, \ \forall d \in \mathcal{D}.\label{const: Line med}
\end{align}
Equation \eqref{const: Line high} ensures lines in the highest-risk category are either de-energized 
or undergrounded, removing the risk by one method or the other. Alternatively, \eqref{const: Line med} allows for the lines in the medium-risk category to optionally be de-energized, undergrounded, or left above ground and energized. Regardless, in both categories, undergrounded lines will remain energized. 

We introduce a budget $B$ (in millions of dollars) to limit the resources available for line undergrounding. We use a fixed cost of $\$7$ million per mile for undergroundinglines in the network \citep{cpuc_undergrounding}. The total cost of undergrounding line $\ell \in \mathcal{L}$ is precomputed as $\phi_\ell^{\text{ug}}$, defined above. This introduces the following constraint that enforces that the total cost of undergrounding must not exceed budget $B$:
\begin{align} \label{const: budget}
    \sum_{\ell \in \mathcal{L}^{\text{harden}}} \phi_\ell^{\text{ug}} y_\ell \leqslant B.
\end{align}
We can reformulate equation~\eqref{const: switching threshold} to account for both the risk reduction associated with de-energizing lines as well as with line undergrounding:
\begin{align} \label{const: undergrounding threshold}
    \sum_{\ell \in \mathcal{L}} r_{\ell,d} (z_{\ell,d} - y_\ell) \leqslant R_{PSPS}, \qquad \forall d \in \mathcal{D}.
\end{align}

\subsection{Baseline Model} \label{sec: baseline}
The baseline model reflects the scenario where undergrounding decisions are made \textit{without} considering the vulnerability or demographics of affected populations.
That is, the objective of the baseline model is to minimize total load shed in the network. Let $P^l$ be the total demand in the network over all time periods:
\begin{align}
    P^l = \sum_{d \in \mathcal{D}}\sum_{t \in \mathcal{T}}\sum_{n \in \mathcal{N}} p_{n,t,d}^l, \nonumber
\end{align}
and let $P^s$ be the total load shed in the network over all time periods:
\begin{align}
    P^s = \sum_{d \in \mathcal{D}}\sum_{t \in \mathcal{T}}\sum_{n \in \mathcal{N}} p_{n,t,d}^{s}. \nonumber
\end{align}
Now, let $S$ represent the fraction of load shed in the network, given by:
\begin{align} \label{objective func}
    S = \frac{P^s}{P^l}.
\end{align}
We define the baseline model as:
\begin{equation}\tag{BL}
\begin{aligned}
&\min\limits_{p^g, \theta, f, p^{s}, z, y} \ \eqref{objective func} \ \ \text{s.t.} \ \eqref{eq: dcots},~\eqref{const: Line high} - \eqref{const: undergrounding threshold}. 
\end{aligned}
\label{obj: base-PSPS}
\end{equation}

\subsection{Fairness Models}\label{subsec: equity considerations}
In this subsection we present optimization models that select lines to underground \textit{while considering} the vulnerability and/or demographics of affected populations.
Let $\mathcal{M}$ be the set of considered demographic groups (e.g., income, race).
Let $\gamma_{n,m}^{\text{grp}}$ represent the percentage of the population at bus $n \in \mathcal{N}$ that belongs to demographic group $m \in \mathcal{M}$.
Let $\gamma_{n}^{\text{vuln}}$ be the percentage of population at bus $n$ that belongs to vulnerable populations, which is defined based on the considered metric (e.g., CEJST or SVI).

We assume that the load requested at a bus is allocated proportionally to the groups at that bus. 
That is, if group A is 20\% of the population at bus $n$, we attribute 20\% of the load at bus $n$ to group A.
We define the total load demanded by each group $m \in \mathcal{M}$ across all buses and all considered time periods as $P^{l}_{m}$:
\begin{align}
    P^{l}_{m} &= \sum_{d \in \mathcal{D}}\sum_{t \in \mathcal{T}}\sum_{n \in \mathcal{N}} \gamma_{n,m}^{\text{grp}} p_{n,t,d}^l \nonumber.
\end{align}
Similarly, we assume that load shed at a bus proportionally affects the groups at that bus.
Let $P^{s}_{m}$ for all groups $m \in \mathcal{M}$ be the total load shed experienced by group $m$ across all buses and all considered time periods:
\begin{align}
    P^{s}_{m} &= \sum_{d \in \mathcal{D}}\sum_{t \in \mathcal{T}}\sum_{n \in \mathcal{N}} \gamma_{n,m}^{\text{grp}} p_{n,t,d}^{s}. \nonumber
\end{align}
The demographics and vulnerability of the population at each bus in the network are determined using census data, where we assign each census tract to the nearest load bus.
See Appendix \ref{appendix: bus features} for details on the methodology we employ to attribute demographic data to each bus.

Next, we introduce three categories of optimization model that aim to incorporate fairness considerations into undergrounding decisions:
\begin{enumerate}[noitemsep,nolistsep]
    \item {\sc Policy} constraint models, which preferentially allocate resources to groups with \textit{vulnerability} characteristics by imposing proportionate policy-level constraints,
    \item {\sc Equity} objective models, which seek to minimize disparate impacts of load shedding across different \textit{demographic} groups through specially designed objective functions, and
    \item {\sc Policy} constraint and {\sc Equity} objective model, which combines the two above goals.
\end{enumerate}
Table \ref{table: model summary} can be used as reference for the considered models in this study.

\subsubsection{{\sc Policy} Constraint Models}\label{sec:j40 constraint}
As discussed in Section \ref{sec: introduction}, the goal of the Justice40 initiative is \textcolor{red}{for 40\% of the ``benefit" induced by federal climate investment to be accessible to vulnerable groups}. 
\textcolor{red}{The {\sc Policy} constraint models are intended to enforce this prioritization.
\sout{The aim of the {\sc Policy} constraint models are to ensure 40\% of benefits are realized by vulnerable groups.
We enforce this allocation through the inclusion of an additional optimization constraint.}}
In the context of a wildfire risk mitigation scenario, the ``benefit'' to vulnerable groups could be interpreted as (1) the budget allocated to vulnerable groups, or (2) the load shed reduction experienced by vulnerable groups.
Consequently, there are two versions of the {\sc Policy} constraint model: one that enforces the benefit to vulnerable groups through a budget constraint, and another that enforces it through a load shed reduction constraint.

We first introduce the version of the {\sc Policy} constraint model that ensures 40\% of the total budget aids vulnerable populations.
Recall from Section \ref{sec:undergrounding} that $\phi_\ell^{\text{ug}}$ is the cost to underground line $\ell \in \mathcal{L}$.
We assume that half of this investment is attributed to each of the populations at the terminal buses.
That is, the dollar amount of benefit is split equally between the populations at the two buses of the line, $n_{\ell, \text{to}}$ and $n_{\ell, \text{fr}}$.
We then define the budget spent on vulnerable populations as:
\begin{align}
    \phi^{\text{vuln}} = \sum_{\ell \in \mathcal{L}} y_\ell \frac{\phi_\ell^{\text{ug}}}{2}\left(\gamma_{n_{\ell, \text{to}}}^{\text{vuln}} + \gamma_{n_{\ell, \text{fr}}}^{\text{vuln}}\right) \nonumber.
\end{align}
To ensure 40\% of the budget is allocated to vulnerable populations  based on the Justice40 initiative, we enforce the following constraint:
\begin{align}\label{const: j40 bud}
    \phi^{\text{vuln}} \geq 0.4\cdot B.
\end{align}
Now, we can introduce the first version of the {\sc Policy} constraint model, which ensures that 40\% of the budget is allocated to vulnerable populations:
\begin{equation}
\begin{aligned}
&\min\limits_{p^g, \theta, f, p^{s}, z, y} \ \eqref{objective func} \ \ \text{s.t.} \ \eqref{eq: dcots},~\eqref{const: Line high} - \eqref{const: undergrounding threshold}, \eqref{const: j40 bud}. 
\end{aligned}
\label{obj: base-j40_bud}
\end{equation}

Next, we introduce the second version of the {\sc Policy} constraint model that ensures 40\% of the total load shed reduction affects vulnerable populations.
We define the total load shed seen by vulnerable populations over the considered time period as:
\begin{align}
    P^{s, \text{vuln}} = \sum_{d \in \mathcal{D}}\sum_{t \in \mathcal{T}}\sum_{n \in \mathcal{N}} p_{n,d,t}^{s} \ \gamma_{n}^{\text{vuln}} \nonumber.
\end{align}
Let $P^s_{BL-M0}$ and $P^{s, \text{vuln}}_{BL-M0}$ represent the total load shed and the load shed experienced by vulnerable populations \textcolor{red}{respectively} under the baseline model \eqref{obj: base-PSPS} with no budget, where the subscript \textit{BL-M0} references the descriptors in Table \ref{table: model summary}.
Recall from Section \ref{sec: baseline} that $P^s$ is the total load shed in the network.
Now, we can introduce the optimization constraint to enforce that 40\% of load shed reduction is seen by vulnerable populations:
\begin{align}\label{const: j40 ls}
    \frac{P^{s, \text{vuln}}_{BL-M0} - P^{s, \text{vuln}}}{P^s_{BL-M0} - P^s} \geq 0.4.
\end{align}
The second version of the {\sc Policy} constraint model that ensures that 40\% of the load shed reduction affects vulnerable populations is then defined as:
\begin{equation}
\begin{aligned}
&\min\limits_{p^g, \theta, f, p^{s}, z, y} \ \eqref{objective func} \ \ \text{s.t.} \ \eqref{eq: dcots},~\eqref{const: Line high} - \eqref{const: undergrounding threshold}, \eqref{const: j40 ls}. 
\end{aligned}
\label{obj: base-j40_ls}
\end{equation}

\subsubsection{{\sc Equity} Objective Model}
The {\sc Equity} objective models aim to minimize the disparity in proportional load shed seen across groups, i.e. minimizing the maximum proportion of load shed seen by any one group. To define the {\sc Equity} objective, we define the auxiliary variable, $\alpha$, denoting the maximum proportion of load demanded that is shed across each group. That is, for each group $m \in M$, we constrain
\begin{align}\label{const:alph_m}
    \alpha \geq \frac{\text{P}^{s}_{m}}{\text{P}^{l}_{m}}.
\end{align}

%
We define the following \gls*{mip} to minimize the maximum percent load shed across groups: 
\begin{align}\label{obj: mm_pct}
    &\min\limits_{p^g, \theta, f, p^{s}, z, y} \ \alpha \ \ \text{s.t.} \ \eqref{eq: dcots},~\eqref{const: Line high} - \eqref{const: undergrounding threshold}, ~\eqref{const:alph_m}.
\end{align}

This method can result in extraneous load shed. For example, load shed occurring at a bus that has no associated population belonging to the most affected group is not accounted for by the objective. To maintain a focus on improving outcomes for the most impacted groups, a post-processing step minimizes total load shed, \eqref{objective func}, in the network after setting the binary decisions found from models E-M6 through E-M10 as discussed in Table~\ref{table: model summary}. This post-processing re-optimizes over the continuous variables.

\subsubsection{{\sc Policy} Constraints with {\sc Equity} Objectives}
We include two additional models that incorporate both {\sc policy} constraints as well as {\sc equity} objectives. 
With the inclusion of the {\sc policy} constraint considering budget, the objective becomes~\eqref{obj: mm_pct_w_j40_bud-PSPS}:
\begin{equation}
    \begin{aligned}
        &\min\limits_{p^g, \theta, f, p^{s}, z, y} \ \alpha \ \ \text{s.t.} \ \eqref{eq: dcots},~\eqref{const: Line high} - \eqref{const: undergrounding threshold},~\eqref{const: j40 bud}, ~\eqref{const:alph_m},
    \end{aligned}
    \label{obj: mm_pct_w_j40_bud-PSPS}
\end{equation}
Next, we consider a constraint on proportional load shed reduction along with the {\sc equity} objective,~\eqref{obj: mm_pct_w_j40_ls-PSPS}
\begin{equation}
    \begin{aligned}
        &\min\limits_{p^g, \theta, f, p^{s}, z, y} \ \alpha \ \ \text{s.t.} \ \eqref{eq: dcots},~\eqref{const: Line high} - \eqref{const: undergrounding threshold},~\eqref{const: j40 ls}, ~\eqref{const:alph_m}.
    \end{aligned}
    \label{obj: mm_pct_w_j40_ls-PSPS}
\end{equation}

\subsection{Model Descriptions}
Results presented in the body of this paper cover 11 different models each for 5 different budgets, representing various combinations of {\sc Policy} constraints and {\sc Equity} objectives. A summary of models is given in Table~\ref{table: model summary}.
\begin{enumerate}[nolistsep,noitemsep]
    \item Load shed outcomes under no intervention, given by model \eqref{obj: base-PSPS} when the budget parameter $B$ is set to zero. We will call this our \textit{no-budget} model, and is denoted by BL-M0.
    \item Load shed outcomes when allocating budget for undergrounding without considering any vulnerability or demographic factors (i.e., no {\sc Policy} constraints or {\sc Equity} objective). This model is given by \eqref{obj: base-PSPS} for nonzero budgets. In tables and plots, this model is denoted by BL-M1 to emphasize that this model is the baseline for a specified budget $B$. In this results section, we look at budget $B = 1000$ million USD, but the Appendix shows results for other budgets.
    \item Load shed outcomes when allocating budget for undergrounding when using an objective that minimizes total load shed, with {\sc Policy} constraints to preferentially allocate 40\% of the budget \eqref{obj: base-j40_bud} (M2, M4) or 40\% of the load shed reduction \eqref{obj: base-j40_ls} (M3, M5) relative to the no-budget model to tracts designated as disadvantaged by the \mbox{\gls{cejst}} and \gls{svi}, respectively.
    \item Load shed outcomes when allocating budget for undergrounding when using the {\sc Equity} objective, either without \eqref{obj: mm_pct} or with (\eqref{obj: mm_pct_w_j40_bud-PSPS} or \eqref{obj: mm_pct_w_j40_ls-PSPS}) {\sc Policy} constraints by budget (\eqref{const: j40 bud} in E-M7 and E-M9) and load shed reduction (\eqref{const: j40 ls} in E-M8 and E-M10), respectively. These models use either the \mbox{\gls{cejst}} (E-M7 and E-M8) or \mbox{\gls{svi}} (E-M9 and E-M10) datasets. Models of this type are given the prefix ``E-" for equity, appended by the model label M6-M10 for a pre-specified budget; see Table \ref{table: model summary} for specifics of each model type.
\end{enumerate}

\section{Texas Wildfire Case Study}\label{sec:case_study}
In this section, we discuss the demographic distribution and network configuration specific to this case study. The case study models a synthetic grid covering a similar area to the \gls{ercot} transmission network in Texas. 
 
\subsection{Synthetic Texas Network}\label{subsec: ercot}
We use the synthetic Texas7k transmission network test case, developed by the Texas A\&M PERFORM group~\citep{xu2017creation}. This test case provides a realistic approximation of the area covered by the \gls*{ercot} \citep{xu2017creation,xu2018synthetic} while not disclosing critical energy infrastructure information~\citep{ferc_ceii}.

\begin{figure*}[th!]
\centering
\begin{subfigure}{.31\textwidth}
  \centering
    \includegraphics[width=.99\linewidth]{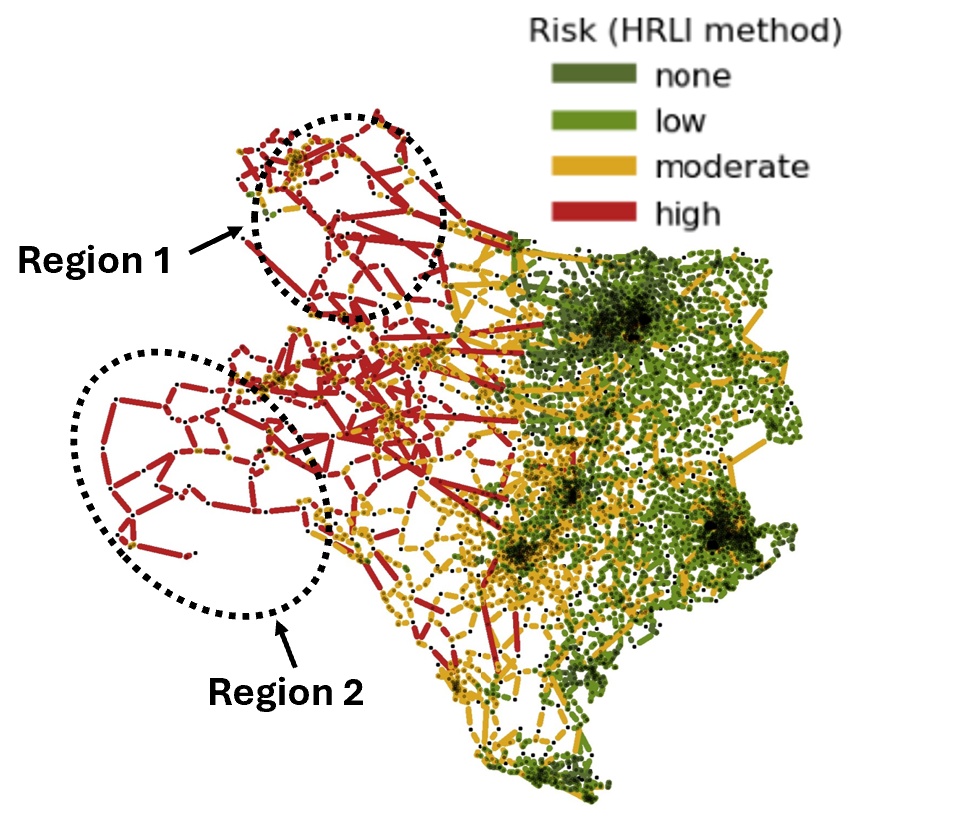} 
    \caption{Risk per line in the Texas7k network}
  \label{fig: subfig: risk per line Texas}
\end{subfigure}
\begin{subfigure}{.31\textwidth}
  \centering
  \includegraphics[width=.99\linewidth]{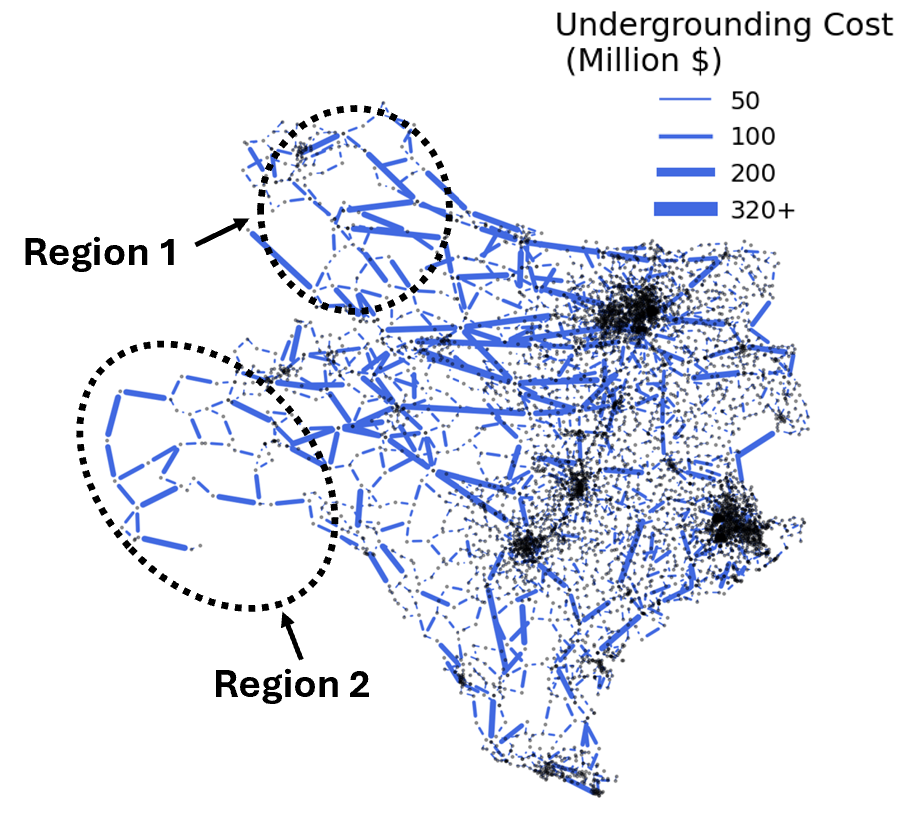}  
  \caption{Cost to underground each line}
  \label{fig: subfig: cost to underground}
\end{subfigure}
\begin{subfigure}{.35\textwidth}
  \centering
  \includegraphics[width=.99\linewidth]{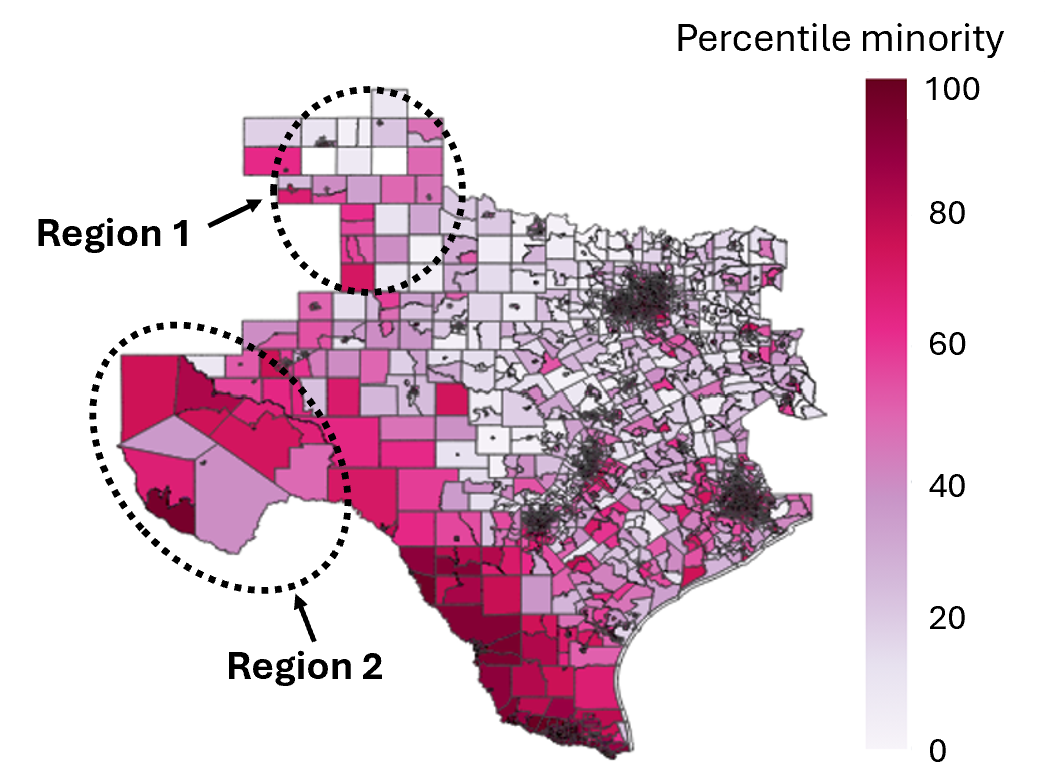}  
  \caption{Percent of tract that is a racial minority}
  \label{fig: subfig: racial minority}
\end{subfigure}
\caption{This figure shows some vulnerability characteristics of the synthetic Texas7k network where the circled regions show the overlap in vulnerabilities between these three metrics (wildfire risk, high cost to underground, and racial minority status). These areas also have a lower likelihood of being selected for power line undergrounding due to lower population density (not pictured).}
\label{fig: Texas7k profile}
\end{figure*}

\subsection{Simulating Features of the Texas 7k Network}
The Texas7k network gives us information about the synthetic transmission network structure (e.g., network topology, loads, line impedances, etc.). To understand the wildfire ignition risk and demographic profile of different parts of the Texas7k network, we augment the dataset with three other types of data:
\begin{enumerate}[nolistsep,noitemsep]
    \item Demographic data from the 2010 US Census \citep{COP_data,demographic_data}.
    \item Vulnerability index data at the US Census tract level 
    \citep{CDC_SVI_data,ScreeningTool_2022}.
    \item \gls{usgs} wildfire risk data \citep{USGS_WFPI}.
\end{enumerate}
Using these datasets, we assign demographic information and vulnerability metrics for each bus in the network. For this problem, we \textcolor{red}{generally} define the groups $\mathcal{M}$ as the five census defined racial and ethnic groups; Hispanic, White, Black, Indigenous, and Asian populations.
For each line in the network, we assign a risk of igniting a wildfire if the line is energized for each day in the simulation period.
Reference the Appendix for more details on demographic data (\ref{subapp: census data}), CEJST (\ref{subappendix: j40}), SVI (\ref{subappendix: svi}), and the USGS data (\ref{appendix: usgs data}). 

This augmented dataset enables visualizing risk, demographic, and vulnerability data for the synthetic Texas7k transmission network.
Figure \ref{fig: Texas7k profile} shows vulnerability characteristics of the synthetic network where the circled regions highlight the overlap in vulnerability between these four metrics (wildfire risk, high cost to underground, low-income status, and racial minority status). We can see that the western- and northern-most parts of Texas are likely most susceptible to power outages from \gls*{psps} events (i.e., high wildfire risk, but lower ability to cope due to lower income) as well as lower likelihood of being selected for power line undergrounding due to lower population density (not pictured) and the high cost to underground the lines. Hence, we can anticipate some difficulty allocating load shed relief to populations in these northwest regions. This hypothesis will be corroborated in \S\ref{sec:results}.
\begin{figure*}[ht]
    \centering
    \begin{subfigure}[t]{.32\textwidth}
        \centering 
        \captionsetup{width=.95\linewidth}
        \includegraphics[width=.95\textwidth]{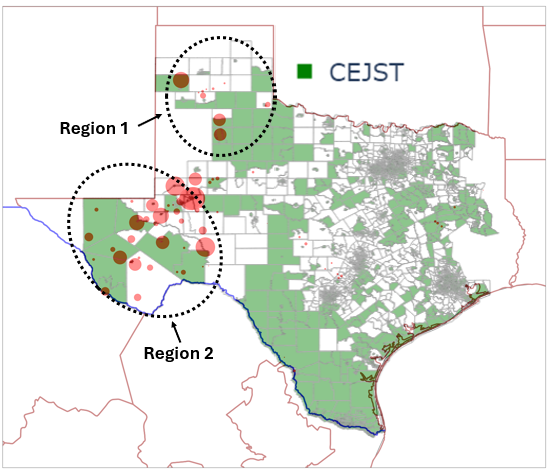}
        \caption{Approximately 40\% of all Texas census tracts are vulnerable per the {\sc CEJST}, and only 37\% of Texas census tracts experiencing load shed are {\sc CEJST} tracts.}
        \label{fig:subfig:j40_ls}
    \end{subfigure}
    \begin{subfigure}[t]{.32\textwidth}
        \centering 
        \captionsetup{width=.95\linewidth}
        \includegraphics[width=.95\textwidth]{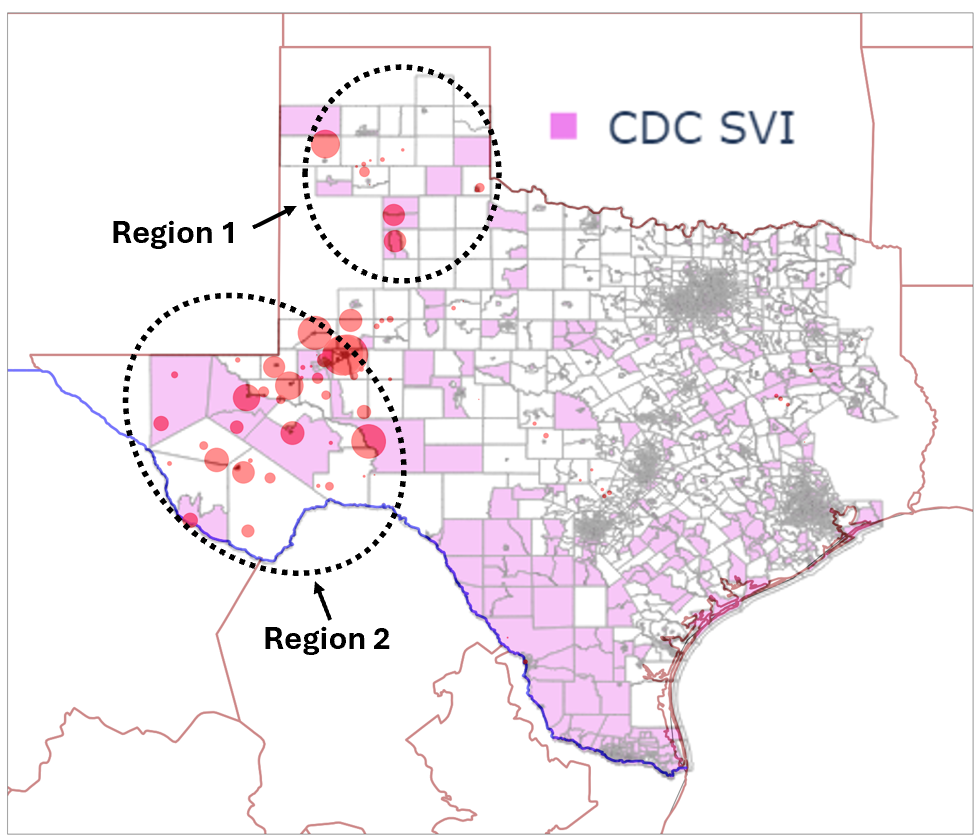}
        \caption{Approximately 38\% of all Texas census tracts are vulnerable per the {\sc SVI}, and only 28\% of Texas census tracts experiencing load shed are {\sc SVI} tracts.}
        \label{fig:subfig:svi_ls}
    \end{subfigure}
    \begin{subfigure}[t]{.32\textwidth}
        \centering 
        \captionsetup{width=.95\linewidth}
        \includegraphics[width=.95\textwidth]{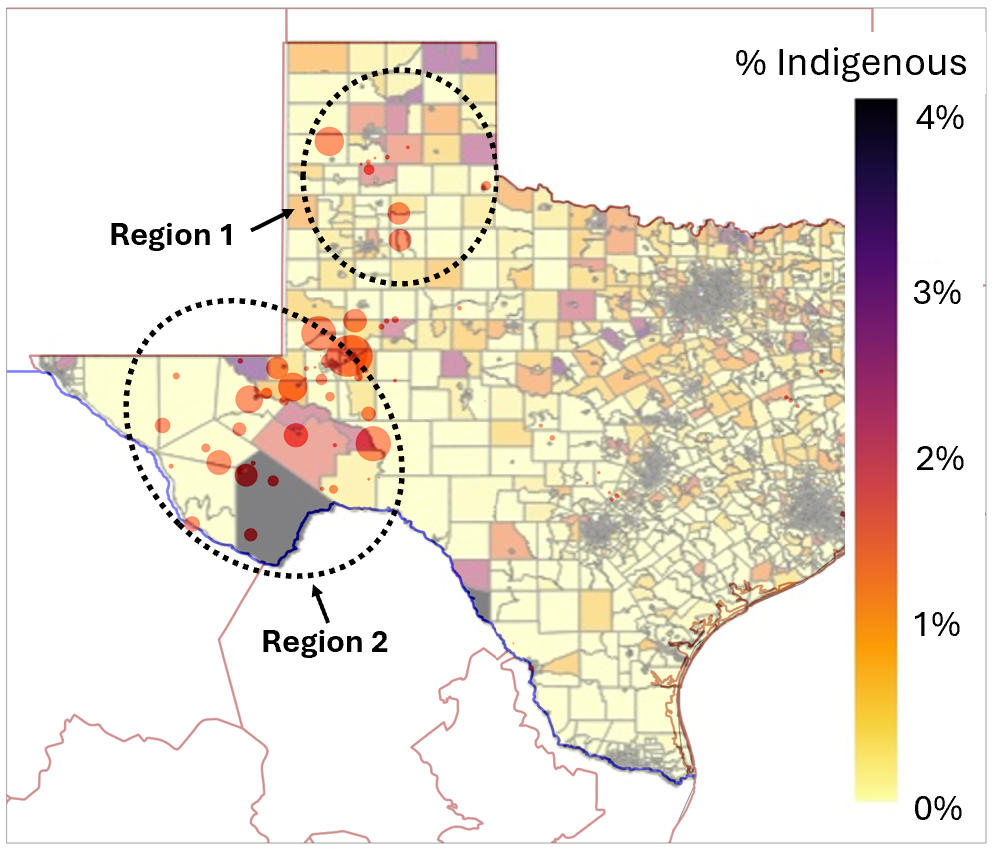}
        \caption{Census tracts with relatively larger (although still minority) Indigenous populations are high-risk, and have load shed. \newline}
        \label{fig:subfig:indig_ls}
    \end{subfigure}
    \vspace{-.3 cm}
    \caption{Load shed (in the absence of undergrounding decisions and equity considerations, BL-M0), visualized with red bubbles on the map, occurs almost exclusively in the vulnerable areas highlighted in Figure \ref{fig: Texas7k profile}. 
    }
    \label{fig: svi plots with load shedding}
\end{figure*}

\section{Results}\label{sec:results}
In this section, we discuss the impacts of equitable resource allocation policies on reducing load shed to vulnerable populations, particularly for Hispanic and Indigenous populations. We will consider four model types, discussed in each of the next four subsections, \S \ref{subsection: ls under no intervention}-\S\ref{subsec: group fair}. A summary of the models discussed below can be seen in Table \ref{table: model summary}.

For this case study, we show that models including the {\sc Equity} object more effectively reduces the proportional burden of load shed on Indigenous populations compared to using {\sc policy} constraints alone. 

\begin{sidewaystable*}[!htbp]
\centering
\begin{tabular}{|
>{\columncolor[HTML]{EFEFEF}}l |l|ll|l|l|}
\hline
\multicolumn{1}{|c|}{\cellcolor[HTML]{EFEFEF}}                                 & \multicolumn{1}{c|}{\cellcolor[HTML]{EFEFEF}}                                     & \multicolumn{2}{c|}{\cellcolor[HTML]{EFEFEF}\textbf{Policy Constraint}} & \multicolumn{1}{c|}{\cellcolor[HTML]{EFEFEF}} & \multicolumn{1}{c|}{\cellcolor[HTML]{EFEFEF}}                         \\ \cline{3-4} 
\multicolumn{1}{|c|}{\multirow{-2}{*}{\cellcolor[HTML]{EFEFEF}\textbf{Model}}} & \multicolumn{1}{c|}{\multirow{-2}{*}{\cellcolor[HTML]{EFEFEF}\textbf{Objective (\(\min\))}}} & \multicolumn{1}{c|}{\cellcolor[HTML]{EFEFEF}\textbf{Vulnerability Index}} & \multicolumn{1}{c|}{\cellcolor[HTML]{EFEFEF}\textbf{Constraint Type}} & \multicolumn{1}{|c|}{\multirow{-2}{*}{\cellcolor[HTML]{EFEFEF}\textbf{Budget}}} & \multicolumn{1}{|c|}{\multirow{-2}{*}{\cellcolor[HTML]{EFEFEF}\textbf{Eq.}}}   \\ \hline
\textbf{BL-M0}                                                         
    & total load shed \eqref{objective func}                                                                   &  \multicolumn{1}{l|}{N/A}                                                  & ~N/A & 0.0 & \eqref{obj: base-PSPS}\\ \hline

\textbf{BL-M1}                                                         
    & total load shed \eqref{objective func}                                                                   &  \multicolumn{1}{l|}{N/A}                                                  & ~N/A  & 1000.0 & \eqref{obj: base-PSPS}                                                                      \\ \hline
\textbf{M2}                                                                 & total load shed \eqref{objective func}                                                                          & \multicolumn{1}{l|}{CEJST}                                            & \begin{tabular}[c]{@{}l@{}}Prop. Budget Expenditure \eqref{const: j40 bud}\end{tabular}  & 1000.0 & \eqref{obj: base-j40_bud}   \\ \hline
\textbf{M3}                                                                    & total load shed\eqref{objective func}                                                                           & \multicolumn{1}{l|}{CEJST}                                            & \begin{tabular}[c]{@{}l@{}}Proportional Load Shed Reduction~\eqref{const: j40 ls}\end{tabular} & 1000.0 & \eqref{obj: base-j40_ls}   \\ \hline
\textbf{M4}                                                                    & total load shed \eqref{objective func}                                                                           & \multicolumn{1}{l|}{SVI}                                 & \begin{tabular}[c]{@{}l@{}}Prop. Budget Expenditure~\eqref{const: j40 bud}\end{tabular}  & 1000.0 & \eqref{obj: base-j40_bud}   \\ \hline
\textbf{M5}                                                                    & total load shed \eqref{objective func}                                                                          & \multicolumn{1}{l|}{SVI}                                 & \begin{tabular}[c]{@{}l@{}}Proportional Load Shed Reduction~\eqref{const: j40 ls}\end{tabular} & 1000.0  & \eqref{obj: base-j40_ls}  \\ \hline

\textbf{E-M6}                                                             & \(\max\) \% load shed across groups \eqref{const:alph_m}                                                   & \multicolumn{1}{l|}{N/A}                                                  & ~N/A      & 1000.0 & \eqref{obj: mm_pct}                                                                     \\ \hline
\textbf{E-M7}                                                             & \(\max\) \% load shed across groups \eqref{const:alph_m}                                                   & \multicolumn{1}{l|}{CEJST}                                            & \begin{tabular}[c]{@{}l@{}}Prop. Budget Expenditure~\eqref{const: j40 bud}\end{tabular} & 1000.0  & \eqref{obj: mm_pct_w_j40_bud-PSPS}   \\ \hline
\textbf{E-M8}                                                            & \(\max\) \% load shed across groups \eqref{const:alph_m}                                                 & \multicolumn{1}{l|}{CEJST}                                            & \begin{tabular}[c]{@{}l@{}}Proportional Load Shed Reduction~\eqref{const: j40 ls}\end{tabular} & 1000.0  & \eqref{obj: mm_pct_w_j40_ls-PSPS}  \\ \hline
\textbf{E-M9}                                                            & \(\max\) \% load shed across groups \eqref{const:alph_m}                                                    & \multicolumn{1}{l|}{SVI}                                 & \begin{tabular}[c]{@{}l@{}}Prop. Budget Expenditure~\eqref{const: j40 bud}\end{tabular}  & 1000.0 & \eqref{obj: mm_pct_w_j40_bud-PSPS}   \\ \hline
\textbf{E-M10}                                                            & \(\max\) \% load shed across groups \eqref{const:alph_m}                                                    & \multicolumn{1}{l|}{SVI}                                 & \begin{tabular}[c]{@{}l@{}}Proportional Load Shed Reduction~\eqref{const: j40 ls}\end{tabular} & 1000.0  & \eqref{obj: mm_pct_w_j40_ls-PSPS}  \\ \hline
\end{tabular}\caption{Summary of all models considered in the results section of this paper. Appendix results show alternative budget outcomes. All models include the constraints \eqref{eq: dcots}, and~\eqref{const: Line high}--\eqref{const: undergrounding threshold}. Budgets are listed in millions USD. }\label{table: model summary}
\end{sidewaystable*}

\subsection{Baseline Load Shed Under No Intervention}\label{subsection: ls under no intervention}

First, we discuss the results of the \textit{no-budget} baseline model (BL-M0), which shows the regions and populations most subjected to load shed from \gls{psps} events when there is no budget for undergrounding power lines. 
Figure \ref{fig: svi plots with load shedding} shows the load shed patterns from a PSPS event on the network covered by the Texas7k case study for the 5-day period from June 11\textsuperscript{th} through June 15\textsuperscript{th} of 2021\footnote{We choose this week because of the relatively high load and wildfire risk on the network during this week.}.  Figure \ref{fig: svi plots with load shedding} highlights two regions in the north and west of Texas that bear most of the load shed in the network overlaid on choropleths of the vulnerable census tracts in Texas per the \mbox{\gls{cejst}} (Figure \ref{fig:subfig:j40_ls}) and \gls{svi} (Figure \ref{fig:subfig:svi_ls}).

\begin{figure}[ht]
    \centering
    \includegraphics[width=0.7\linewidth]{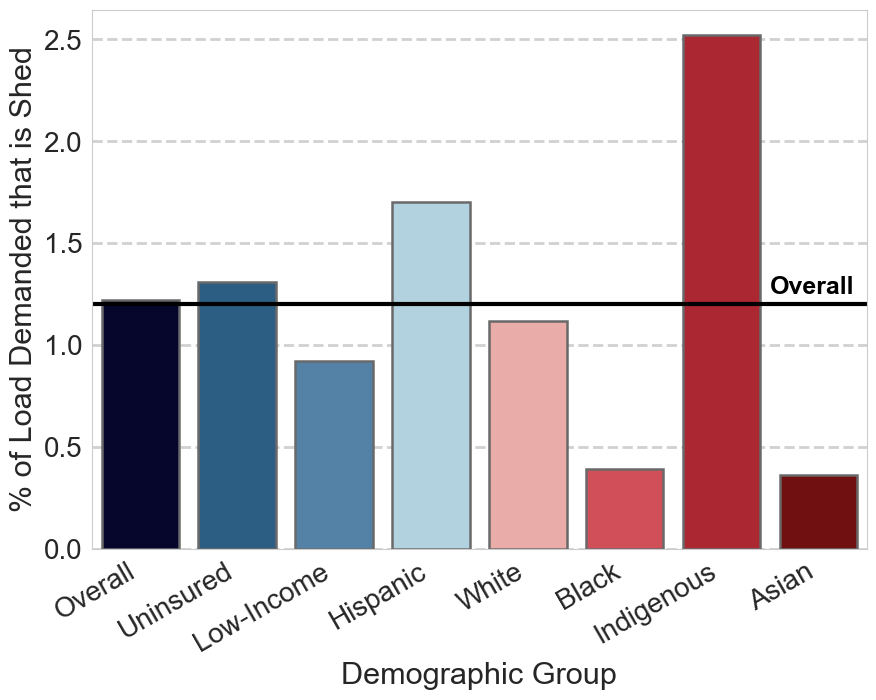}
    \caption{In the \textit{no-budget} baseline case (BL-M0) we observe that Indigenous, Hispanic, and uninsured populations face above-average percentages of load shed.}
    \label{fig:baseline load shed}
\end{figure}

Figure \ref{fig:baseline load shed} presents the load shed results for the no-budget baseline model (BL-M0), which show that Indigenous, Hispanic, and uninsured populations face above-average percentages of load shed. In particular, Indigenous populations experience \textcolor{red}{more than \sout{nearly}} twice as much load shed as the average Texas resident (2.52\% vs 1.22\%).
Figure \ref{fig:subfig:indig_ls} shows that census tracts with relatively larger (although still minority) Indigenous populations tend to live in areas experiencing load shed. 

One might be surprised to see that the low-income group in Figure \ref{fig:baseline load shed} experiences a lower percentage of load shed than the overall percentage of load shed; however, this finding is not an anomaly.
\citet{shah2023inequitable} studied the correlative effects of race and income-level on the likelihood that a Texas resident experienced a power outage during 2021 Winter Storm Uri, and the authors found that ethnic minority status has greater impact on load shed than income group. \citet{shah2023inequitable} also found that at every income level, minority groups experienced disproportionately high amounts of power loss compared to non-minority groups. Our own simulations of power loss by income group, a sample of which is shown in Appendix \ref{appendix:income and ls}, reflect the same trend as \citet{shah2023inequitable} that income has very little effect on the likelihood of experiencing power loss. Hence, in this paper, we primarily discuss discrepancies in load shed by racial group.

\subsection{Load Shed Under Increasing Budgets for Undergrounding}\label{subsec: ls increased budgets}
\begin{figure*}[htbp]
\centering
\begin{subfigure}[t]{.325\textwidth}
  \centering
    \includegraphics[width=.99\linewidth]{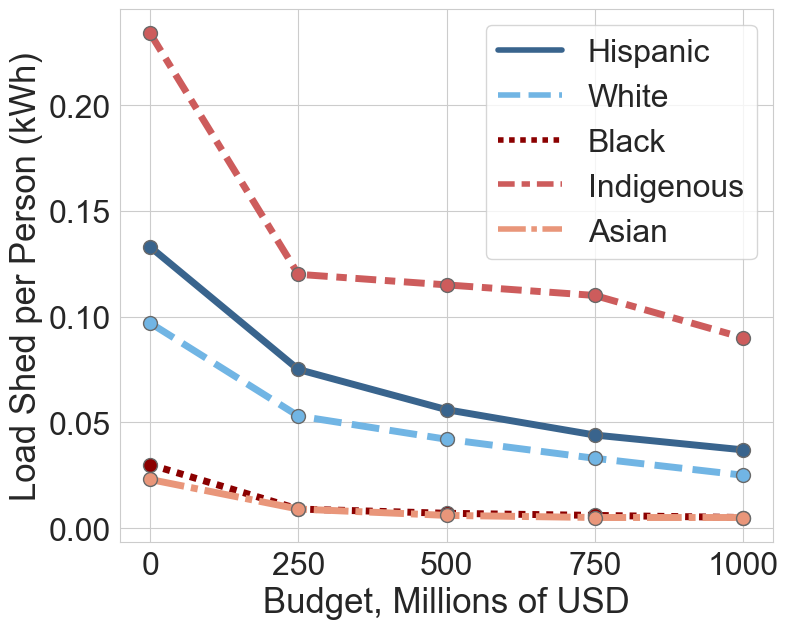} 
    \caption{Normalized load shed by group}
  \label{fig: subfig: norm ls across budgets}
\end{subfigure}
\begin{subfigure}[t]{.325\textwidth}
  \centering
  \includegraphics[width=.99\linewidth]{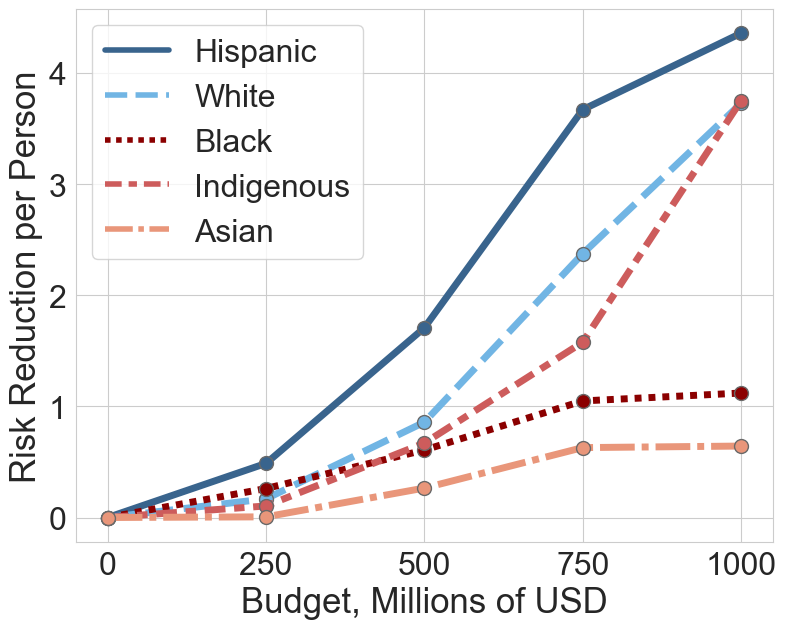}  
  \caption{Normalized risk reduction by group}
  \label{fig: subfig: fig: subfig: norm risk redux across budgets}
\end{subfigure}
\begin{subfigure}[t]{.325\textwidth}
  \centering
  \includegraphics[width=.99\linewidth]{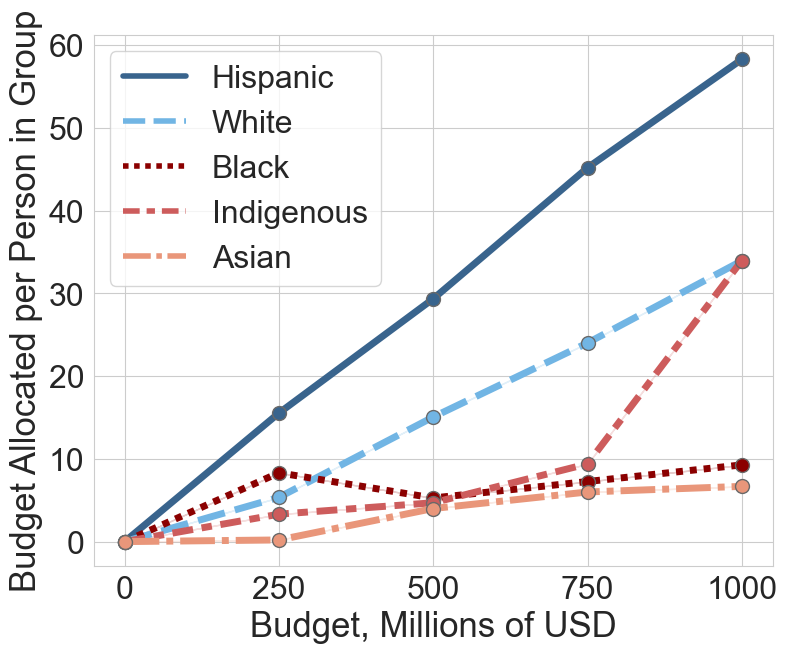}  
  \caption{Normalized budget allocation (in USD) per person in group}
  \label{fig: subfig: norm cost allocated across budgets}
\end{subfigure}
\caption{Trends in load shed, risk reduction, and budget allocation per person in each group  under BL-M1 as a function of the total budget allocated for undergrounding.}
\label{fig: no_equity_summaries}
\end{figure*}
Now, we consider how load shed patterns change when we allocate budget for power line undergrounding according to the \textit{baseline} model \eqref{obj: base-PSPS}. Ideally, if budget allocations were equitable, the groups that were subject to the most load shed in the \textit{no-budget} case (BL-M0) should see relatively higher percentages of budget allocated to them, as well as higher decreases in load shed and wildfire risk as budget increases. For this case study, we would expect to see this for Indigenous and Hispanic populations.

In Figure~\ref{fig: no_equity_summaries}, we show how the per capita load shed, risk reduction (as compared to the no-budget case), and budget-allocated by group changes as the total available budget increases. We observe that Indigenous and Hispanic populations experience the highest load shed on a per capita basis (Figure \ref{fig: subfig: norm ls across budgets}), across all budgets. However, as the amount of budget allocated increases, Hispanic populations see steady decreases in load shed per person and the highest per capita increases in risk reduction and budget allocation. In contrast, Indigenous populations have a large initial drop in per capita load shed when \$250 million is allocated, but  do not experience substantial reductions in load shed per person until the \$1 billion budget is allocated. 
Furthermore, Indigenous groups neither experience substantial reductions in wildfire risk (Figure \ref{fig: subfig: fig: subfig: norm risk redux across budgets}) nor commensurate levels of budget allocation (Figure \ref{fig: subfig: norm cost allocated across budgets}) relative to their risk and load shed profile until \$1 billion is allocated. This indicates that Indigenous groups live in census tracts that are either (a) too expensive for power line undergrounding or (b) reside in census tracts where they constitute a minority of the population---and thus are not prioritized by the baseline model---or both (a) and (b). These hypotheses are corroborated by Figures \ref{fig: subfig: cost to underground}  and \ref{fig:subfig:indig_ls}, which show that the census tracts with the highest portion of Indigenous populations still only see less than 5\% of their population comprised of Indigenous populations. These same census tracts (in Region 2 of the figures) have many transmission lines that cost upwards of 100 million USD to underground. Further discussion about the challenge of using vulnerability indices when optimizing for equitable resource allocation can be found in \S \ref{sec: curse of agg}.

\subsection{Impacts of Justice40-Style Resource Allocation Using Policy Constraints}\label{sec: quant results}
The previous section showed that Indigenous communities in Texas face disproportionately high load shed under the baseline model, even as increased budget is allocated to reduce load shed. 
Policies like the Justice40 initiative seek to ensure some minimal level (40\%) of realized benefit (e.g., relief from load shed or budget allotment) is allocated to vulnerable populations. In this section, we analyze whether constraints modeling the Justice40 initiative rectify the current imbalance in load shed experienced across demographic groups after investments have been made.  
In other words, {\it do the {\sc Policy} constraints result in the intended effect?} 

\begin{figure*}[h]
 \centering
 \begin{subfigure}[t]{.47\textwidth}
        \centering
        \captionsetup{width=.9\linewidth}
        \includegraphics[width=\linewidth]{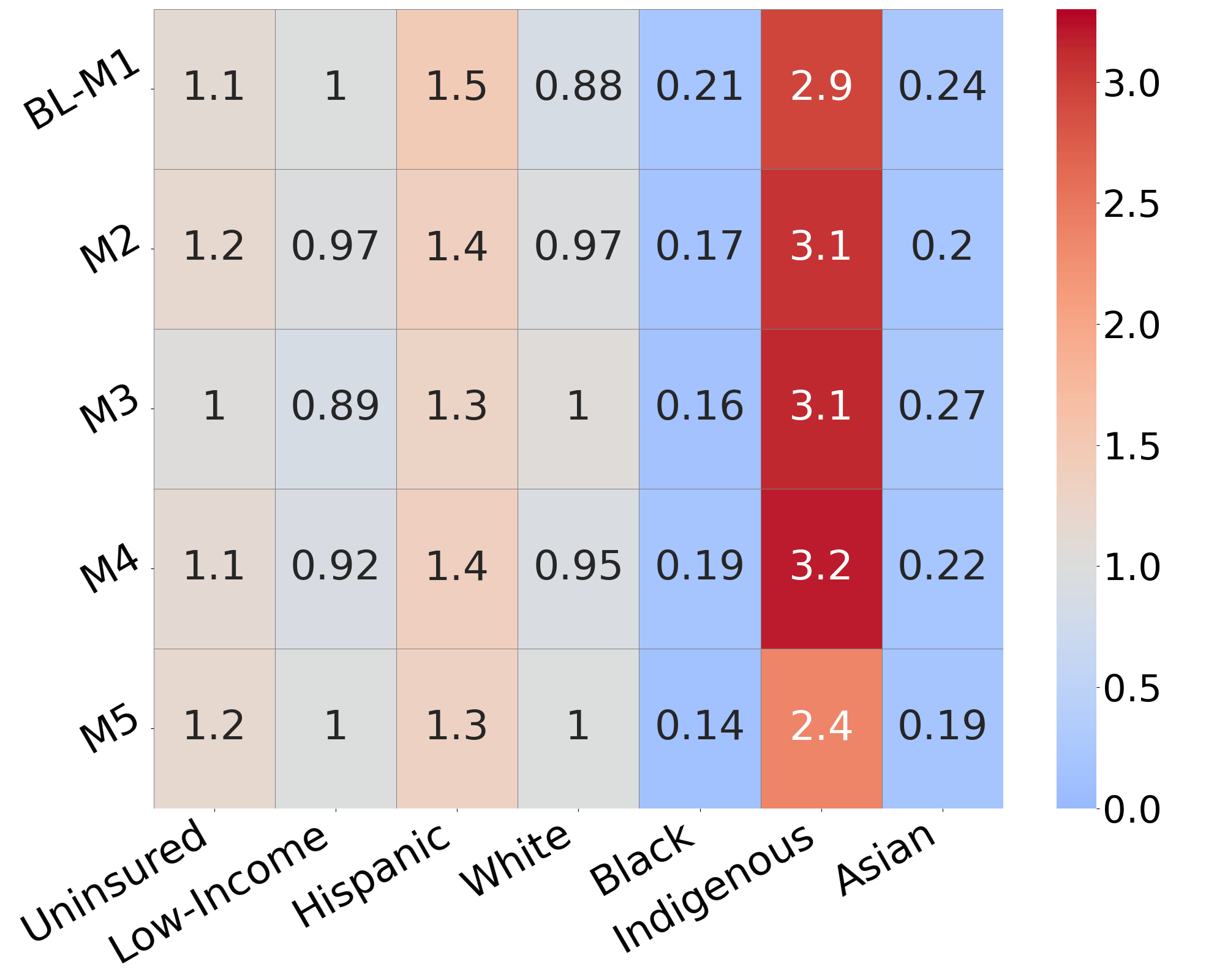}
        \caption{Relative unfairness in the percent of load demanded that is shed at the \$1 billion budget for simulations NOT using the {\sc Equity} objective.}
        \label{fig: unfairness B1000 - no equity}
    \end{subfigure}
    \hfill 
    \begin{subfigure}[t]{.47\textwidth}
        \centering
        \captionsetup{width=.9\linewidth} 
        \includegraphics[width=\linewidth]{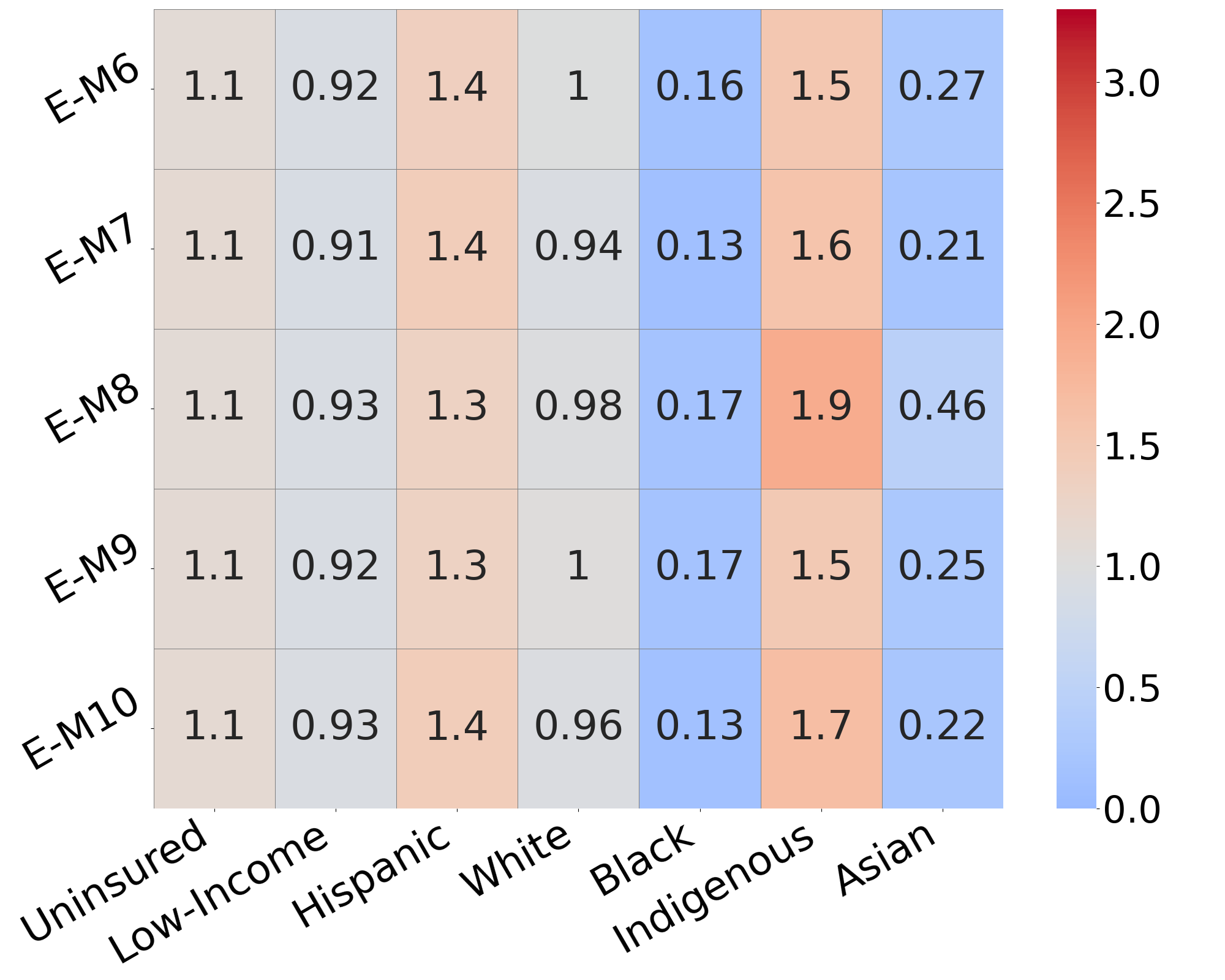}
        \caption{Relative unfairness in the percent of load demanded that is shed at the \$1 billion budget for simulations using the {\sc Equity} objective.}
        \label{fig: unfairness B1000 - equity}
    \end{subfigure}
    \caption{Relative unfairness is the ratio of the percent of a group's load shed to the percent of load shed over individuals in the entire network. The larger the relative unfairness score, the more unfair the load shed outcomes are for that group. We can see that the maximum unfairness values decrease when considering the MMF framework with the {\sc Equity} objective (b) versus the {\sc Policy} constraints alone (a).}\label{fig: unfairness-ls}
\end{figure*}

Optimization models \eqref{obj: base-j40_bud} and \eqref{obj: base-j40_ls} employ {\sc Policy} constraints by budget allocation \eqref{const: j40 bud} and load shed reduction \eqref{const: j40 ls}, respectively. 
Table \ref{table: results - no mmf} shows the results for these models for a \$1 billion budget. Recall that model BL-M0 in the first row represents the baseline no-budget case, and is given as a reference. Each row in the table represents one model, which is defined by a vulnerability index, constraint type pair, where the constraint type indicates whether we implement a {\sc Policy} constraint by budget, reduction in load shed over the no-budget model, or have no additional constraint. The vulnerability index column indicates which index is used to designate the vulnerable communities when a {\sc Policy} constraint is in effect. Each cell gives the percent of a given group's load demanded that is shed. We note that load shed values above 1\% have red text to indicate a rate of load shed that we highlight as particularly elevated. From these load shed values, we can compute an ``unfairness ratio," which is the ratio of that group's percent of load shed to the overall population's percent of load shed. The unfairness ratios for models with {\sc Policy} constraints and the baseline ``minimize total load shed" objective are shown in Figure \ref{fig: unfairness B1000 - no equity}.

\begin{sidewaystable}[!htbp]
\small
\centering
\begin{tabular}{|l|ll|llllllll|}
\hline
\rowcolor[HTML]{EFEFEF} 
\multicolumn{1}{|c|}{\cellcolor[HTML]{EFEFEF}}                            & \multicolumn{2}{c|}{\cellcolor[HTML]{EFEFEF}\textbf{Policy Constraint}}                                                                                                                                                                 & \multicolumn{8}{c|}{\cellcolor[HTML]{EFEFEF}\textbf{Percent of Load Demanded that is Shed}}                                                                                                                               \\ \cline{2-11} 
\rowcolor[HTML]{EFEFEF} 
\multicolumn{1}{|c|}{\multirow{-2}{*}{\cellcolor[HTML]{EFEFEF}\textbf{}}} & \multicolumn{1}{c|}{\cellcolor[HTML]{EFEFEF}\textbf{\begin{tabular}[c]{@{}c@{}}Vulnerability \\ Index\end{tabular}}} & \multicolumn{1}{c|}{\cellcolor[HTML]{EFEFEF}\textbf{\begin{tabular}[c]{@{}c@{}}Constraint \\ Type\end{tabular}}} & \multicolumn{1}{c|}{\cellcolor[HTML]{EFEFEF}\textbf{Overall}}                     & \multicolumn{1}{c|}{\cellcolor[HTML]{EFEFEF}\textbf{Uninsured}}                   & \multicolumn{1}{c|}{\cellcolor[HTML]{EFEFEF}\textbf{Low-Income}} & \multicolumn{1}{c|}{\cellcolor[HTML]{EFEFEF}\textbf{Hispanic}}                   & \multicolumn{1}{c|}{\cellcolor[HTML]{EFEFEF}\textbf{White}}                       & \multicolumn{1}{c|}{\cellcolor[HTML]{EFEFEF}\textbf{Black}} & \multicolumn{1}{c|}{\cellcolor[HTML]{EFEFEF}\textbf{Indigenous}}                  & \multicolumn{1}{c|}{\cellcolor[HTML]{EFEFEF}\textbf{Asian}} \\ \hline
\rowcolor[HTML]{ECF4FF} 
\cellcolor[HTML]{DAE8FC}\textbf{BL-M0}                                    & \multicolumn{1}{l|}{\cellcolor[HTML]{DAE8FC}None}                                                                    & \cellcolor[HTML]{DAE8FC}None                                                                                     & \multicolumn{1}{l|}{\cellcolor[HTML]{ECF4FF}{\color[HTML]{FE0000} \textbf{1.22}}} & \multicolumn{1}{l|}{\cellcolor[HTML]{ECF4FF}{\color[HTML]{FE0000} \textbf{1.31}}} & \multicolumn{1}{l|}{\cellcolor[HTML]{ECF4FF}0.92}                & \multicolumn{1}{l|}{\cellcolor[HTML]{ECF4FF}{\color[HTML]{FE0000} \textbf{1.7}}} & \multicolumn{1}{l|}{\cellcolor[HTML]{ECF4FF}{\color[HTML]{FE0000} \textbf{1.12}}} & \multicolumn{1}{l|}{\cellcolor[HTML]{ECF4FF}0.39}           & \multicolumn{1}{l|}{\cellcolor[HTML]{ECF4FF}{\color[HTML]{FE0000} \textbf{2.52}}} & 0.36                                                        \\ \hline
\cellcolor[HTML]{EFEFEF}\textbf{BL-M1}                                    & \multicolumn{1}{l|}{\cellcolor[HTML]{EFEFEF}None}                                                                    & \cellcolor[HTML]{EFEFEF}None                                                                                     & \multicolumn{1}{l|}{0.33}                                                         & \multicolumn{1}{l|}{0.37}                                                         & \multicolumn{1}{l|}{0.33}                                        & \multicolumn{1}{l|}{0.48}                                                        & \multicolumn{1}{l|}{0.29}                                                         & \multicolumn{1}{l|}{0.07}                                   & \multicolumn{1}{l|}{0.97}                                                         & 0.08                                                        \\ \hline
\cellcolor[HTML]{EFEFEF}\textbf{M2}                                       & \multicolumn{1}{l|}{\cellcolor[HTML]{EFEFEF}CEJST}                                                                   & \cellcolor[HTML]{EFEFEF}Budget                                                                                   & \multicolumn{1}{l|}{0.35}                                                         & \multicolumn{1}{l|}{0.41}                                                         & \multicolumn{1}{l|}{0.34}                                        & \multicolumn{1}{l|}{0.48}                                                        & \multicolumn{1}{l|}{0.34}                                                         & \multicolumn{1}{l|}{0.06}                                   & \multicolumn{1}{l|}{{\color[HTML]{FE0000} \textbf{1.07}}}                         & 0.07                                                        \\ \hline
\cellcolor[HTML]{EFEFEF}\textbf{M3}                                       & \multicolumn{1}{l|}{\cellcolor[HTML]{EFEFEF}CEJST}                                                                   & \cellcolor[HTML]{EFEFEF}Load Shed                                                                                & \multicolumn{1}{l|}{0.45}                                                         & \multicolumn{1}{l|}{0.46}                                                         & \multicolumn{1}{l|}{0.40}                                        & \multicolumn{1}{l|}{0.58}                                                        & \multicolumn{1}{l|}{0.47}                                                         & \multicolumn{1}{l|}{0.07}                                   & \multicolumn{1}{l|}{{\color[HTML]{FE0000} \textbf{1.41}}}                         & 0.12                                                        \\ \hline
\cellcolor[HTML]{EFEFEF}\textbf{M4}                                       & \multicolumn{1}{l|}{\cellcolor[HTML]{EFEFEF}SVI}                                                                     & \cellcolor[HTML]{EFEFEF}Budget                                                                                   & \multicolumn{1}{l|}{0.37}                                                         & \multicolumn{1}{l|}{0.42}                                                         & \multicolumn{1}{l|}{0.34}                                        & \multicolumn{1}{l|}{0.51}                                                        & \multicolumn{1}{l|}{0.35}                                                         & \multicolumn{1}{l|}{0.07}                                   & \multicolumn{1}{l|}{{\color[HTML]{FE0000} \textbf{1.18}}}                         & 0.08                                                        \\ \hline
\cellcolor[HTML]{EFEFEF}\textbf{M5}                                       & \multicolumn{1}{l|}{\cellcolor[HTML]{EFEFEF}SVI}                                                                     & \cellcolor[HTML]{EFEFEF}Load Shed                                                                                & \multicolumn{1}{l|}{0.37}                                                         & \multicolumn{1}{l|}{0.43}                                                         & \multicolumn{1}{l|}{0.37}                                        & \multicolumn{1}{l|}{0.49}                                                        & \multicolumn{1}{l|}{0.37}                                                         & \multicolumn{1}{l|}{0.05}                                   & \multicolumn{1}{l|}{0.88}                                                         & 0.07                                                        \\ \hline
\end{tabular}
\caption{Percentage load shed across different groups, when not using the {\sc Equity} objective and different types of {\sc Policy} constraints when a \$1 billion budget is allocated. The blue-highlighted row showing the results of model BL-M0 represents the no-budget baseline case (i.e., no undergrounding) for reference. The red cells in the table indicate that the \% load shed is above a threshold of 1\%.}
\label{table: results - no mmf}
\end{sidewaystable}

\begin{sidewaystable}[!htbp]
\small
\centering
\begin{tabular}{|
>{\columncolor[HTML]{EFEFEF}}l |
>{\columncolor[HTML]{EFEFEF}}l 
>{\columncolor[HTML]{EFEFEF}}l |cccccccc|}
\hline
\multicolumn{1}{|c|}{\cellcolor[HTML]{EFEFEF}\textbf{}} & \multicolumn{2}{c|}{\cellcolor[HTML]{EFEFEF}\textbf{Policy Constraint}}                                                                  & \multicolumn{8}{c|}{\cellcolor[HTML]{EFEFEF}\textbf{Percent of Load Demanded that is Shed}}                                                                \\ \hline
\multicolumn{1}{|c|}{\cellcolor[HTML]{EFEFEF}\textbf{}} & \multicolumn{1}{c|}{\cellcolor[HTML]{EFEFEF}\textbf{\begin{tabular}[c]{@{}c@{}}Vulnerability \\ Index\end{tabular}}} & \multicolumn{1}{c|}{\cellcolor[HTML]{EFEFEF}\textbf{\begin{tabular}[c]{@{}c@{}}Constraint \\ Type\end{tabular}}} & \multicolumn{1}{c|}{\cellcolor[HTML]{EFEFEF}\textbf{Overall}} & \multicolumn{1}{c|}{\cellcolor[HTML]{EFEFEF}\textbf{Uninsured}} & \multicolumn{1}{c|}{\cellcolor[HTML]{EFEFEF}\textbf{Low-Income}} & \multicolumn{1}{c|}{\cellcolor[HTML]{EFEFEF}\textbf{Hispanic}} & \multicolumn{1}{c|}{\cellcolor[HTML]{EFEFEF}\textbf{White}} & \multicolumn{1}{c|}{\cellcolor[HTML]{EFEFEF}\textbf{Black}} & \multicolumn{1}{c|}{\cellcolor[HTML]{EFEFEF}\textbf{Indigenous}} & \cellcolor[HTML]{EFEFEF}\textbf{Asian} \\ \hline
\textbf{E-M6}                                           & \multicolumn{1}{l|}{\cellcolor[HTML]{EFEFEF}None}                                                                    & None                                                                                                             & \multicolumn{1}{c|}{0.49}                                     & \multicolumn{1}{c|}{0.54}                                       & \multicolumn{1}{c|}{0.45}                                        & \multicolumn{1}{c|}{0.67}                                      & \multicolumn{1}{c|}{0.49}                                   & \multicolumn{1}{c|}{0.08}                                   & \multicolumn{1}{c|}{0.75}                                        & 0.13                                   \\ \hline
\textbf{E-M7}                                           & \multicolumn{1}{l|}{\cellcolor[HTML]{EFEFEF}CEJST}                                                                   & Budget                                                                                                           & \multicolumn{1}{c|}{0.47}                                     & \multicolumn{1}{c|}{0.53}                                       & \multicolumn{1}{c|}{0.43}                                        & \multicolumn{1}{c|}{0.67}                                      & \multicolumn{1}{c|}{0.44}                                   & \multicolumn{1}{c|}{0.06}                                   & \multicolumn{1}{c|}{0.75}                                        & 0.10                                   \\ \hline
\textbf{E-M8}                                           & \multicolumn{1}{l|}{\cellcolor[HTML]{EFEFEF}CEJST}                                                                   & Load Shed                                                                                                        & \multicolumn{1}{c|}{0.41}                                     & \multicolumn{1}{c|}{0.45}                                       & \multicolumn{1}{c|}{0.38}                                        & \multicolumn{1}{c|}{0.54}                                      & \multicolumn{1}{c|}{0.40}                                   & \multicolumn{1}{c|}{0.07}                                   & \multicolumn{1}{c|}{0.79}                                        & 0.19                                   \\ \hline
\textbf{E-M9}                                           & \multicolumn{1}{l|}{\cellcolor[HTML]{EFEFEF}SVI}                                                                     & Budget                                                                                                           & \multicolumn{1}{c|}{0.53}                                     & \multicolumn{1}{c|}{0.59}                                       & \multicolumn{1}{c|}{0.49}                                        & \multicolumn{1}{c|}{0.70}                                      & \multicolumn{1}{c|}{0.54}                                   & \multicolumn{1}{c|}{0.09}                                   & \multicolumn{1}{c|}{0.79}                                        & 0.13                                   \\ \hline
\textbf{E-M10}                                          & \multicolumn{1}{l|}{\cellcolor[HTML]{EFEFEF}SVI}                                                                     & Load Shed                                                                                                        & \multicolumn{1}{c|}{0.45}                                     & \multicolumn{1}{c|}{0.51}                                       & \multicolumn{1}{c|}{0.42}                                        & \multicolumn{1}{c|}{0.64}                                      & \multicolumn{1}{c|}{0.43}                                   & \multicolumn{1}{c|}{0.06}                                   & \multicolumn{1}{c|}{0.75}                                        & 0.10                                   \\ \hline
\end{tabular}
\caption{Percentage load shed across different groups, for all combinations of {\sc Equity} objective (group-level protections) and {\sc Policy} constraints with a \$1 billion budget.}
\label{table: results - mmf}
\end{sidewaystable}

No-budget model BL-M0 and baseline model for \$1 billion budget BL-M1 in Table \ref{table: results - no mmf} seek to minimize the total load shed in the network subject to power flow and budget constraints only (i.e., there are no {\sc Policy} constraints). 
Between BL-M0 and BL-M1, there is a substantial reduction in load shed across every group when budget is allocated for undergrounding, although the level of benefit is not uniform. For example, while the overall population sees their load shed decrease by a factor of nearly four, Indigenous load shed does not even decrease by a factor of three. 

Now, consider models M1-M5. Under BL-M1, Indigenous populations experience 0.97\% of load shed, which is almost three times higher than that of the overall population at 0.33\%.\footnote{We can compare these percentages with the System Average Interruption Duration Index (SAIDI) which measures the average number of minutes of power loss seen by a costumer in a year for different states. When looking at Texas in 2021, we see roughly 1500 minutes of power loss or 0.2\% of the year. However, Winter Storm Uri impacted Texas during 2021, causing widespread outages. When looking at the 10 year average SAIDI for Texas from 2013 to 2022, we see a value of 0.074\% of the year \cite{saidi2022table}.} For models M2-M5, we test various optimization models that minimize the total load shed in the network subject to {\sc Policy} constraints ensuring that 40\% of the benefit goes to vulnerable tracts, designated by either CEJST or SVI. For each of these, with the exception of M5, Table \ref{table: results - no mmf} shows a consistently higher percent of load shed for Indigenous populations compared to BL-M1, thus, not having a considerable impact on protecting this vulnerable population. Even for model M5, the Indigenous load shed reduction is minimal, with an improvement of only 0.09\% over the baseline. We also observe that, in general, no group is substantially ``better off" (in terms of percent of load shed) after implementing {\sc Policy} constraints---in fact, they are usually \textit{worse} off.

\begin{figure*}[ht]
     \centering
     \begin{subfigure}{.99\textwidth}
        \centering
        \captionsetup{width=.95\linewidth}
        \includegraphics[width=.97\linewidth]{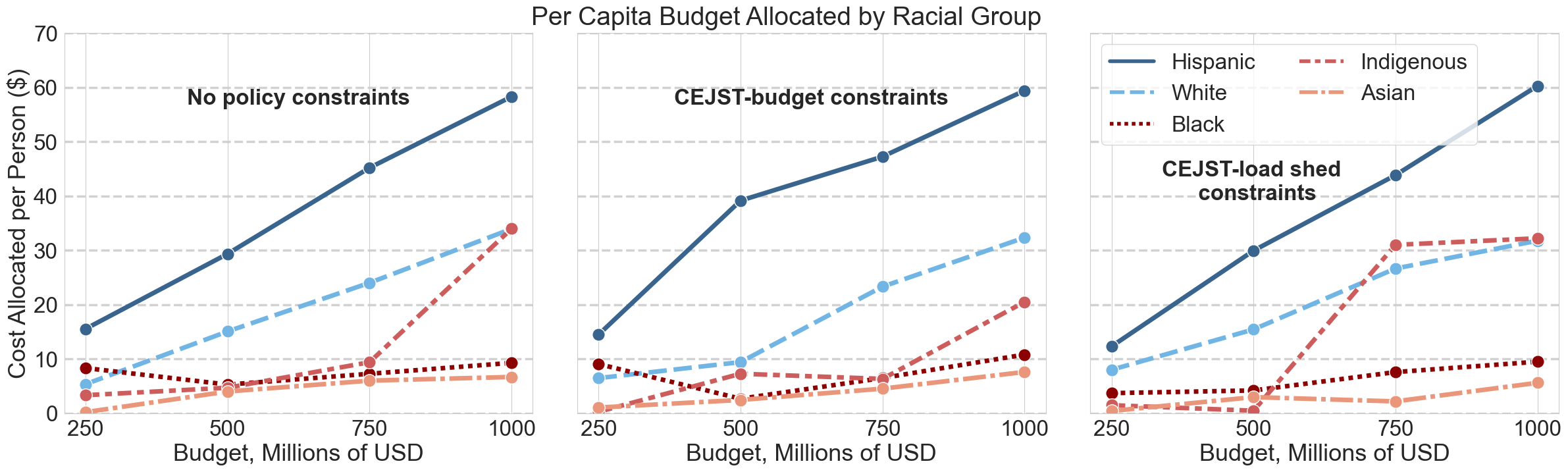}
         \caption{Normalized cost allocated to group in the baseline case (no {\sc policy} constraints) (left), the case when constraining 40\% of the budget to go to Justice40 communities (center), and the case when constraining 40\% of load shed reduction to go to Justice40 communities (right).}
         \label{fig: subfig: cost_allocated}
     \end{subfigure}
     \vspace{2mm}
     \\
     \begin{subfigure}{.99\textwidth}
        \centering
        \captionsetup{width=.95\linewidth}
        \includegraphics[width=.97\linewidth]{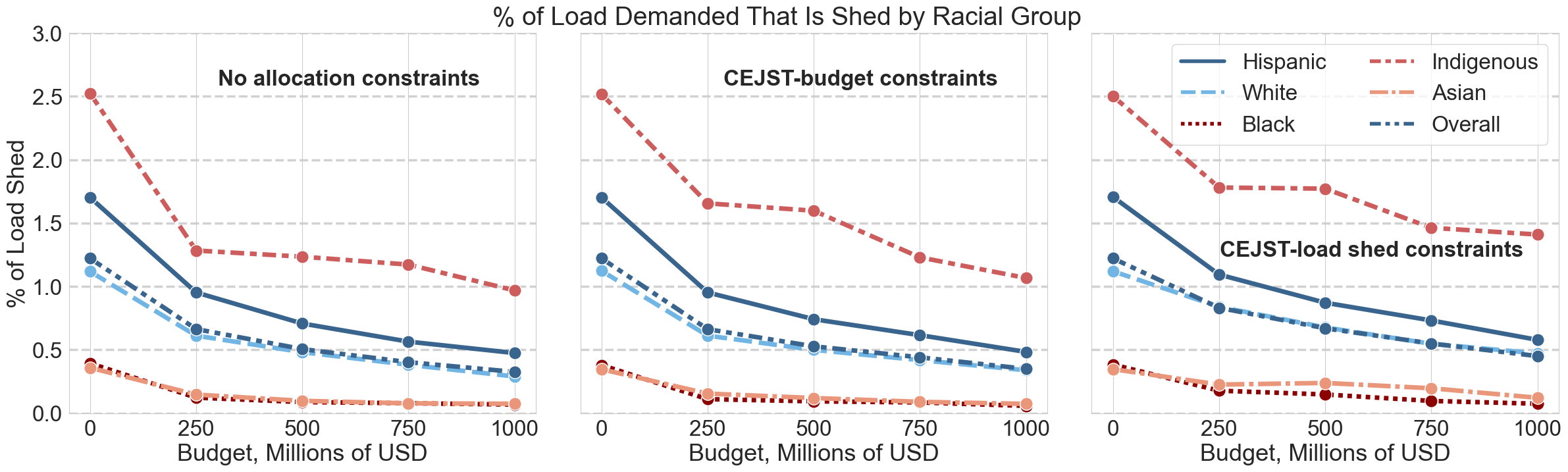}
         \caption{Percent of load demanded that is shed by racial group in the baseline case (no {\sc policy} constraints) (left), the case when constraining 40\% of the budget to go to Justice40 communities (center), and the case when constraining 40\% of load shed reduction to go to Justice40 communities (right).}
         \label{fig: subfig: ls}
     \end{subfigure}
     \caption{The budget allocation increase to Indigenous populations (top right) does not translate to load shed reduction (bottom right). In fact, the bottom subfigure shows load shed trends by racial group remain relatively consistent across each of these policy-constrained cases.}
     \label{fig: unfairness-budget}
 \end{figure*}

Furthermore, Figure \ref{fig: subfig: ls} shows how every racial group---including Indigenous groups---experiences \textit{more} load shed after either Justice40-style {\sc Policy} constraint is implemented, even though the budget allocation \textit{does} become increasingly allocated to Indigenous groups. 
It is likely that the discrepancy between budget allocation and load shed reduction stems from {\sc Policy} constraints channeling funds towards census tracts with lower population density or tracts where the cost to underground power lines is higher, reducing the total amount of lines that can be undergrounded. Consequently, although the per capita investment for Indigenous groups may be more substantial, the per capita reduction in load shed remains unimproved. 

\subsection{Group Fair Frameworks Using Equity Objective}
\label{subsec: group fair}
In this case study, it is clear that protecting Indigenous populations from load shed is challenging due to their minority status in each census tract and their residence in areas where undergrounding is expensive and wildfire ignition risk is high.\footnote{Note that low income communities may be more likely to reside in high wildfire risk areas given lower housing prices in wildfire-prone areas \citep{mueller2009repeated}.} 
This motivates the study of \textit{group-level protections}, which may help ensure that each population group of interest is able to receive the appropriate level of resources. This framework provides a lower bound on the percentage load shed experienced by the proportionally most-affected group, which in our study is Indigenous populations.
\begin{figure}[h]
    \centering
    \includegraphics[width=0.75\linewidth]{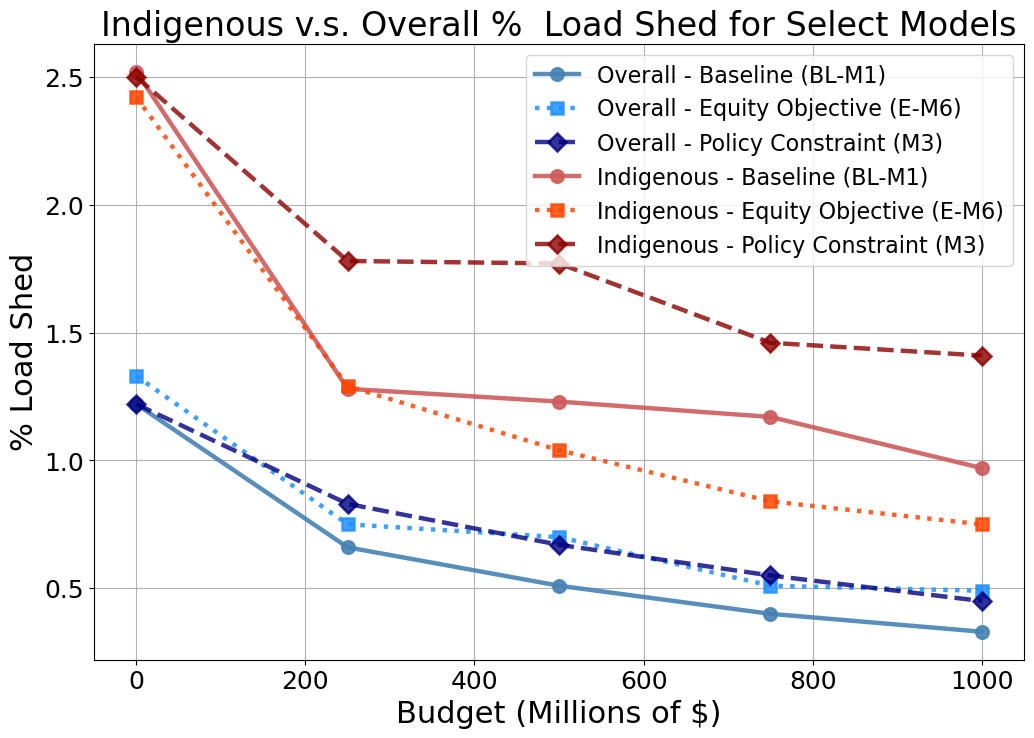}
    \caption{{\sc Policy} constraints increase the percent of load shed experienced by Indigenous households over BL-M1, whereas the {\sc Equity} objective decreases the percent of load shed experienced by Indigenous households, particularly at higher budgets. The overall percentage of load shed increases relatively little when implementing either {\sc Policy} constraints or the {\sc Equity} objective.}
    \label{fig:revision-plot}
\end{figure}
 
 Optimization models \eqref{obj: mm_pct}, \eqref{obj: mm_pct_w_j40_bud-PSPS}, and \eqref{obj: mm_pct_w_j40_ls-PSPS} employ group-level protections by implementing a percentage-based \gls*{mmf} framework using the {\sc Equity} objective defined by \eqref{const:alph_m}. 
 By using this percentage-based \gls*{mmf} framework, we account for the total load demanded by each group, which prevents minority groups within census tracks from being de-prioritized for investment allocation. We summarize the results of using any combination of {\sc Equity} objective and  {\sc Policy} constraint for a \$1 billion budget in Table \ref{table: results - mmf}. \textcolor{red}{We also provide a visual comparison of load shed statistics across variance budgets between Indigenous groups and the overall Texas population in Figure \ref{fig:revision-plot}.} Note that the ``E" prefix in the model names E-M6 through E-M10 correspond to the ``{\sc Equity}'' objective.
 
\textcolor{red}{First, we discuss the overall trends from the 1 billion USD budget case in Table \ref{table: results - mmf}}. We observe from models E-M6 through E-M10 that using {\sc Equity} objectives with and without {\sc Policy} constraints reduces the percentage of load shed across {\it all} populations to be under 0.80\%. In particular, between the baseline BL-M1 and the models using the {\sc Equity} objective, we see a minimum of a 19\% decrease in percent of load demanded that is shed. However, prioritizing relief to Indigenous communities comes at a cost; the average load shed across all individuals increases from 0.33\% to, at best, 0.41\%, and all other non-Indigenous groups see increases in percent of load shed. One could make the case that the increases in load shed are relatively minor since each non-Indigenous group experiences less than 0.7\% of their load shed under the {\sc Equity} objective, which is a relatively small amount of power loss on the network. However, it is ultimately the task of a policymaker to determine what level of overall load shed increase is tolerable to reduce disparities across groups.

In Figure \ref{fig: unfairness B1000 - equity}, we show the unfairness ratios when incorporating the {\sc Equity} objective. While Indigenous and Hispanic groups still face disproportionately high levels of load shed relative to the overall population, these unfairness ratios are significantly lower than the unfairness ratios when not using the {\sc Equity} objective, shown in Figure \ref{fig: unfairness B1000 - no equity}. When we consider both the overall percent of load shed in Table \ref{table: results - mmf} 
and the unfairness ratios in Figure \ref{fig: unfairness B1000 - equity}, we see that all of the models show varying degrees of trade-off between the overall load shed and the unfairness ratio experienced by Indigenous groups.

In terms of providing power loss protections to Indigenous populations, we show a significant benefit of using {\sc Equity} objectives in conjunction with {\sc Policy} constraints. While Tables \ref{table: results - no mmf} and \ref{table: results - mmf} only show the results for the \$1 billion budget, we note that the \gls*{mmf} framework, which promotes the most load shed relief for Indigenous groups under the \$1 billion budget, also promotes load shed relief compared to the baseline at the \$500 million and \$750 million budgets. \textcolor{red}{The relative improvement in Indigenous load shed outcomes can be seen more directly in Figure \ref{fig:revision-plot} for a subset of the models. For tabular results like those in Tables \ref{table: results - no mmf} amd \ref{table: results - mmf} for the  \$500 million case, see Figure~\ref{fig: pct loadshed} in the Appendix.} 
Hence, we find that in order to see meaningful reductions in Indigenous load shed while maintaining reasonably low load shed for other groups, our study identifies two requirements: (1)~a sufficiently high budget, in our case, at least \$500 million \sout{(See Appendix \ref{appendix: extended results})}, and (2)~a \gls*{mmf} framework that considers minimizing the maximum percentage of a group's load demanded that is shed. This latter requirement is necessary to place Indigenous groups on equal priority with different racial, ethnic, and other groups which make up a higher percentage of Texas' total population.

\subsection{The Curse of Aggregation in Vulnerability Index Creation}\label{sec: curse of agg}
Our analysis reveals that adding {\sc Policy} constraints generally leads to increased load shed outcomes for \textit{all} racial and socioeconomic groups, and does not alleviate the disproportionately high load shed experienced by Indigenous Texans. This is important because Indigenous populations experience nearly double the poverty rate of Texas overall, and would be ideally flagged as disadvantaged by a vulnerability index. Yet, we do not observe any benefit from these {\sc Policy} constraints. Why?

\begin{figure}[t!]
    \centering
    \includegraphics[width=0.7\linewidth]{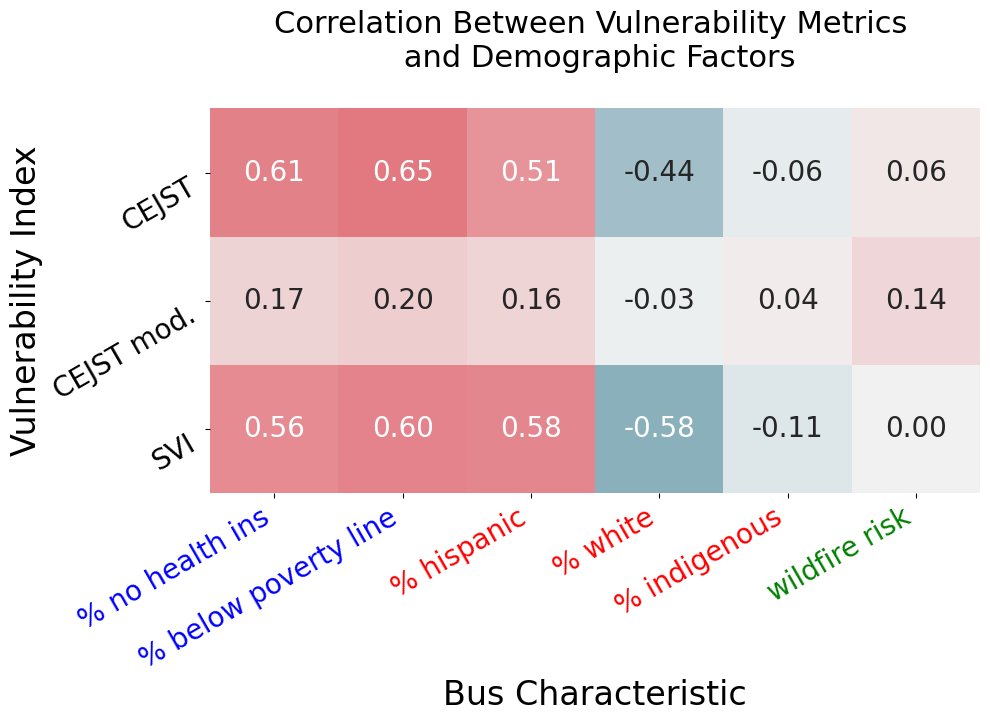}
    \caption{Vulnerability metrics correlate with socioeconomic indicators (blue text), but do not necessarily correlate with racial/ethnicity variables (red text) or indicators of high wildfire ignition risk (green text). Blue cells indicate negative correlation values while red cells indicate positive values.}
    \label{fig:corrplot}
\end{figure}

The challenge with using any vulnerability index, including the \mbox{\gls{cejst}} and \gls{svi} metrics, lies in the necessity of data aggregation, which homogenizes populations within a census tract and can obscure the vulnerability profiles of minority subpopulations when they are surrounded by a non-disadvadvantaged majority. We observe that this is the case when we consider which census tracts in Texas fail to meet the qualifications to be considered vulnerable by these metrics.
Figure \ref{fig:corrplot} shows that there is little to no correlation between either (a) the percentage of a tract that is Indigenous or (b) the wildfire ignition risk at the census tract during a high-risk period on the likelihood that the tract is categorized as vulnerable. 
Since Indigenous populations make up a small fraction of the total population of each census tract, their relative disadvantage\footnote{Indigenous groups make up less than 1\% of the Texas population, but experience poverty at nearly double the rate of White Texans. This makes aggregation a likely explanation for the discrepancy in budget allocation and lack of characterization of Indigenous groups as belonging to Justice40 census tracts \cite{ACS_2022_poverty,ACS_2022_race,ACS_2022_Hispanic}} may be eclipsed by the fact that a majority of the census tract consists of groups that are generally \textit{not} disadvantaged (i.e., higher income, not a racial minority). This is what we define to be ``the curse of aggregation:" the loss of minority characteristics of a region's demography after aggregating population data to create a profile of that region.

We see that while northern and western Regions 1 and 2 in Figure \ref{fig: svi plots with load shedding} \textit{do}  experience the disproportionate fraction of the power loss experienced during this period, as we might have assumed from our observations in \S\ref{subsec: ercot}, the \mbox{\gls{cejst}} and \gls{svi} metrics do \textit{not} identify these tracts experiencing power outages as being particularly vulnerable. A natural question is whether this is due to (1) the fact that the \mbox{\mbox{\gls{cejst}}} indices are based on vulnerability from a national standpoint instead of considering statewide trends or (2) the fact that ~\mbox{\gls{cejst}} defines general climate vulnerability as opposed to wildfire vulnerability, specifically. To test this, we modified the \mbox{\mbox{\gls{cejst}}} criteria to only consider Texas percentiles of wildfire spread risk, as opposed to national percentiles of many types of climate risks (e.g., flood risk, agriculture loss rate, wildlife loss, etc). Even this ``modified ~\mbox{\gls{cejst}}" index fails to capture the vulnerability of the census tracts in the highlighted Regions~1 and~2, as shown in Figure \ref{fig:corrplot}.
%
Table \ref{table: misses} shows a subset of census tracts that the \mbox{\gls{cejst}}, modified \mbox{\gls{cejst}}, and \gls{svi} criteria fail to designate as vulnerable, despite being above the 50th percentile for low income and at very high percentiles for ignition risk (over 88\%), as well as the fraction of the census tract that is Indigenous (over 85\%). Indeed, there are cumulative network effects that neither the \gls{svi}, \mbox{\gls{cejst}}, nor modified \mbox{\gls{cejst}} can capture in their entirety.
\begin{table}[ht!]
\begin{center}
\begin{tabular}{ |c|C{1.78cm}|C{1.75cm}|C{1.75cm}| }
\hline
 \bf{GIDTR} & \bf{percentile below poverty line} & \bf{percentile Indigenous} & \bf{percentile ignition risk}\\ \hline
 110093 & 65 & 94 & 88\\ \hline
 120013 & 53 & 85 & 94\\ \hline
 120013 & 57 & 94 & 94\\ \hline
 220006 & 52 & 99 & 99\\ \hline
 220020 & 52 & 99 & 99\\ \hline
 220077 & 77 & 99 & 99\\ \hline
\end{tabular}
\end{center}
\caption{A subset of census tracts categorized as ``not vulnerable" for each of the three indices.}\label{table: misses}
\end{table}

This ``curse of aggregation" highlights a potential limitation of social vulnerability indices, which has also been corroborated by other works (e.g., \cite{cutter2024origin, javanmard2023impacts, hinojos2023social}). This work demonstrates an observation made by Cutter in \cite{cutter2024origin}: ``Clearly, the differing contexts of social vulnerability become averages when using aggregated data, and the single [social vulnerability index] score tends towards oversimplification.'' That is, our work demonstrates that Indigenous communities become ``averaged" with non-indigenous neighbors, and the index cannot identify their unique attributes.

In contrast, we find that {\sc Equity} objectives taking group size into account may ensure that minority subpopulations within census tracts are not overlooked. By design, such \gls*{mmf} objectives balance ``fair'' load shed outcomes for the considered set of groups; however, as discussed in Section \ref{subsec: group fair}, this may come at the cost of increased total load shedding. Despite this trade-off, the outcomes simulated from decisions made by a group-fair \gls{mmf} model can inform decision-makers on the conceivable bounds on possible load shed improvements for vulnerable communities, even if the model is not implemented in practice. Stakeholders, policymakers, and other decision-makers are ultimately responsible for determining what balance of competing objectives should be prioritized. 

\section{Conclusion}\label{sec: policy insights}
In this paper, we analyzed the performance of a Justice40 framework for making power line undergrounding and de-energization decisions on a synthetic transmission network modeling the \gls{ercot} region of Texas. We modeled this multi-criteria decision with the following considerations: (i) to have minimal (or close to minimal) total load shed in the network, (ii) so that the total wildfire ignition risk in the network remains within set limits, and (iii) to allocate benefits fairly across various groups. While considerations (i) and (ii) are relatively straightforward, achieving a ``fair" allocation policy (iii) was significantly more challenging.

First, we analyzed the load shed outcomes when implementing Justice40-style {\sc Policy} constraints on a \gls*{mip} model of the transmission network with the objective of minimizing total network load shed. Specifically, we considered constraints that proportionally allocate 40\% of the total budget to vulnerable communities or proportionally allocate 40\% of the total load shed reduction to vulnerable communities as defined by two vulnerability indices: the \gls{cejst} and the \gls{svi}. The optimal solution to these programs generally led to worse load shed outcomes for \textit{all} racial and socioeconomic groups. In particular, these {\sc Policy} constraints often fail to protect Indigenous populations who experience nearly double the poverty rate and double the anticipated load shed as the average rate of poverty and load shed in Texas. We posit that there are two key causes for the misalignment between the \textit{intent} of the Justice40 initiative with and the \textit{realized benefit} (or lack thereof) to communities experiencing high wildfire risk after implementing {\sc Policy} constraints. First, there is information loss due to data aggregation, and second, generalized vulnerability criteria lack context-specificity to be effective across different types of climate vulnerabilities and investments. 

Aggregation of demographic data is necessary to compute a ``vulnerability index" for a census tract. However, doing so homogenizes the demographic profile of the census tract instead of capturing its diversity. Hence, a disadvantaged group that makes up a minority of the census tract is likely to be overlooked if the majority of the census tract is not disadvantaged. This is likely why Indigenous groups were not more likely to live in disadvantaged census tracts than other groups, despite higher rates of poverty. In contrast, Hispanic groups in poor areas (a much larger group in Texas) could be appropriately identified as vulnerable by these vulnerability indices.

Furthermore, many of these vulnerability indices are intentionally general, as to be able to be utilized for a variety of different resource allocation projects. However, as this case study has demonstrated, the use of general vulnerability indices may not have the intended effect under specific policies. 
One option may be to define context-specific \textcolor{red}{vulnerability indices}, both by location and the type of risk that is being mitigated. However, even when vulnerability is defined in this way, the curse of aggregation may still apply, leaving some groups at disproportately high risk. In our case, when we evaluated a modified \mbox{\gls{cejst}} index which only considered Texas wildfire spread risk (as opposed to national percentiles of a plethora of different, potentially unrelated climate risks), we saw similarly poor load shed outcomes for Indigenous groups. In short, creating such a context-specific definition of vulnerability may improve its efficacy, but doing so is not necessarily sufficient for adequately prioritizing at-risk populations. 

In this case study, the use of aggregated vulnerability indices did not appropriately allocate resources to vulnerable minority populations, which motivated the study of explicit \textit{group protections.} The percentage-based {\sc Equity} objective is group-size conscious, which prevents minority subpopulations within census tracts from being overlooked during the optimization routine in a way that vulnerability indices cannot. Furthermore, such an objective, by construction, balances ``fair" load shed outcomes with total load shed in the network and does so without any intrinsic bias about which group should be prioritized. That is, the model does not know \textit{a priori} that Indigenous groups need to be prioritized, the objective simply tries to minimize the maximum percentage4 of load shed experienced by any group. While group-size-conscious allocation mechanisms may not be advantageous in certain contexts, the inability of the \mbox{\gls{cejst}} and \mbox{\gls{svi}} metrics to identify Indigenous disadvantage makes this particular context of budget allocation for power line undergrounding a feasible candidate for a group-size-conscious allocation policy. Furthermore, we emphasize that an \gls{mmf} framework is not the only {\sc Equity} objective that could be used in this model; other objectives like $L_p-$norm fairness or ordered norms can also be used to balance the percentage load shed across various groups \cite{bonald2015multi,kalai1975other,bertsimas2011price, gupta2023lp}. \textcolor{red}{However, the MMF framework, which minimizes the $L_{\infty}$ norm, would decrease the percentage load shed on the most impacted group maximally out of these choices of $L_p$-norm objectives.} 

Though our work underscores how considering ethnicity and race using the MMF framework may help ensure certain groups are not disproportionally impacted by power loss, the legality of the use of explicit group-level protections is unclear, particularly for public investments. For some states, explicit use of race when using the \mbox{\gls{cejst}} is illegal (e.g., California Proposition 209) \cite{walls2024implementation}, whereas in others, there is some gray area due to the disparate level of benefit of federal funding when \textit{not} considering race. In the United States, there are some cases that lend legal plausibility to such approaches, e.g., American Indigenous communities \textit{living on federally recognized tribal lands} are often given special protections. 
Stakeholders, policymakers, and other decision-makers are responsible for which protections or populations should be prioritized, and the feasibility of these choices.
This work is intended to aid these decision-makers in their analysis process, and identify trade-offs that exist when working towards equitably allocating benefits. 

In summary, our analysis reveals the limitations and potential of vulnerability indices like CEJST and SVI in effectively directing resources to high-risk, vulnerable communities—in this case, Indigenous populations—within climate resilience projects. Aggregating demographic data into indices risks overlooking smaller disadvantaged groups within census tracts that may not fit broad vulnerability criteria, such as Indigenous communities. This gap signals the need to consider context-sensitive, group-conscious approaches like the MMF framework. 
Our work thus contributes empirical insights into the trade-offs of different vulnerability frameworks, equipping policymakers and stakeholders with an informed basis to advance more equitable climate resilience initiatives.

\appendix

\section{Data Sources}
A number of data sets were used for this project. Information on each data source is described below. The methodology used to map these data sets to the synthetic network is described in Appendix \ref{appendix: bus features}.
    \subsection{US Census Data}\label{subapp: census data}
    We use the census tract definitions that were in effect from 2010-2020. 
    To obtain data by census tract for total population, median income, and number of individuals in each racial group, we used the 2019 Planning Database, version 2 \citep{demographic_data}. To get the latitude and logitude of the center of population for each census tract, we used center of population data from the 2010 Decennial census \citep{COP_data}, which was the most recently available center of population data for this districting. 
    
    \subsection{The Justice40 Initiative}\label{subappendix: j40}
    The Justice40 initiative defines ``disadvantaged communities'' \citep{exec_order2021} at the US Census Tract level on the basis of income, energy access, housing access, and environmental burden, based on data provided by the \mbox{\gls{cejst}}. In general, according to the US Council on Environmental Quality \cite{ScreeningTool_2022}, to qualify for the Justice40-designation, a census tract must either (a)~``meet the thresholds for at least one of the tool’s categories of burden,'' or (b)~be ``within the boundaries of a federally recognized tribe.''
    Note that the latter condition \textit{does not apply} to indigenous populations living outside recognized tribal borders. 
    The former condition generally requires that a census tract be above the 65th percentile nationwide for the percentage of the population which is considered low-income and be above the 90th percentile nationwide for one of the many different types of climate burden. 
    Since the \gls{cejst} considers national percentiles and Texas has higher climate risk and poverty levels than much of the United States, we see that nearly 50\% of census tracts served by \mbox{\gls*{ercot}} are categorized as vulnerable. Furthermore, this metric involves climate risks that are not limited to wildfire risk, so in this study, we also computed a modified \mbox{\gls{cejst}} indexing that categorizes a tract as vulnerable if the tract is at or above the 50th percentile of Texas census tracts for the percentage of the tract which is low-income and at or above the 75th percentile of Texas census tracts for the wildfire risk subcategory of the Justice40 dataset. With this definition, only around 11\% of census tracts are considered vulnerable. It is also worth noting that the way that the \gls{cejst} computes wildfire risk is based on wildfire spread models \citep{ScreeningTool_2022,Firststreet.org_2023}, instead of the ignition risk models used in this study \citep{USGS_WFPI}.
    
    \subsection{The CDC/ATSDR Social Vulnerability Index}\label{subappendix: svi}
    The \gls*{cdc} and \gls*{atsdr} have a joint metric of social vulnerability given in terms of an \gls*{svi} \citep{CDC_SVI_data}. The goal of this \gls*{svi} is to designate communities which may have additional difficulty coping with a disaster event. This \gls*{svi} considers four main ``themes'' of risk: socioeconomic status, household characteristics, racial and ethnic minority status, and housing type/transportation. In this study, we classify a census tract as vulnerable if the census tract is is at or above the 75th percentile of burden out of all the census tracts in Texas for at least one of these four themes. We use the 2010 \gls*{svi} dataset to remain consistent with the 2010 census tracts \citep{CDC_SVI_data}. 

    \subsection{USGS Wildland Fire Potential Index \& Wildfire Risk Values}\label{appendix: usgs data}
    The \gls*{usgs} \gls*{wfpi} is a data set compromised of unitless risk values ranging from 0 to 247 for each 1~km by 1~km ``pixel'' of the United States. These values are updated daily along with a 7-day forecast of expected risk values. The USGS bases this data on the following~\citep{USGS_WFPI}:
    \vspace{-.3cm}
    \setlength{\columnsep}{4pt} 
    \begin{multicols}{2}
    \begin{itemize}[noitemsep, nolistsep]
        \item Maximum Live Ratio,
        \item Dead Fuel Moisture,
        \item Fuel Model,
        \item Wind Reduction Factor,
        \item Vegetation Index,
        \item Relative Greenness,
        \item Dead Fuel Moisture,
        \item Wind Speed,
        \item Rain,
        \item Dry Bulb Temperature.
    \end{itemize}
    \end{multicols}
    \vspace{-.3cm}
\noindent The USGS WFPI provides a proxy for the risk of ignition from electric infrastructure since higher WFPI values have, historically, correlated to larger fires and fires that have spread to burn more area~\citep{USGS_WFPI}. 

Historically, wildfire season in the western United States typically spans from late summer to early fall; however, recent wildfire seasons have been lengthening \citep{USDA_wildfire_season}.
Therefore, our analyses use data from June~1 to October~31, which we will refer to as the wildfire season. We assign a unitless wildfire risk value $r_{\ell, d}$ for each line $\ell \in \mathcal{L}$ for each day $d \in \mathcal{D}$ in the considered network in the wildfire season over three years (2019, 2020, and 2021). To find this value, we find the average pixel risk, $\bar{r}_p$ by taking the mean of all pixel values on all lines from the data used. We find the standard deviation on this data as well, $\sigma_p$. We then define a high-risk pixel to be an pixel with a value more than one standard deviation above the mean:
\begin{equation}
    r^h_{p,l,d} = 
    \left\{
    \begin{array}{lr}
        r_{p,l,d}, & \text{if } r_{p,l,d} \geq \bar{r_p} + \sigma_p\\
        0, & \text{if } r_{p,l,d} < \bar{r}_p + \sigma_p
    \end{array}
    \right\}.
\end{equation}
We calculate the risk value $r_{\ell, d}$ by integrating the high-risk pixel values along each line. This method balances the risk contributions from both long line lengths and underlying risk of the terrain. ``Line length" is a characteristic that has been correlated with higher ignition risk \citep{waseem2021} but this processing avoids a situation where a long line with relatively low risk along the entire length appears much riskier than a shorter line with points of much higher ignition risks. This is a method adapted from~\citep{piansky2025hicss}.

For each day simulated, we first determine if the wildfire threat is high enough to necessitate de-energizing lines via a threshold on the total risk during that day. 
Let $R_d$ be the total wildfire risk the network poses if all lines $\ell \in \mathcal{L}$ are energized on day $d$, i.e.: $R_d = \sum_{\ell \in \mathcal{L}} r_{\ell,d}$. 
In our assessment methodology, operators are required to reduce the total risk of the network by making line de-energization decisions during any day for which $R_d \geqslant R_{\text{PSPS}}$, where $R_{\text{PSPS}}$ is a specified system-wide de-energization threshold. 
Conversely, if $R_d < R_{\text{PSPS}}$, then the risk the network poses is not great enough to require the widespread de-energization of lines. For the purposes of this paper, $R_{\text{PSPS}}$ is set to $6\times 10^8$. Results and figures indicate what overall threshold was applied to the network. 

Two more thresholds are used for the network to split $\mathcal{L}_d$ in to $\mathcal{L}^{\text{high}}_d$, $\mathcal{L}^{\text{med}}_d$, and $\mathcal{L}^{\text{low}}_d$. These thresholds, $R_{\text{high}}$ and $R_{\text{low}}$, are used to indicate the highest acceptable risk before lines must be de-energized or undergrounded and the lowest risk below which lines are not considered for undergrounding or de-energization:
\begin{equation}
\begin{aligned}
    \mathcal{L}^{\text{high}}_d &\triangleq \{ \ell \in \mathcal{L} \ | \ r_{\ell, d} \geq R_{\text{high}} \}\\
    \mathcal{L}^{\text{med}}_d & \triangleq \{ \ell \in \mathcal{L} \ | \ R_{\text{low}} \leq r_{\ell, d} < R_{\text{high}} \}\\
    \mathcal{L}^{\text{low}}_d & \triangleq \{ \ell \in \mathcal{L} \ | \ r_{\ell,d} < R_{\text{low}} \}\\
\end{aligned}
\end{equation}
For the results in this paper, $R_{\text{high}}$ and $R_{\text{low}}$ are set to $1\times10^6$ and $1\times 10^0$, respectively. Note that this means all lines with risk values greater than one are allowed to be de-energized or undergrounded. These values were chosen to allow enough lines to be candidates for undergrounding such that the MIP produces non-trivial solutions.

\section{Mapping Demographic Features to Buses}\label{appendix: bus features}
%
%
\begin{table*}[h!]
\centering
\caption{List of disadvantaged census tracts that every vulnerability index categorized as not vulnerable.}
\label{table:full_table_misses}
\begin{tabular}{lcccc}
\toprule
GIDTR & \makecell{percentile\\uninsured} & \makecell{percentile\\impoverished} & \makecell{percentile\\indigenous} & \makecell{percentile\\ignition risk} \\
\midrule
110093 & 69.8 & 64.8 & 94.5  & 88.3 \\
110965 & 73.9 & 51.0 & 76.9  & 75.1 \\
120013 & 32.9 & 52.7 & 84.9  & 95.7 \\
120063 & 32.9 & 52.7 & 84.9  & 84.3 \\
120112 & 52.8 & 76.2 & 94.4  & 89.6 \\
120171 & 55.7 & 66.5 & 90.4  & 83.9 \\
140066 & 9.2  & 80.1 & 75.1  & 82.5 \\
210028 & 20.3 & 58.4 & 97.5  & 89.6 \\
210029 & 50.0 & 56.8 & 94.2  & 94.4 \\
210030 & 69.4 & 62.2 & 93.2  & 91.7 \\
210053 & 50.5 & 62.4 & 96.2  & 78.3 \\
210075 & 30.7 & 61.2 & 77.5  & 80.4 \\
210207 & 25.7 & 59.4 & 89.6  & 84.2 \\
210293 & 59.8 & 61.0 & 94.5  & 88.7 \\
220006 & 78.5 & 52.4 & 100.0 & 99.3 \\
220020 & 78.5 & 52.4 & 100.0 & 99.7 \\
220026 & 54.6 & 56.5 & 95.7  & 95.7 \\
220057 & 29.9 & 57.9 & 98.0  & 97.2 \\
220077 & 25.5 & 76.7 & 99.4  & 99.3 \\
220109 & 38.0 & 52.9 & 92.0  & 95.4 \\
240084 & 44.4 & 64.3 & 75.3  & 94.7 \\
240121 & 67.0 & 64.5 & 86.4  & 93.7 \\
240148 & 67.0 & 64.5 & 86.4  & 91.8 \\
240160 & 54.1 & 62.8 & 86.9  & 84.8 \\
240162 & 12.3 & 74.4 & 98.1  & 93.2 \\
240193 & 18.7 & 51.6 & 91.2  & 84.7 \\
\bottomrule
\end{tabular}
\end{table*}
We  match census tracts to buses with load based on the distance between the bus and the population center of the census tract. Let $\tracts$ denote the set of census tracts and $\buses$ represent the set of buses with nonzero load. Let $d_{cn}$ represent the distance between the center of population of census tract $c$ and bus $n$. Each census tract $c \in \tracts$ has a feature vector $f_c$. 
We take a three-pass approach.  
\begin{enumerate}[noitemsep,nolistsep]
    \item For every census tract $c \in \tracts$, we initialize the radius $r_c$ as the minimum distance from $c$ to any other bus in the transmission  network. If a bus $n$ is within the radius $r_c$ for any $c \in \tracts$, we say that the bus has been \textit{assigned}.
    \item For any bus $n$ that has not yet been assigned, we find the closest census tract $c$ to $n$. Let this closest distance be given by $r_n$. We then update $r_c \leftarrow \max\{r_c, r_n\}$. Now, bus $n$ has been assigned.
    \item For every tract $c \in \tracts$, we consider the subset of buses $\buses_{r_c}^c$ within a distance $r_c$ from $c$. We divide the population of $c$ between each bus $n \in \buses_{r_c}^c$ proportionally based on their relative distance from $c$. That is, the fraction of $f_c$ that is assigned to bus $n' \in \buses_{r_c}^c$ is $a_{cn'} = \frac{d_{cn'}}{\sum_{i \in \buses_{r_c}^c}d_{ci}}.$
\end{enumerate}

\noindent At the termination of this algorithm, we have a sparse matrix, $A$, where each entry $a_{ij}$ is the ``fraction'' of each census tract $i$ that is assigned to each bus $j$. Finally, we have that the demographic feature vector for bus $n \in \buses$ is given by $f_n = \sum_{c \in \tracts}f_c \cdot a_{cn}$.

\subsection{``Missed'' Vulnerable Tracts}\label{appendix:missed tracts}
Table \ref{table: full table, misses} gives the set of Texas census tracts which are at or above the 75th percentile for the fraction of the population that is indigenous, at or above the 75th percentile for wildfire ignition risk (derived from \gls*{wfpi} data), and at or above the 50th percentile for number of people below the poverty line, but that were \textit{not} characterized as vulnerable by the \gls{cejst}, modified \gls{cejst}, and \gls{svi} criteria. We also show the percentile without health insurance for reference. One might note from Table \ref{table: full table, misses} that the impoverished percentile is not often \textit{too} high; indeed, the maximum percentile impoverished in the table is 80.1 for GIDTR 140066, and most values are between the 50th and 60th percentiles. Given that indigenous poverty rates are nearly double that of the overall population \cite{ACS_2022_poverty,ACS_2022_race,ACS_2022_Hispanic}, this indicates that indigenous populations are almost always the minority of the census tract that they live in, and the majority population of those tracts is not disadvantaged.  


\subsection{Income and Power Outages in Texas}\label{appendix:income and ls}
\begin{figure}
    \centering
    \includegraphics[width=0.6\linewidth]{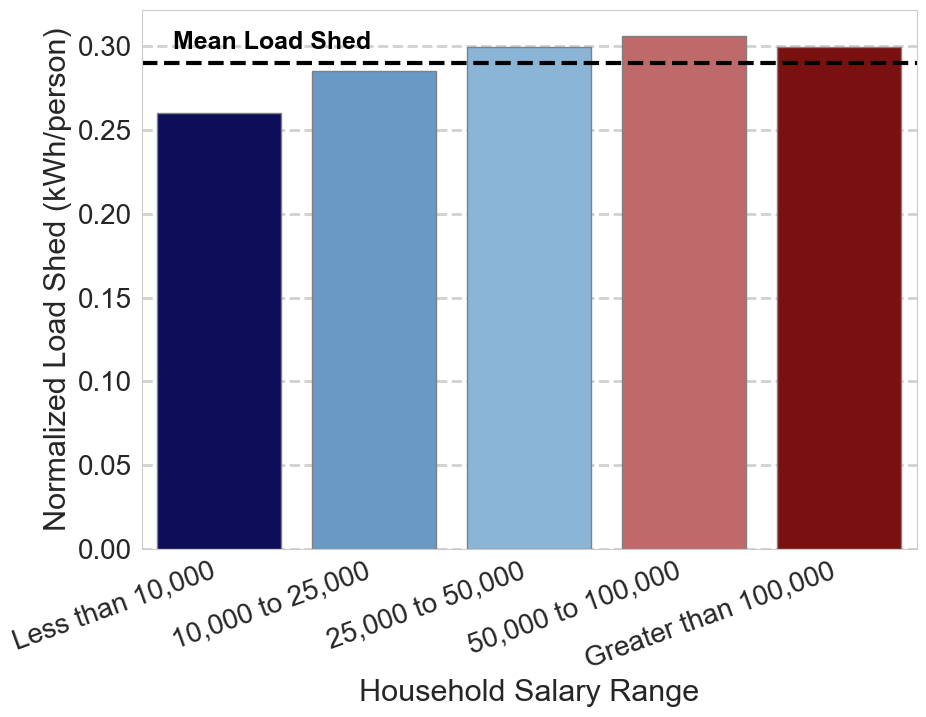}
    \caption{Income has very little impact on the percent of load demanded that is shed in the baseline, no budget case (BL-M0). In fact, there appears to be a slight correlative effect wherein higher income groups see a higher percentage of their load shed, although the increase is very slight.}
    \label{fig:ls by inc group}
\end{figure}

\begin{figure}
    \centering
    \includegraphics[width=0.6\linewidth]{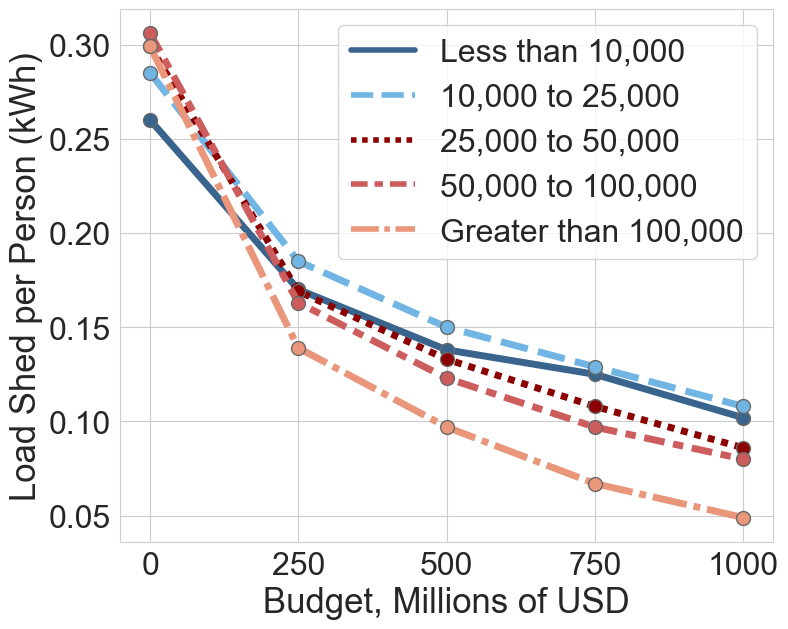}
    \caption{As budget increases, the discrepancy in percent of load demanded that is shed between income groups remains small; however, the trend reverses itself where after a \$250 million allocation, the wealthiest groups now have marginally less of their load that is shed.}
    \label{fig:ls by inc group over budget}
\end{figure}

Per our simulations, income level has very little effect on load shed outcomes in Texas. Figure \ref{fig:ls by inc group} shows that there is very little discrepancy in the percent of load demanded that is shed between income groups for the baseline, no budget case. In fact, the small discrepancy that \textit{does} exist most favorably impacts the lower-income groups.
While Figure \ref{fig:ls by inc group over budget} shows this slight trend reverse once budget is allocated, the discrepancy still remains quite small. The finding that minority status is a greater indicator of the likelihood of a power outage than income mirrors those of \cite{shah2023inequitable}.

\section{Optimization Results}\label{appendix: extended results}

\begin{figure*}[t]
    \centering
    \begin{subfigure}[t]{.49\textwidth}
      \centering
      \includegraphics[width=.98\linewidth]{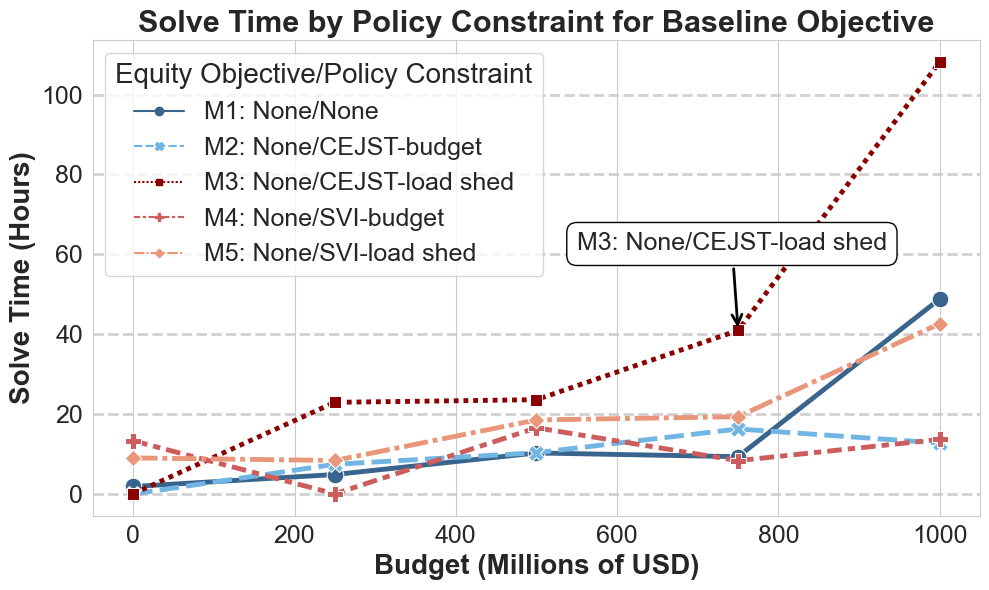}  
      \caption{Solution times by budget for simulations with {\sc policy} constraints only.}
      \label{fig: subfig: baseline_times}
    \end{subfigure}
    \hfill
    \begin{subfigure}[t]{.49\textwidth}
      \centering
      \includegraphics[width=.98\linewidth]{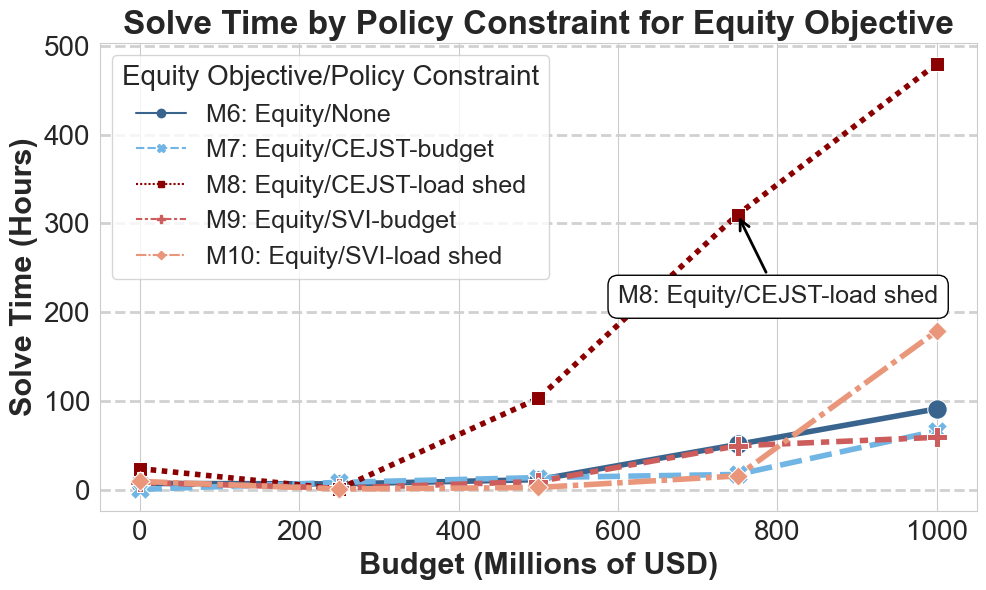}  
      \caption{Solution times in hours by budget for simulations with the {\sc equity} objective and {\sc policy} constraints.}
      \label{fig: subfig: equity_times}
    \end{subfigure}
    \\
    \begin{subfigure}[t]{.49\textwidth}
      \centering
      \includegraphics[width=.98\linewidth]{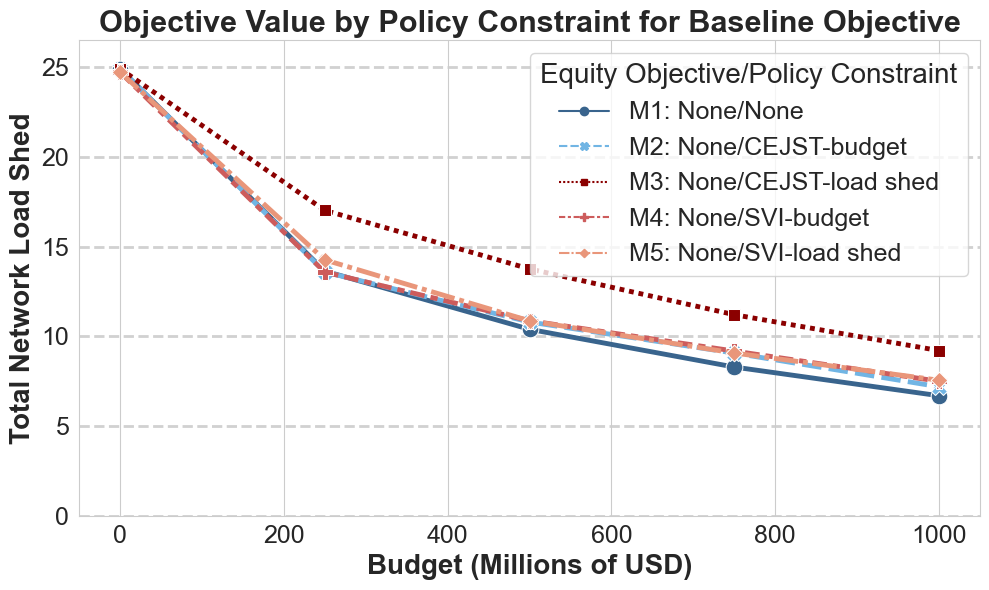}  
      \caption{Objective values by budget for simulations with {\sc policy} constraints only.}
      \label{fig: subfig: baseline_objs}
    \end{subfigure}
    \hfill
    \begin{subfigure}[t]{.49\textwidth}
      \centering
      \includegraphics[width=.98\linewidth]{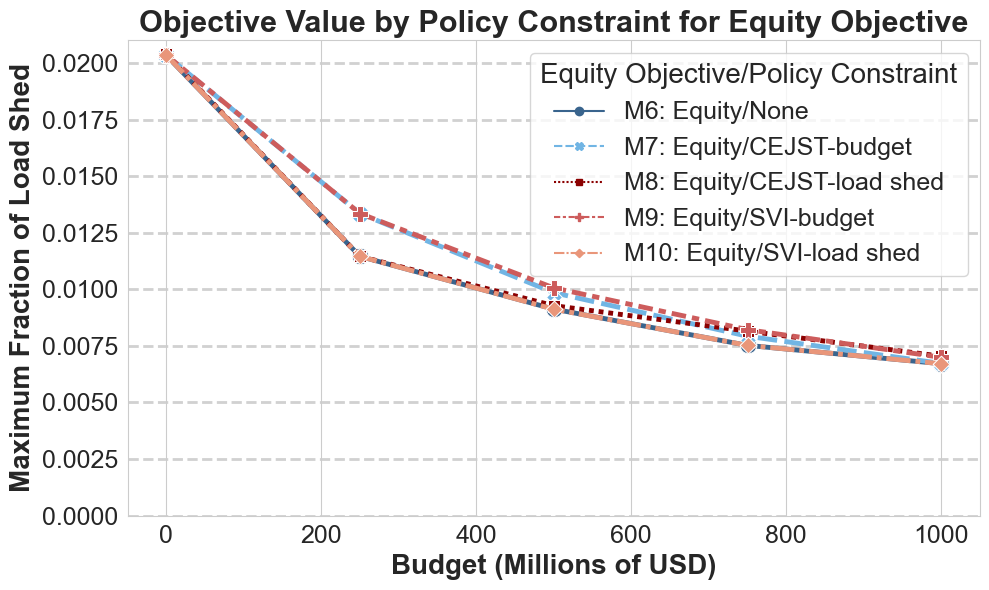}  
      \caption{Objective values by budget for simulations with the {\sc equity} objective and {\sc policy} constraints.}
      \label{fig: subfig: equity_objs}
    \end{subfigure}
    \caption{Solve times and objective values for each optimization model.}
    \label{fig: computational_results}
    \vspace{.5cm}
\end{figure*}

\subsection{Optimization Software, Set Up, and Solve Time}\label{appendix: optimization outcomes}
Optimization problems were solved using Gurobi~10.0.0~\citep{gurobi}. To implement the optimization formulations, we use Julia 1.8.0~\citep{bezanson2017julia} with JuMP~1.18.1~\citep{DunningHuchetteLubin2017} along with the data input functionality of PowerModels.jl 0.21.0~\citep{coffrin2018}. Simulations were completed on the Partnership for an Advanced Computing Environment (PACE) at the Georgia Institute of Technology~\citep{PACE}. For nonzero budgets, we warm-started the simulation with the results from the same model on the previous budget. All simulations are run for 5 days or a MIP gap of 1\%. Any simulations that are outside the 1\% MIP gap are run for an additional 5 days with a warm-start of the last found incumbent. Any simulations that are still outside the 1\% MIP gap are run for an additional 10 days, again warm-started from the last-found incumbent, or until they reach a 1\% MIP gap. This 10-day computation time is done with Gurobi's MIP focus parameter set to 3 to prioritize improvement in the best bound. After 20 days of computation time per scenario, MIP gaps are reported, with all scenarios finishing within a 5\% MIP gap. Solution times are shown in Figures~\ref{fig: subfig: baseline_times} and \ref{fig: subfig: equity_times}. 

While not all simulations converge to a 1\% MIP gap, we note the monotonically decreasing objective across budgets; see Figures~\ref{fig: subfig: baseline_objs} and \ref{fig: subfig: equity_objs}. Note the objective value for the baseline objective displays total network load shed while the {\sc equity} objective values portray the maximum percentage of demanded load that is shed for a given group, resulting in different scales. All simulations with a baseline objective solve to within a 1\% MIP gap. One combinations of the {\sc equity} objective with the CEJST load shed constraint (E-M8) converges to within a 5\% MIP gap. All other cases with the {\sc equity} objective converge to within a 1\% MIP gap. 

\subsection{Load Shed Results}\label{subsection: extended results} 
\begin{figure*}[t]
    \centering
    \begin{subfigure}[t]{.325\textwidth}
        \centering
        \includegraphics[width=.99\linewidth]{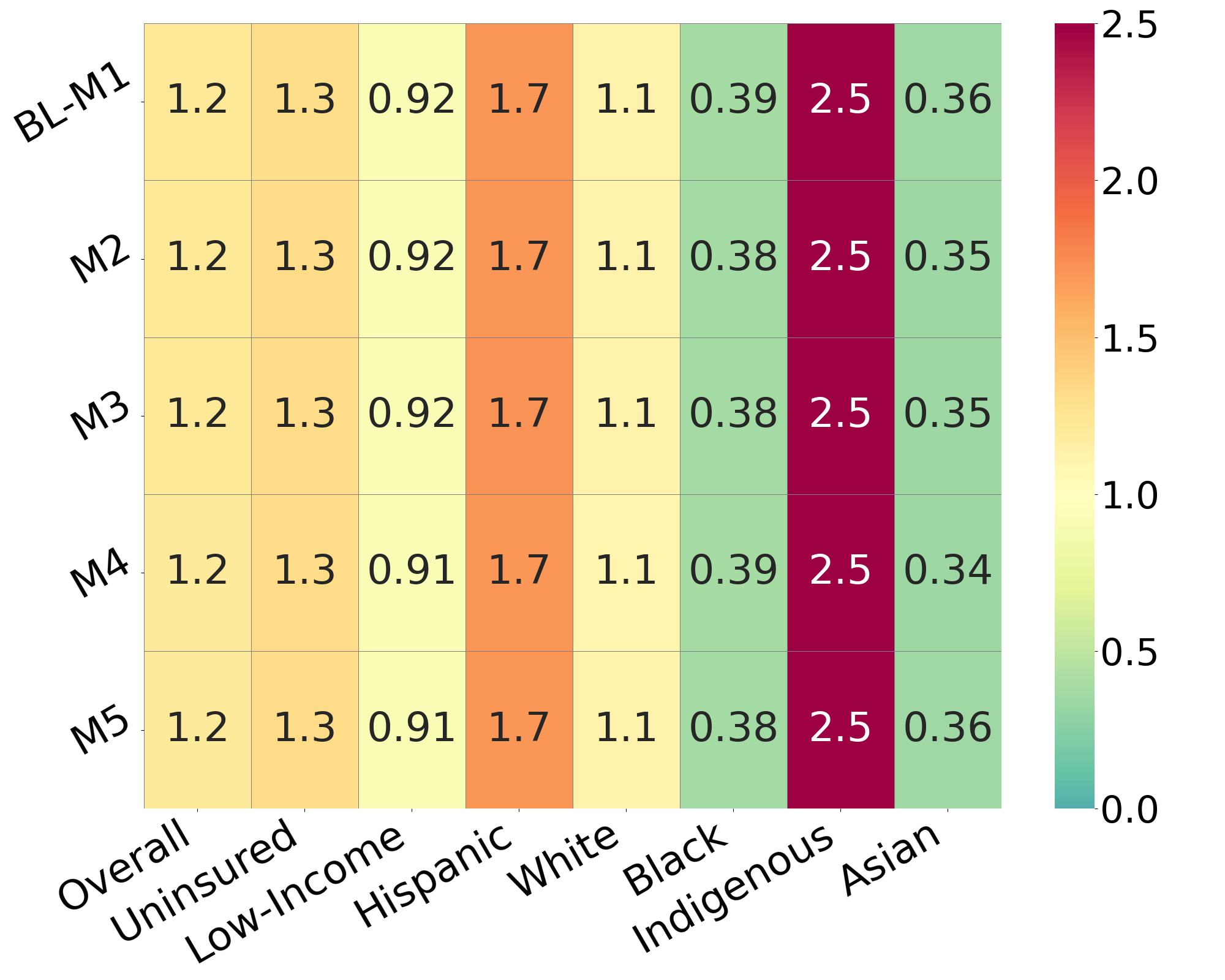}
        \caption{No {\sc Equity} Objective, \$0 budget}
        \label{fig: subfig: ls, B0, Baseline}
    \end{subfigure}
    \begin{subfigure}[t]{.325\textwidth}
        \centering
        \includegraphics[width=.99\linewidth]{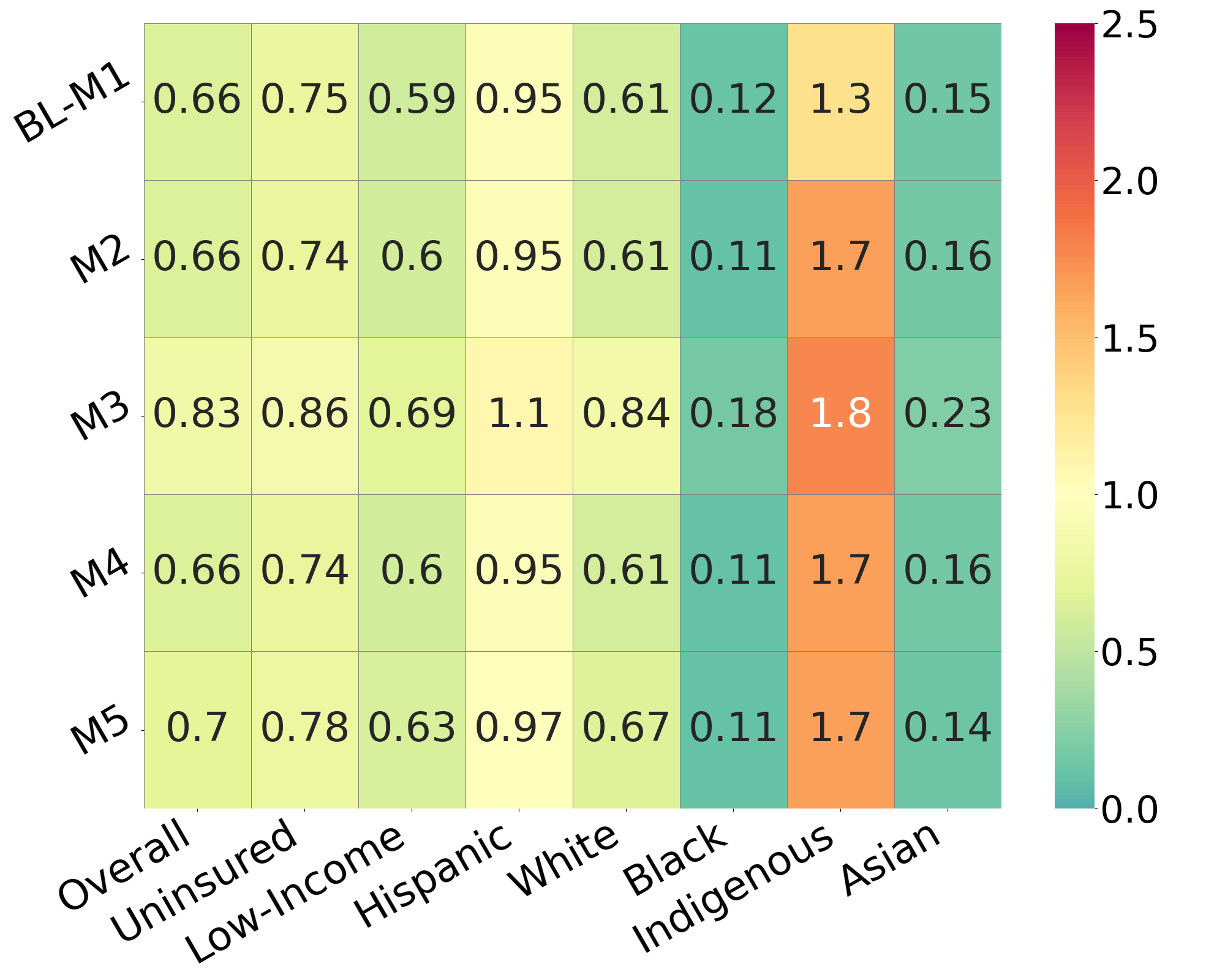}
        \caption{No {\sc Equity} Objective, \$250 million budget}
        \label{fig: subfig: ls, B250, Baseline}
    \end{subfigure}
    \begin{subfigure}[t]{.325\textwidth}
        \centering
        \includegraphics[width=.99\linewidth]{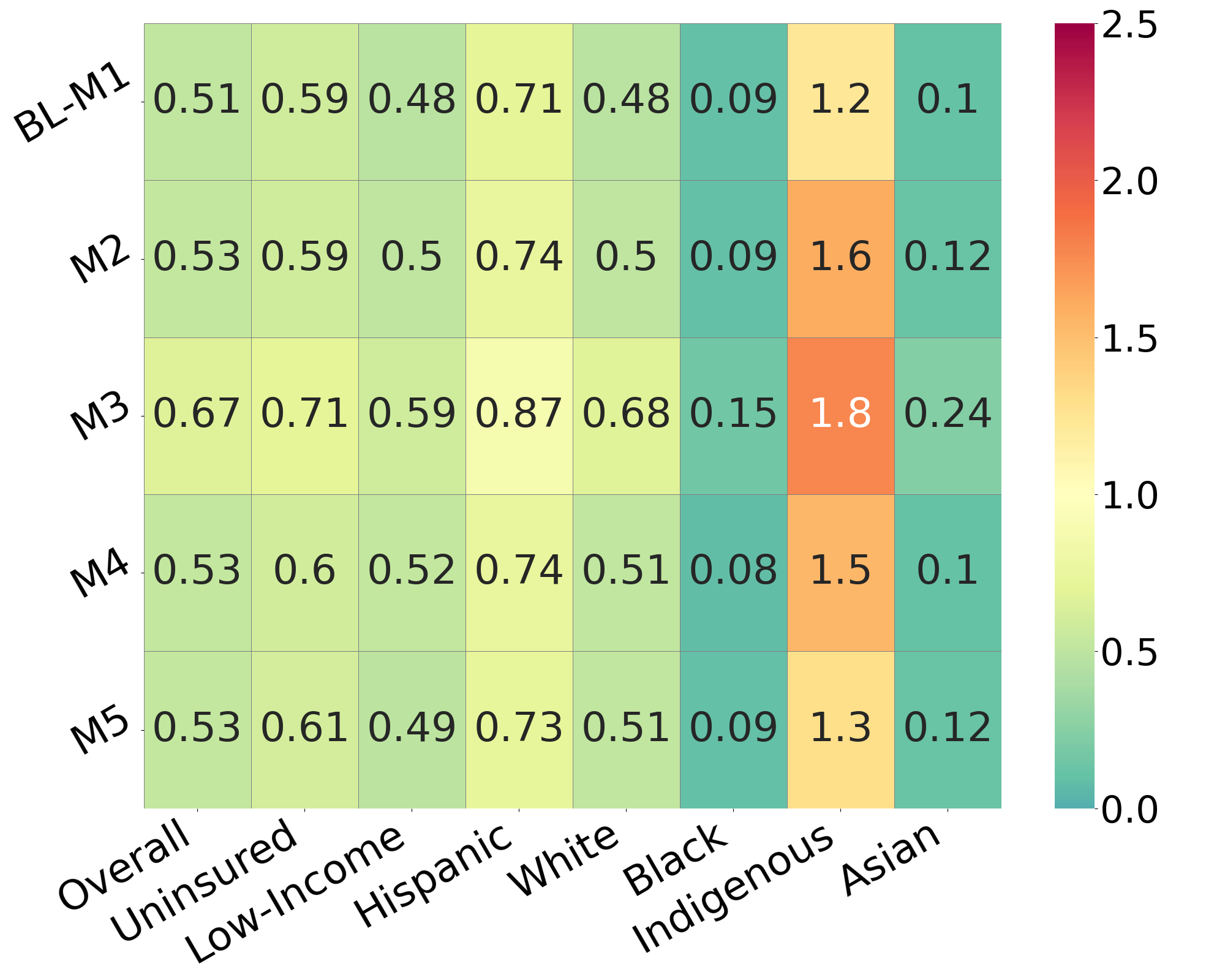}
        \caption{No {\sc Equity} Objective, \$500 million budget}
        \label{fig: subfig: ls, B500, Baseline}
    \end{subfigure}\\
    \begin{subfigure}[t]{.325\textwidth}
        \centering
        \includegraphics[width=.99\linewidth]{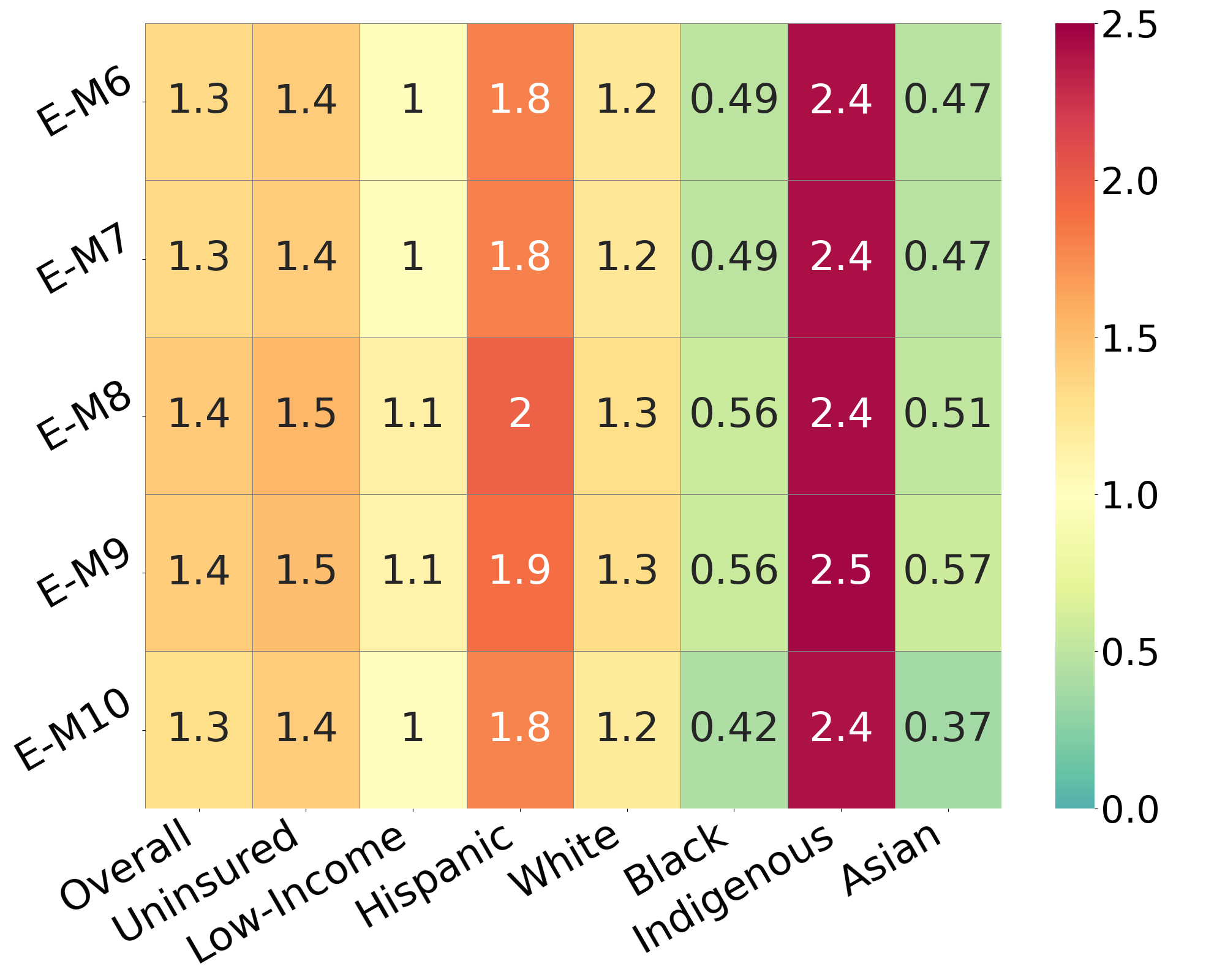}
        \caption{{\sc Equity} Objective, \$0 budget}
        \label{fig: subfig: ls, B0, MMF}
    \end{subfigure}
    \begin{subfigure}[t]{.325\textwidth}
        \centering
        \includegraphics[width=.99\linewidth]{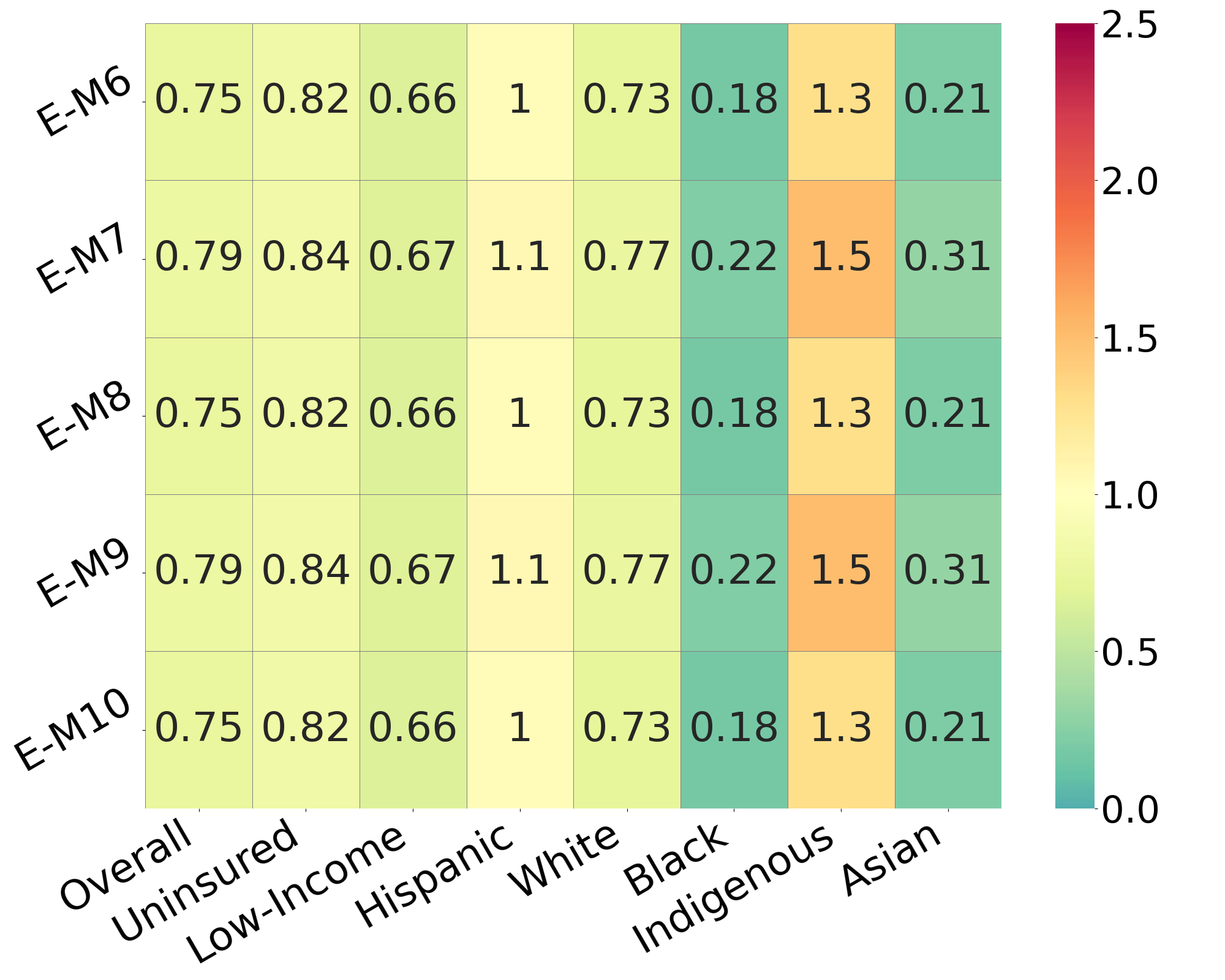}
        \caption{{\sc Equity } Objective, \$250 million budget}
        \label{fig: subfig: ls, B250, MMF}
    \end{subfigure}
    \begin{subfigure}[t]{.325\textwidth}
        \centering
        \includegraphics[width=.99\linewidth]{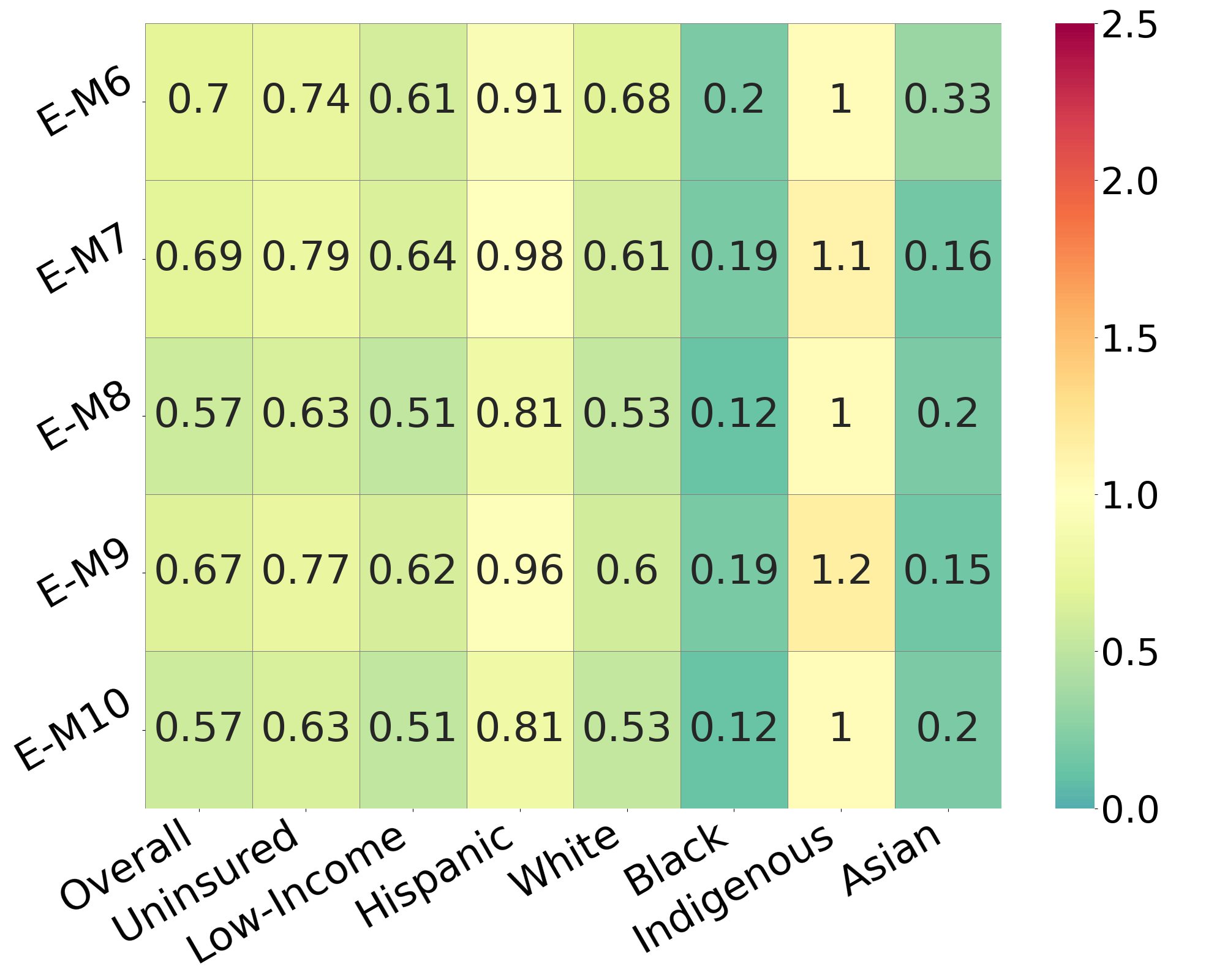}
        \caption{{\sc Equity } Objective, \$500 million budget}
        \label{fig: subfig: ls, B500, MMF}
    \end{subfigure}
    \caption{Percent of load demanded that is shed across budgets of \$0, \$250 million, and \$500 million. The first row shows {\sc Policy} constraints alone, and the second shows {\sc Policy} constraints and the {\sc Equity} objective.}
    \label{fig: pct loadshed}
\end{figure*}

\begin{figure*}[t]
    \centering
    \begin{subfigure}[t]{.325\textwidth}
        \centering
        \includegraphics[width=.99\linewidth]{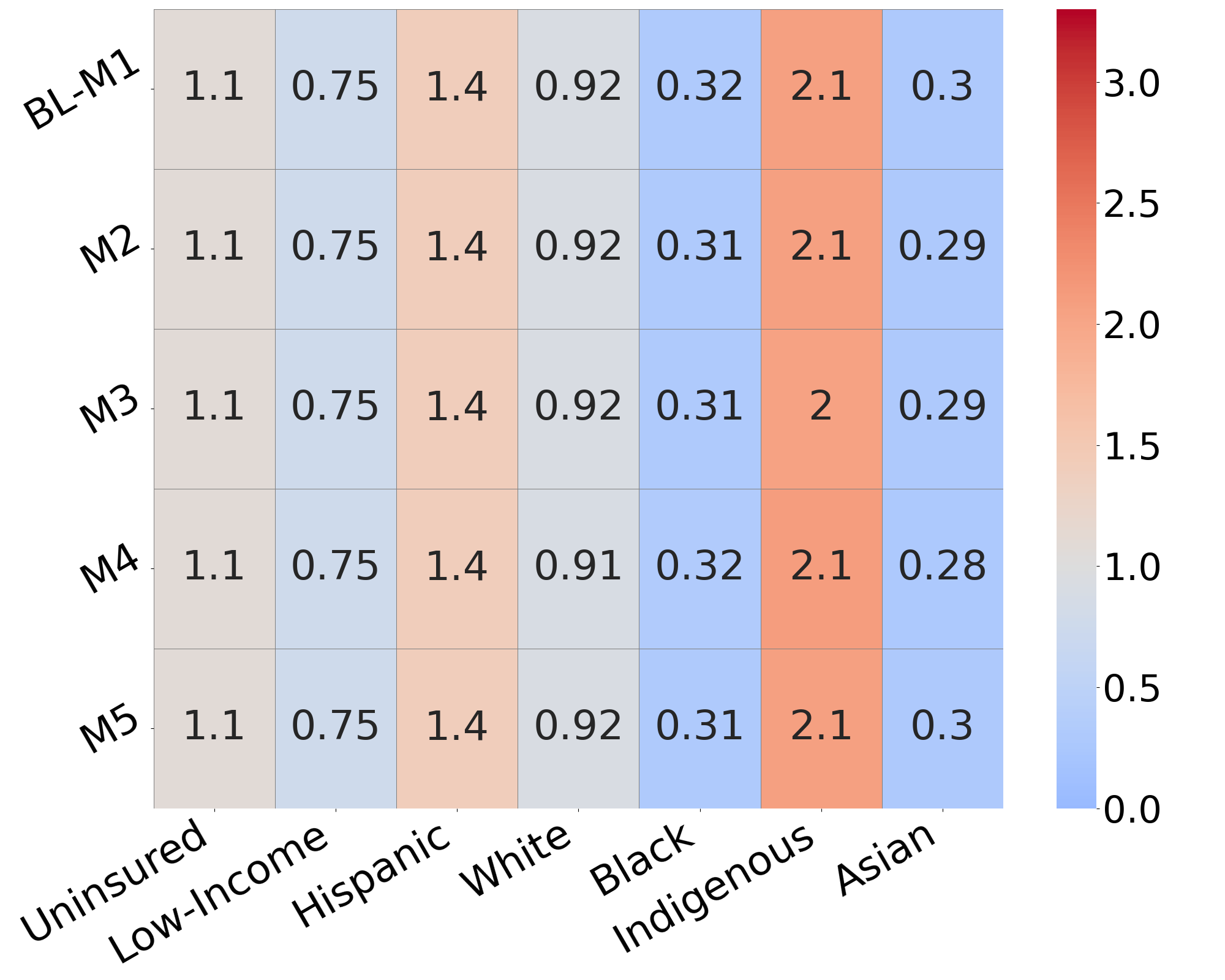}
        \caption{Relative unfairness in simulations with no {\sc Equity } Objective, \$0 budget}
        \label{fig: subfig: unfairness, B0, Baseline}
    \end{subfigure}
    \begin{subfigure}[t]{.325\textwidth}
        \centering
        \includegraphics[width=.99\linewidth]{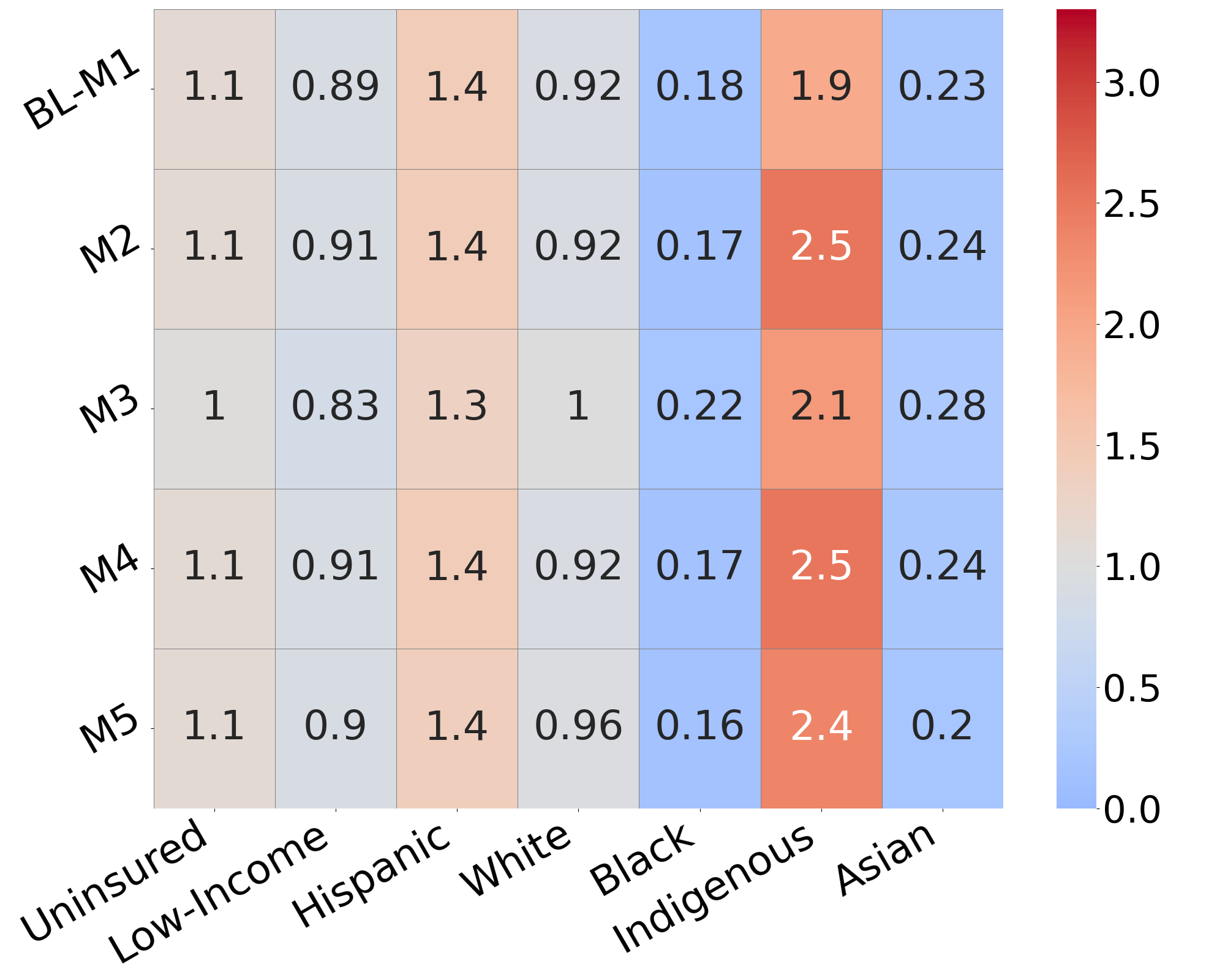}
        \caption{Relative unfairness in simulations with no {\sc Equity } Objective, \$250 budget}
        \label{fig: subfig: unfairness, B250, Baseline}
    \end{subfigure}
    \begin{subfigure}[t]{.325\textwidth}
        \centering
        \includegraphics[width=.99\linewidth]{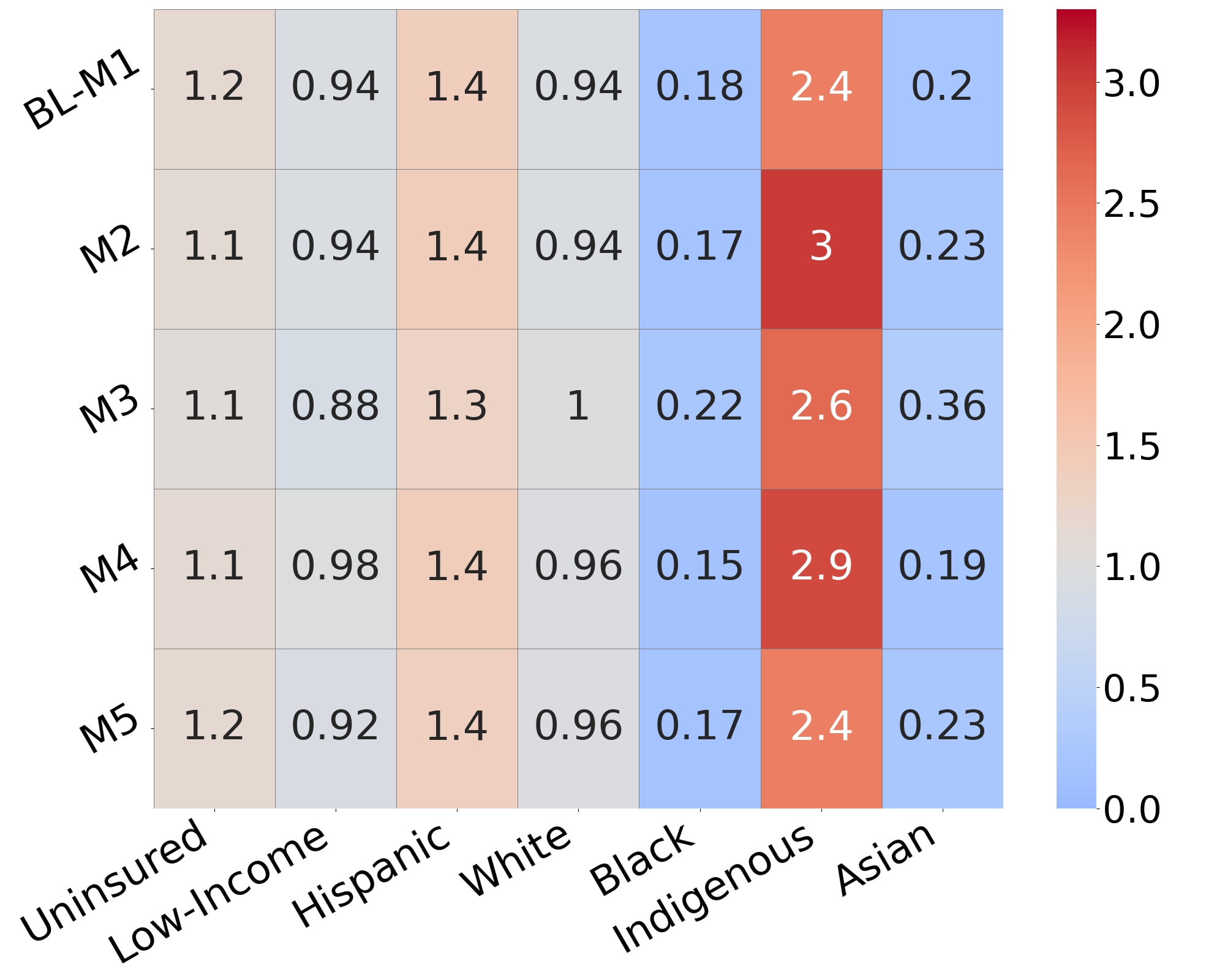}
        \caption{Relative unfairness in simulations with no {\sc Equity } Objective, \$500 budget}
        \label{fig: subfig: unfairness, B500, Baseline}
    \end{subfigure}\\
    \begin{subfigure}[t]{.325\textwidth}
        \centering
        \includegraphics[width=.99\linewidth]{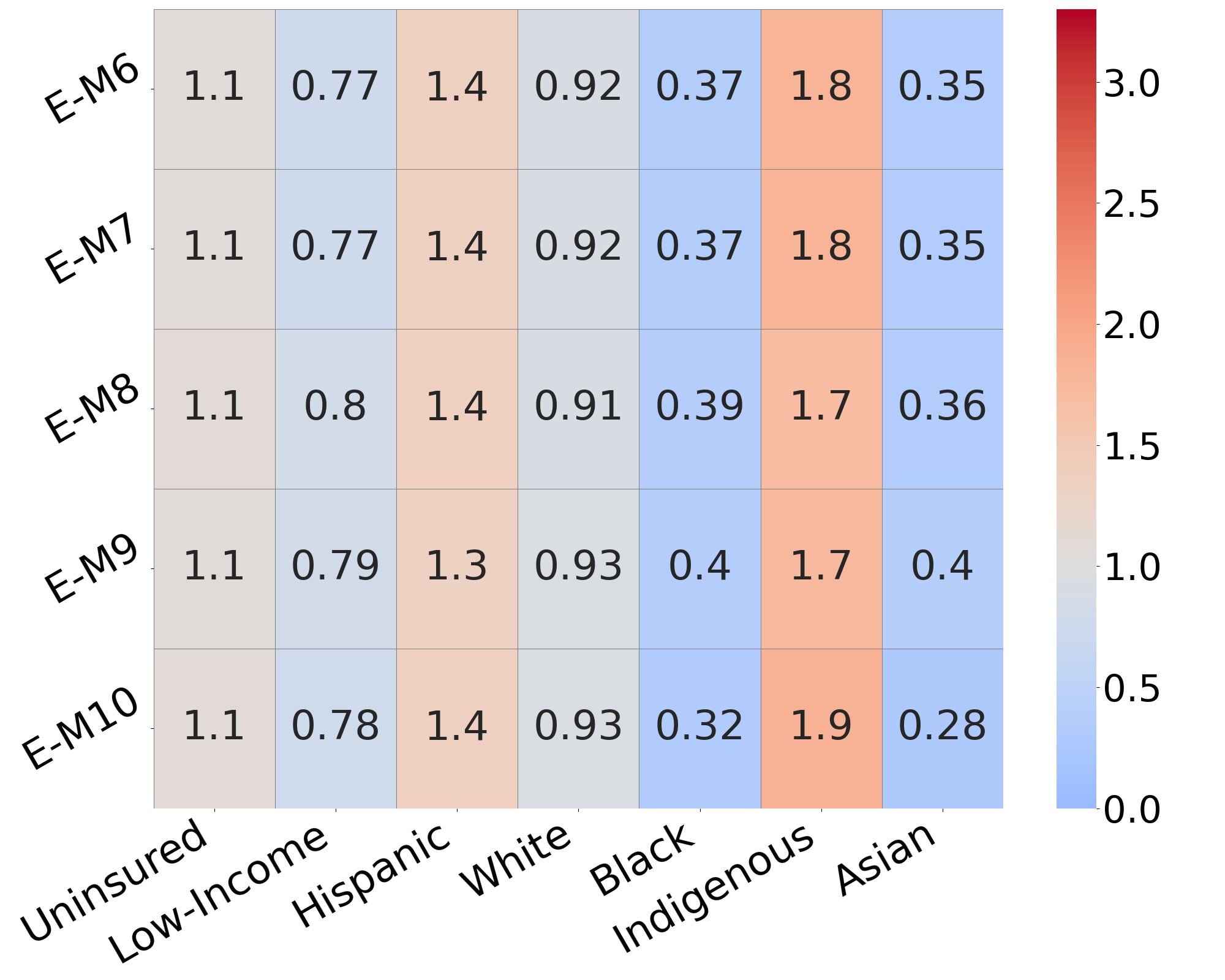}
        \caption{Relative unfairness in simulations with the {\sc Equity } Objective, \$0 budget}
        \label{fig: subfig: unfairness, B0, MMF}
    \end{subfigure}
    \begin{subfigure}[t]{.325\textwidth}
        \centering
        \includegraphics[width=.99\linewidth]{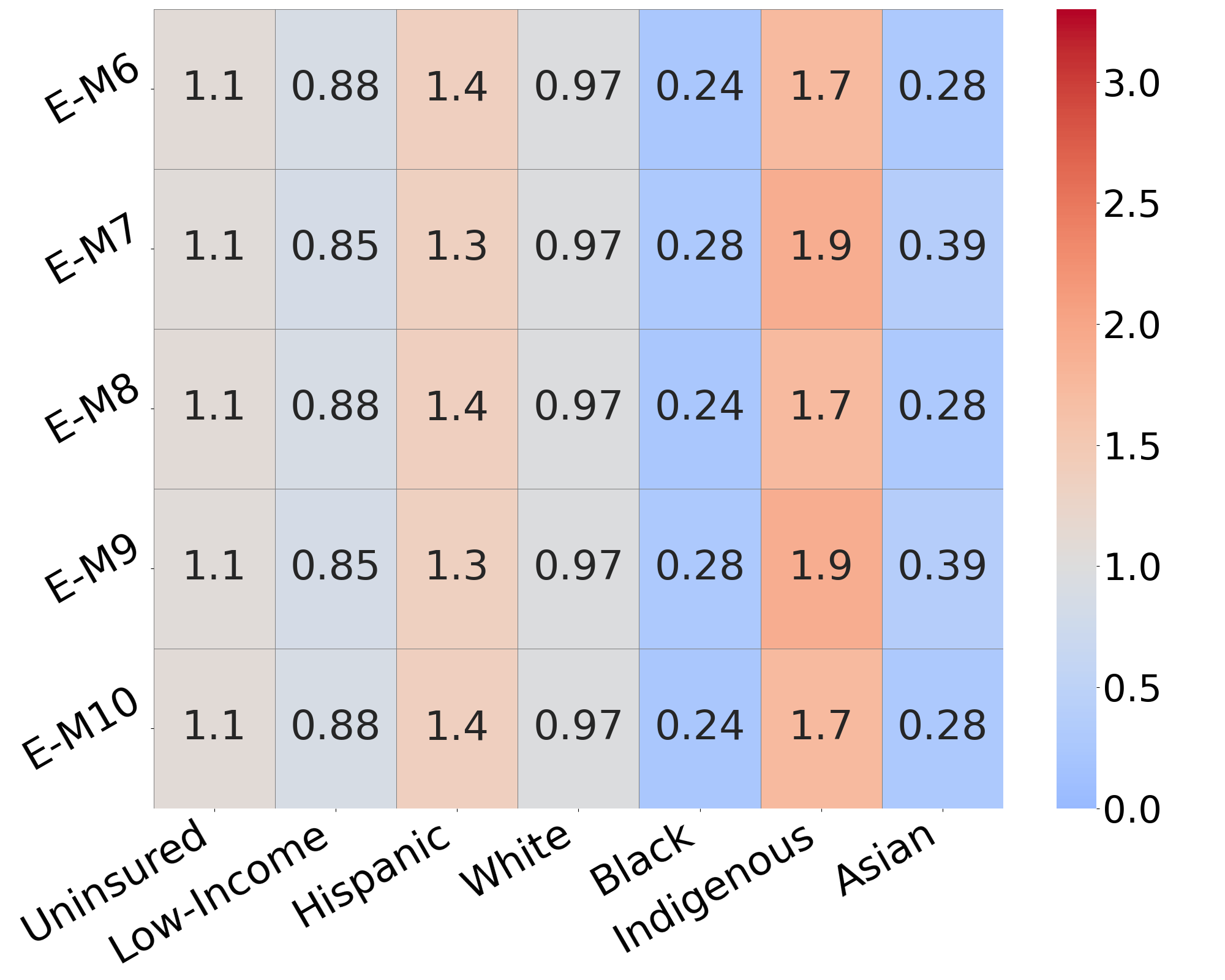}
        \caption{Relative unfairness in simulations with the {\sc Equity } Objective, \$250 budget}
        \label{fig: subfig: unfairness, B250, MMF}
    \end{subfigure}
    \begin{subfigure}[t]{.325\textwidth}
        \centering
        \includegraphics[width=.99\linewidth]{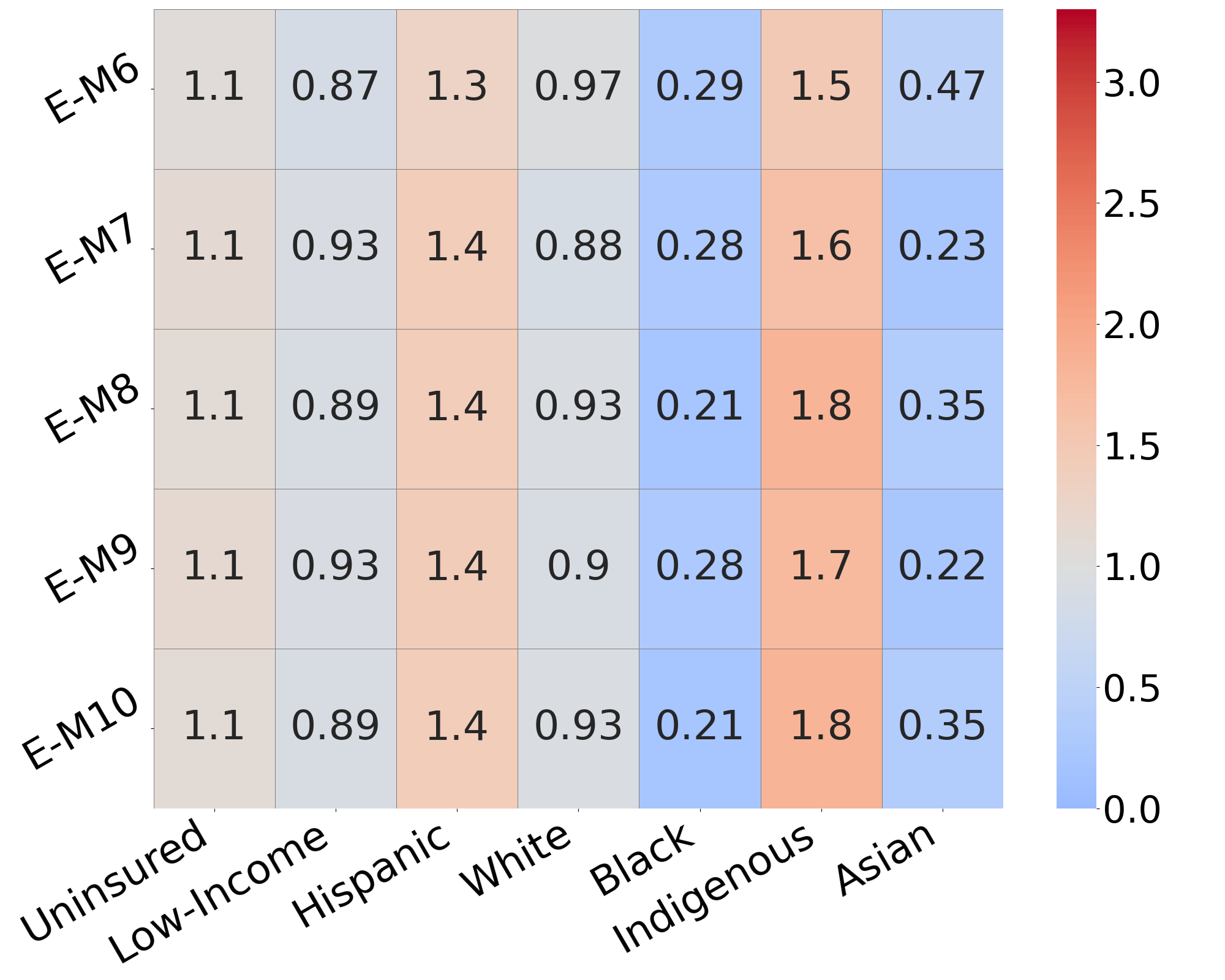}
        \caption{Relative unfairness in simulations with the {\sc Equity } Objective, \$500 budget}
        \label{fig: subfig: unfairness, B500, MMF}
    \end{subfigure}
    \caption{Relative unfairness in the percent of load demanded that is shed across budgets of \$0, \$250 million, and \$500 million. The first row shows {\sc Policy} constraints alone, and the second shows {\sc Policy} constraints and the {\sc Equity} objective.}
    \label{fig: unfairness appendix}
\end{figure*}

In this section, we discuss the all load shed results for each budget we considered, up to \$1 billion allocated in \$250 million increments for each combination of {\sc Policy} constraint and {\sc Equity} objective. 
It is nearly impossible to ascertain how equitable the switching and power line undergrounding decisions are without some degree of normalization for different group sizes. In this subsection, we discuss the percent of load demanded that is shed by each group and give a ``relative unfairness" metric. This metric computes the ratio of the percent of load shed experienced by the group to the percent of load shed experienced by the overall population. We designate an ``unfair" outcome as one in which a group experiences more than 1.1 times the overall percent load shed for a budget of \$1 billion or more or more than 1.3 times if the budget is under \$1 billion, and we bold and color the text in those cells red to call attention to these unfair outcomes. Figure \ref{fig: pct loadshed} gives the percent of load demanded that is shed for budgets of \$0, \$250 million, and \$500 million, and Figure \ref{fig: unfairness appendix} shows the relative unfairness in the percent of load demanded that is shed for those budgets. 

In Figure \ref{fig: subfig: unfairness, B0, Baseline}, we observe that when there is no budget allocated, indigenous and Hispanic groups experience unfair load shedding outcomes across the board. When considering the baseline case (BL-M1), we see that indigenous and Hispanic groups see about 2 and 1.4 times the average overall percent of load demanded that is shed, respectively. Asian and black communities experience disproportionately low levels of load shed, likely due to residence in urban areas with lower wildfire ignition risk, and low-income communities are also not at higher risk of experiencing load shed, again likely due to city poverty. If we set a 1\% threshold as an ``acceptable" percentage of load shed, we see that the overall population load shed is above this percentage threshold, and, in particular, uninsured, Hispanic, white, and indigenous groups are above this percentage threshold (Figure \ref{fig: subfig: ls, B0, Baseline}).

Figures \ref{fig: subfig: ls, B250, Baseline}, \ref{fig: subfig: ls, B250, MMF}, \ref{fig: subfig: unfairness, B250, Baseline}, and \ref{fig: subfig: ls, B250, MMF} show how these results change as we begin allocating budget for power line undergrounding. When we only consider policy constraints, we see that this small budget allocation decreases overall percent load shed by about 46\%. Furthermore, all groups except the indigenous group see their load shed percentage drop below the 1\% threshold. White and Hispanic groups seem to experience the most benefit from this investment with over 50\% reductions in percent of load demanded that is shed. Indigenous groups experience relatively less benefit, which is why the relative unfairness that they experience increases between the \$0 and \$250 million budgets when using only {\sc Policy} constraints. Now, we consider when {\sc Equity} objective is incorporated. Relative to the {\sc Policy}-constraint-only case and the 0-budget case, relative unfairness decreases for both indigenous and Hispanic groups when we introduce the {\sc Equity} objective. However, while relative unfairness is lower than that in the {\sc Policy}-constraint only case, the percent of load shed experienced is higher. For example, the overall load shed in Model M6 ({\sc Equity} objective, no {\sc Policy} constraints), leads to only a 16\% decrease in the percent of load demanded that is shed relative to the baseline model in the 0-budget case.

When budget increases to \$500 million, \$750 million, and eventually \$1 billion, indigenous load shed outcomes continue to be relatively unfair under {\sc Policy} constraints alone; even at the \$1 billion budget, use of policy constraints alone allow for relative unfairness ratios for indigenous populations sometimes over 3 times that of the overall population. In contrast, unfairness ratios under the {\sc Equity} objective are sub-2 as soon as a \$250 million budget is allocated. More importantly, by an allocation of \$500 million, indigenous percentage of load shed decreases relative to the baseline case, meaning that these groups are finally receiving real benefit from these investments \textit{without pushing other groups' percent load shed above 1\%}. This is why we argue that a \gls*{mmf} framework that minimizes the maximum percent of a group's load demanded that is shed coupled with a \textit{sufficient} budget allows for the controlling of wildfire risk \textit{and} the fair reduction of load shed from emergency power shutoffs to all considered groups. These trends continue for the \$750 million and \$1 billion cases.

\section*{Acknowledgements}
M.P. and S.G. acknowledge generous support from the NSF CAREER Grant 2239824. D.K. Molzahn and R. Piansky acknowledge support from the NSF AI Institute for Advances in Optimization (AI4OPT), \#2112533. We thank Gurobi for providing the use of academic licenses.

\bibliographystyle{plainnat}

\bibliography{refs}

\begin{thebibliography}{72}
\providecommand{\natexlab}[1]{#1}
\providecommand{\url}[1]{\texttt{#1}}
\expandafter\ifx\csname urlstyle\endcsname\relax
  \providecommand{\doi}[1]{doi: #1}\else
  \providecommand{\doi}{doi: \begingroup \urlstyle{rm}\Url}\fi

\bibitem[Ard et~al.(2016)Ard, Colen, Becerra, and Velez]{ard2016two}
Kerry Ard, Cynthia Colen, Marisol Becerra, and Thelma Velez.
\newblock Two mechanisms: {T}he role of social capital and industrial pollution exposure in explaining racial disparities in self-rated health.
\newblock \emph{International journal of environmental research and public health}, 13\penalty0 (10):\penalty0 1025, 2016.

\bibitem[Beam(2023)]{Beam_2023}
Adam Beam.
\newblock P{G}\&{E}’s plan to bury power lines and prevent wildfires faces opposition because of high rates, Oct 2023.

\bibitem[Bertsimas et~al.(2011)Bertsimas, Farias, and Trichakis]{bertsimas2011price}
Dimitris Bertsimas, Vivek~F Farias, and Nikolaos Trichakis.
\newblock The price of fairness.
\newblock \emph{Operations {R}esearch}, 59\penalty0 (1):\penalty0 17--31, 2011.

\bibitem[Bezanson et~al.(2017)Bezanson, Edelman, Karpinski, and Shah]{bezanson2017julia}
Jeff Bezanson, Alan Edelman, Stefan Karpinski, and Viral~B Shah.
\newblock Julia: {A} fresh approach to numerical computing.
\newblock \emph{SIAM Review}, 59\penalty0 (1):\penalty0 65--98, 2017.

\bibitem[Blunt(2023)]{Blunt_2023}
Katherine Blunt.
\newblock {P}{G}\&{E} scraps tree-trimming program once seen as key to fire prevention.
\newblock \emph{The Wall Street Journal}, Aug 2023.

\bibitem[Bonald and Roberts(2015)]{bonald2015multi}
Thomas Bonald and James Roberts.
\newblock Multi-resource fairness: {O}bjectives, algorithms and performance.
\newblock In \emph{Proceedings of the 2015 ACM SIGMETRICS International Conference on Measurement and Modeling of Computer Systems}, pages 31--42, 2015.

\bibitem[{California Public Utility Commision}(2019)]{cpuc_undergrounding}
{California Public Utility Commision}.
\newblock {CPUC} undergrounding programs description, 2019.

\bibitem[Campbell et~al.(2016)Campbell, Greenberg, Mankikar, and Ross]{campbell2016case}
Carla Campbell, Rachael Greenberg, Deepa Mankikar, and Ronald~D Ross.
\newblock A case study of environmental injustice: {T}he failure in {F}lint.
\newblock \emph{International journal of environmental research and public health}, 13\penalty0 (10):\penalty0 951, 2016.

\bibitem[{CDC/ATSDR Social Vulnerability Index}(2010)]{CDC_SVI_data}
{CDC/ATSDR Social Vulnerability Index}.
\newblock \emph{{Database Texas}}.
\newblock Centers for Disease Control and Prevention/ Agency for Toxic Substances and Disease Registry/ Geospatial Research, Analysis, and Services Program, 2010.

\bibitem[Climate and Economic Justice Screening Tool()]{ScreeningTool_2022}
Climate and Economic Justice Screening Tool.
\newblock {Climate and Economic Justice Screening Tool}, Nov 2022.

\bibitem[Coffrin et~al.(2018)Coffrin, Bent, Sundar, Ng, and Lubin]{coffrin2018}
Carleton Coffrin, Russell Bent, Kaarthik Sundar, Yeesian Ng, and Miles Lubin.
\newblock {PowerModels.jl}: {A}n open-source framework for exploring power flow formulations.
\newblock In \emph{19th Power Systems Computation Conference (PSCC)}, 2018.

\bibitem[Cutter(2024)]{cutter2024origin}
Susan~L Cutter.
\newblock The origin and diffusion of the {Social Vulnerability Index} ({SoVI}).
\newblock \emph{International Journal of Disaster Risk Reduction}, page 104576, 2024.

\bibitem[Cutter et~al.(2003)Cutter, Boruff, and Shirley]{cutter2003social}
Susan~L. Cutter, Bryan~J. Boruff, and W.~Lynn Shirley.
\newblock Social vulnerability to environmental hazards.
\newblock \emph{Social Science Quarterly}, 84\penalty0 (2):\penalty0 242--261, 2003.

\bibitem[Dunning et~al.(2017)Dunning, Huchette, and Lubin]{DunningHuchetteLubin2017}
Iain Dunning, Joey Huchette, and Miles Lubin.
\newblock {JuMP}: {A} modeling language for mathematical optimization.
\newblock \emph{SIAM Review}, 59\penalty0 (2):\penalty0 295--320, 2017.

\bibitem[{Federal Energy Regulatory Commision}(2020)]{ferc_ceii}
{Federal Energy Regulatory Commision}.
\newblock {Critical Energy/Electric Infrastructure Information (CEII) Regulations}, June 2020.

\bibitem[Fekete(2019)]{fekete2019social}
Alexander Fekete.
\newblock Social vulnerability (re-) assessment in context to natural hazards: {R}eview of the usefulness of the spatial indicator approach and investigations of validation demands.
\newblock \emph{International Journal of Disaster Risk Science}, 10:\penalty0 220--232, 2019.

\bibitem[Firststreet.org(2023)]{Firststreet.org_2023}
Firststreet.org.
\newblock Wildfire model methodology, Jul 2023.

\bibitem[Flanagan et~al.(2011)Flanagan, Gregory, Hallisey, Heitgerd, and Lewis]{flanagan2011social}
Barry~E Flanagan, Edward~W Gregory, Elaine~J Hallisey, Janet~L Heitgerd, and Brian Lewis.
\newblock A social vulnerability index for disaster management.
\newblock \emph{Journal of homeland security and emergency management}, 8\penalty0 (1), 2011.

\bibitem[Fuller(2019)]{Fuller_2019}
Thomas Fuller.
\newblock For the most vulnerable, {C}alifornia blackouts ‘can be life or death’.
\newblock \emph{The New York Times}, Oct 2019.

\bibitem[Ganatra et~al.(2022)Ganatra, Dani, Kumar, Khan, Wadhera, Neilan, Thavendiranathan, Barac, Hermann, Leja, et~al.]{ganatra2022impact}
Sarju Ganatra, Sourbha~S Dani, Ashish Kumar, Safi~U Khan, Rishi Wadhera, Tomas~G Neilan, Paaladinesh Thavendiranathan, Ana Barac, Joerg Hermann, Monika Leja, et~al.
\newblock Impact of social vulnerability on comorbid cancer and cardiovascular disease mortality in the {U}nited {S}tates.
\newblock \emph{Cardio Oncology}, 4\penalty0 (3):\penalty0 326--337, 2022.

\bibitem[Ganz et~al.(2023)Ganz, Duan, and Ji]{ganz2023socioeconomic}
Scott~C Ganz, Chenghao Duan, and Chuanyi Ji.
\newblock Socioeconomic vulnerability and differential impact of severe weather-induced power outages.
\newblock \emph{PNAS {N}exus}, 2\penalty0 (10):\penalty0 pgad295, 2023.

\bibitem[Gronlund(2014)]{gronlund2014racial}
Carina~J Gronlund.
\newblock Racial and socioeconomic disparities in heat-related health effects and their mechanisms: {A} review.
\newblock \emph{Current epidemiology reports}, 1:\penalty0 165--173, 2014.

\bibitem[Gupta et~al.(2023)Gupta, Moondra, and Singh]{gupta2023lp}
Swati Gupta, Jai Moondra, and Mohit Singh.
\newblock Which lp norm is the fairest? approximations for fair facility location across all" p".
\newblock In \emph{Proceedings of the 24th ACM Conference on Economics and Computation}, pages 817--817, 2023.

\bibitem[{Gurobi Optimization, LLC}(2023)]{gurobi}
{Gurobi Optimization, LLC}.
\newblock {Gurobi Optimizer Reference Manual}, 2023.

\bibitem[Hall(2012)]{hall2012out}
Kenneth~L Hall.
\newblock Out of sight, out of mind: {A}n updated study on the undergrounding of overhead power lines.
\newblock \emph{Edison Electric Institute, Washington, DC}, 2012.

\bibitem[Ham and Lee(2022)]{Ham_Lee_2022}
Youngjib Ham and Seulbi Lee.
\newblock Behavior analysis of socially vulnerable households responding to planned power shutoffs.
\newblock \emph{Natural Hazards Center Mitigation Matters Grant Report}, 2022.

\bibitem[Hinojos et~al.(2023)Hinojos, McPhillips, Stempel, and Grady]{hinojos2023social}
Selena Hinojos, Lauren McPhillips, Peter Stempel, and Caitlin Grady.
\newblock Social and environmental vulnerability to flooding: {I}nvestigating cross-scale hypotheses.
\newblock \emph{Applied geography}, 157:\penalty0 103017, 2023.

\bibitem[Hsu et~al.(2023)Hsu, Taneja, Carvallo, and Shah]{Hsu_Taneja_Carvallo_Shah_2023}
Feng~Chi Hsu, Jay Taneja, J~P Carvallo, and Zeal Shah.
\newblock Frozen out in {Texas}: {B}lackouts and inequity, Oct 2023.

\bibitem[Huang et~al.(2023)Huang, Hu, Sang, Lucas, Wong, Wang, Hong, Yao, and Donde]{huang2023review}
Can Huang, Qinran Hu, Linwei Sang, Donald~D Lucas, Robin Wong, Bin Wang, Wanshi Hong, Mengqi Yao, and Vaibhav Donde.
\newblock A review of {Public Safety Power Shutoffs} {(PSPS)} for wildfire mitigation: Policies, practices, models and data sources.
\newblock \emph{IEEE Transactions on Energy Markets, Policy and Regulation}, 1\penalty0 (3):\penalty0 187--197, 2023.

\bibitem[Javanmard et~al.(2023)Javanmard, Lee, Kim, Liu, and Diab]{javanmard2023impacts}
Reyhane Javanmard, Jinhyung Lee, Junghwan Kim, Luyu Liu, and Ehab Diab.
\newblock The impacts of the modifiable areal unit problem {(MAUP)} on social equity analysis of public transit reliability.
\newblock \emph{Journal of Transport Geography}, 106:\penalty0 103500, 2023.

\bibitem[Jesdale et~al.(2013)Jesdale, Morello-Frosch, and Cushing]{jesdale2013racial}
Bill~M Jesdale, Rachel Morello-Frosch, and Lara Cushing.
\newblock The racial/ethnic distribution of heat risk--related land cover in relation to residential segregation.
\newblock \emph{Environmental health perspectives}, 121\penalty0 (7):\penalty0 811--817, 2013.

\bibitem[Kalai and Smorodinsky(1975)]{kalai1975other}
Ehud Kalai and Meir Smorodinsky.
\newblock Other solutions to {N}ash's bargaining problem.
\newblock \emph{Econometrica: Journal of the Econometric Society}, pages 513--518, 1975.

\bibitem[Kody et~al.(2022{\natexlab{a}})Kody, Piansky, Molzahn, and K]{kody2022optimizing}
Alyssa Kody, Ryan Piansky, Molzahn, and Daniel K.
\newblock Optimizing transmission infrastructure investments to support line de-energization for mitigating wildfire ignition risk.
\newblock \emph{IREP Symposium on Bulk Power System Dynamics and Control}, 2022{\natexlab{a}}.

\bibitem[Kody et~al.(2022{\natexlab{b}})Kody, West, and Molzahn]{kody2022load}
Alyssa Kody, Amanda West, and Daniel~K Molzahn.
\newblock Sharing the load: {C}onsidering fairness in de-energization scheduling to mitigate wildfire ignition risk using rolling optimization.
\newblock In \emph{61st IEEE Conference on Decision and Control (CDC)}, pages 5705--5712. IEEE, 2022{\natexlab{b}}.

\bibitem[MacMillan and Englund(2021)]{macmillan2021longer}
Douglas MacMillan and Will Englund.
\newblock Longer, more frequent outages afflict the {US} power grid as states fail to prepare for climate change.
\newblock \emph{The Washington Post}, 2021.

\bibitem[Mitchell(2013)]{mitchell2013power}
Joseph~W Mitchell.
\newblock Power line failures and catastrophic wildfires under extreme weather conditions.
\newblock \emph{Engineering Failure Analysis}, 35:\penalty0 726--735, 2013.

\bibitem[Mueller et~al.(2009)Mueller, Loomis, and Gonz{\'a}lez-Cab{\'a}n]{mueller2009repeated}
Julie Mueller, John Loomis, and Armando Gonz{\'a}lez-Cab{\'a}n.
\newblock Do repeated wildfires change homebuyers’ demand for homes in high-risk areas? {A} hedonic analysis of the short and long-term effects of repeated wildfires on house prices in {Southern California}.
\newblock \emph{The Journal of Real Estate Finance and Economics}, 38:\penalty0 155--172, 2009.

\bibitem[PACE(2017)]{PACE}
PACE.
\newblock \emph{{P}artnership for an {A}dvanced {C}omputing {E}nvironment ({PACE})}, 2017.

\bibitem[Penn(2024)]{Penn_2024}
Ivan Penn.
\newblock Utility-caused wildfires are becoming a national problem.
\newblock \emph{The New York Times}, Mar 2024.

\bibitem[PG\&E(2023)]{PGE}
PG\&E, May 2023.

\bibitem[Piansky et~al.(2025)Piansky, Taylor, Rhodes, Molzahn, Roald, and Watson]{piansky2025hicss}
R.~Piansky, S.~Taylor, N.~Rhodes, D.~K. Molzahn, L.~A. Roald, and Jean-Paul Watson.
\newblock {Quantifying Metrics for Wildfire Ignition Risk from Geographic Data in Power Shutoff Decision-Making}.
\newblock In \emph{58th Hawaii International Conference on System Sciences (HICSS)}, January 2025.

\bibitem[Piansky et~al.(2024)Piansky, Stinchfield, Kody, Molzahn, and Watson]{piansky2024long}
Ryan Piansky, Georgia Stinchfield, Alyssa Kody, Daniel~K. Molzahn, and Jean-Paul Watson.
\newblock Long duration battery sizing, siting, and operation under wildfire risk using progressive hedging.
\newblock 2024.

\bibitem[Ramesh et~al.(2022)Ramesh, Jagger, Zaitchik, Kolivras, Swarup, Deanes, Hallisey, Sharpe, and Gohlke]{ramesh2022flooding}
Balaji Ramesh, Meredith~A Jagger, Benjamin Zaitchik, Korine~N Kolivras, Samarth Swarup, Lauren Deanes, Elaine Hallisey, J~Danielle Sharpe, and Julia~M Gohlke.
\newblock Flooding and emergency department visits: {E}ffect modification by the {CDC/ATSDR Social Vulnerability Index}.
\newblock \emph{International Journal of Disaster Risk Reduction}, 76:\penalty0 102986, 2022.

\bibitem[Rhodes and Roald(2022)]{rhodes2022co}
Noah Rhodes and Line Roald.
\newblock Co-optimization of power line shutoff and restoration for electric grids under high wildfire ignition risk.
\newblock \emph{arXiv e-prints}, pages arXiv--2204, 2022.

\bibitem[Rhodes et~al.(2020)Rhodes, Ntaimo, and Roald]{rhodes2020balancing}
Noah Rhodes, Lewis Ntaimo, and Line Roald.
\newblock Balancing wildfire risk and power outages through optimized power shut-offs.
\newblock \emph{IEEE Transactions on Power Systems}, 36\penalty0 (4):\penalty0 3118--3128, 2020.

\bibitem[Rhodes et~al.(2024)Rhodes, Coffrin, and Roald]{rhodes2024security}
Noah Rhodes, Carleton Coffrin, and Line Roald.
\newblock Security constrained optimal power shutoff for wildfire risk mitigation.
\newblock \emph{IET Generation, Transmission \& Distribution}, 2024.

\bibitem[Rose and Balarajan(2024)]{rose2024utility}
Andy Rose and Brammhi Balarajan.
\newblock Utility company says its facilities ‘appear to have been involved’ in start of {Smokehouse Creek} fire in {Texas}, Mar 2024.

\bibitem[Rufat et~al.(2019)Rufat, Tate, Emrich, and Antolini]{rufat2019valid}
Samuel Rufat, Eric Tate, Christopher~T Emrich, and Federico Antolini.
\newblock How valid are social vulnerability models?
\newblock \emph{Annals of the American Association of Geographers}, 109\penalty0 (4):\penalty0 1131--1153, 2019.

\bibitem[Schaider et~al.(2019)Schaider, Swetschinski, Campbell, and Rudel]{schaider2019environmental}
Laurel~A Schaider, Lucien Swetschinski, Christopher Campbell, and Ruthann~A Rudel.
\newblock Environmental justice and drinking water quality: {A}re there socioeconomic disparities in nitrate levels in {U}{S} drinking water?
\newblock \emph{Environmental Health}, 18:\penalty0 1--15, 2019.

\bibitem[Schweizer(2019)]{USDA_wildfire_season}
Deb Schweizer.
\newblock {Wildfires in All Seasons?}, June 2019.
\newblock Accessed: 2022-02-18.

\bibitem[Shah et~al.(2023)Shah, Carvallo, Hsu, and Taneja]{shah2023inequitable}
Zeal Shah, Juan~Pablo Carvallo, Feng-Chi Hsu, and Jay Taneja.
\newblock The inequitable distribution of power interruptions during the 2021 {T}exas winter storm {U}ri.
\newblock \emph{Environmental Research: Infrastructure and Sustainability}, 3\penalty0 (2):\penalty0 025011, 2023.

\bibitem[Spielman et~al.(2020)Spielman, Tuccillo, Folch, Schweikert, Davies, Wood, and Tate]{spielman2020evaluating}
Seth~E Spielman, Joseph Tuccillo, David~C Folch, Amy Schweikert, Rebecca Davies, Nathan Wood, and Eric Tate.
\newblock Evaluating social vulnerability indicators: {C}riteria and their application to the {Social Vulnerability Index}.
\newblock \emph{Natural hazards}, 100:\penalty0 417--436, 2020.

\bibitem[Strauss and Smith(2009)]{strauss2009construct}
Milton~E Strauss and Gregory~T Smith.
\newblock Construct validity: {A}dvances in theory and methodology.
\newblock \emph{Annual review of clinical psychology}, 5\penalty0 (1):\penalty0 1--25, 2009.

\bibitem[Su et~al.(2024)Su, Mehrani, Dehghanian, and Lejeune]{su2024quasi}
Jinshun Su, Saharnaz Mehrani, Payman Dehghanian, and Miguel~A. Lejeune.
\newblock Quasi second-order stochastic dominance model for balancing wildfire risks and power outages due to proactive public safety de-energizations.
\newblock \emph{IEEE Transactions on Power Systems}, 39\penalty0 (2):\penalty0 2528--2542, 2024.

\bibitem[Taylor et~al.(2023)Taylor, Setyawan, Cui, Zamzam, and Roald]{Taylor2023Managing}
Sofia Taylor, Gabriela Setyawan, Bai Cui, Ahmed Zamzam, and Line~A. Roald.
\newblock Managing wildfire risk and promoting equity through optimal configuration of networked microgrids.
\newblock In \emph{14th ACM International Conference on Future Energy Systems}, e-Energy '23, page 189–199. Association for Computing Machinery, 2023.
\newblock ISBN 9798400700323.

\bibitem[{Texas A\&M}(2014)]{powerfiretexas}
{Texas A\&M}.
\newblock {How do power lines cause wildfires?}, July 2014.

\bibitem[{US Census 2010}(2013)]{COP_data}
{US Census 2010}.
\newblock \emph{2010 Center of Population by Census Tract}.
\newblock US Census Bureau, Mar 2013.

\bibitem[{US Census 2019}(2019)]{demographic_data}
{US Census 2019}.
\newblock \emph{2019 Tract Data}.
\newblock US Census Bureau, 2019.

\bibitem[{US Census Bureau}(2022{\natexlab{a}})]{ACS_2022_Hispanic}
{US Census Bureau}.
\newblock \emph{2022: {A}{C}{S} 1-Year Supplemental Estimates Table {K}200301}.
\newblock {US Census Bureau}, 2022{\natexlab{a}}.

\bibitem[{US Census Bureau}(2022{\natexlab{b}})]{ACS_2022_poverty}
{US Census Bureau}.
\newblock 2022: {A}{C}{S} 1-year supplemental estimates table {S}1703, 2022{\natexlab{b}}.

\bibitem[{US Census Bureau}(2022{\natexlab{c}})]{ACS_2022_race}
{US Census Bureau}.
\newblock 2022: {A}{C}{S} 1-year supplemental estimates table {K}200201, 2022{\natexlab{c}}.

\bibitem[US Energy Information Administration()]{saidi2022table}
US Energy Information Administration.
\newblock {SAIDI} values (minutes per year) of {U.S.} distribution system by state, 2023.

\bibitem[{US Geological Survey}(2021)]{USGS_WFPI}
{US Geological Survey}.
\newblock {Wildland Fire Potential Index}.
\newblock https://www.usgs.gov/fire-danger-forecast/wildland-fire-potential-index-wfpi, 2021.
\newblock Accessed: 2022-02-18.

\bibitem[Walls et~al.(2024)Walls, Hines, and Ruggles]{walls2024implementation}
Margaret Walls, Sofia Hines, and Logan Ruggles.
\newblock Implementation of {Justice40}: {C}hallenges, opportunities, and a status update.
\newblock \emph{RFF Report}, 2024.

\bibitem[Waseem et~al.(2021)Waseem, Bayani, Manshadi, and Tavakol-Davani]{waseem2021}
Muhammad Waseem, Reza Bayani, Saeed~D. Manshadi, and Hassan Tavakol-Davani.
\newblock Quantifying the risk of wildfire ignition by power lines under extreme weather conditions.
\newblock \emph{arXiv:2110.05551}, October 2021.

\bibitem[{White House Office}(2021)]{exec_order2021}
{White House Office}, editor.
\newblock \emph{{U.S. Executive Order 14008}}.
\newblock Number 14008 in 14008. {White House Office}, 2021.

\bibitem[{White House Office}(2023)]{White_House_2023}
{White House Office}, editor.
\newblock \emph{Justice40: {A} whole-of-government initiative}.
\newblock {White House Office}, Apr 2023.

\bibitem[Wolkin et~al.(2022)Wolkin, Collier, House, Reif, Motsinger-Reif, Duca, and Sharpe]{wolkin2022comparison}
Amy Wolkin, Sarah Collier, John~S House, David Reif, Alison Motsinger-Reif, Lindsey Duca, and J~Danielle Sharpe.
\newblock Comparison of national vulnerability indices used by the {Centers for Disease Control and Prevention} for the {COVID-19} response.
\newblock \emph{Public Health Reports}, 137\penalty0 (4):\penalty0 803--812, 2022.

\bibitem[Woo et~al.(2019)Woo, Kravitz-Wirtz, Sass, Crowder, Teixeira, and Takeuchi]{woo2019residential}
Bongki Woo, Nicole Kravitz-Wirtz, Victoria Sass, Kyle Crowder, Samantha Teixeira, and David~T Takeuchi.
\newblock Residential segregation and racial/ethnic disparities in ambient air pollution.
\newblock \emph{Race and social problems}, 11:\penalty0 60--67, 2019.

\bibitem[Xu et~al.(2017)Xu, Birchfield, Shetye, and Overbye]{xu2017creation}
Ti~Xu, Adam~B Birchfield, Komal~S Shetye, and Thomas~J Overbye.
\newblock Creation of synthetic electric grid models for transient stability studies.
\newblock In \emph{10th Bulk Power Systems Dynamics and Control Symposium (IREP 2017)}, pages 1--6, 2017.

\bibitem[Xu et~al.(2018)Xu, Birchfield, and Overbye]{xu2018synthetic}
Ti~Xu, Adam~B. Birchfield, and Thomas~J. Overbye.
\newblock Modeling, tuning, and validating system dynamics in synthetic electric grids.
\newblock \emph{IEEE Transactions on Power Systems}, 33\penalty0 (6):\penalty0 6501--6509, 2018.

\bibitem[Yang et~al.(2022)Yang, Sparrow, Ashtine, Wallom, and Morstyn]{yang2022resilient}
Weijia Yang, Sarah~N Sparrow, Masa{\=o} Ashtine, David~CH Wallom, and Thomas Morstyn.
\newblock Resilient by design: {P}reventing wildfires and blackouts with microgrids.
\newblock \emph{Applied Energy}, 313:\penalty0 118793, 2022.

\end{thebibliography}

\end{document}